\documentclass[11pt, a4paper]{amsart}
\usepackage{amscd,amssymb,eucal,eufrak,mathrsfs,amsmath}
\usepackage{hyperref} %,backref

\input epsf.tex

\input xy
\xyoption{all}

\newtheorem{thm}{Theorem}[section]
\newtheorem{cor}[thm]{Corollary}
\newtheorem{lem}[thm]{Lemma}
\newtheorem{prop}[thm]{Proposition}

\theoremstyle{definition}
\newtheorem{defn}[thm]{Definition}

\theoremstyle{remark}
\newtheorem{remark}[thm]{Remark}
\newtheorem{remarks}[thm]{Remarks}

\newtheorem{examples}[thm]{Examples}
\numberwithin{equation}{section}
\newtheorem{problem}[thm]{Problem}
\newcommand{\delete}[1]{} % Comment out text.

\newcommand{\ben}{\begin{enumerate}}
\newcommand{\een}{\end{enumerate}}

\newcommand{\sk}{\vskip 0.3cm}

\def\Homeo{{\mathrm{Homeo}}\,}

\def\QED{\nobreak\quad\ifmmode\roman{Q.E.D.}\else{\rm Q.E.D.}\fi}

\newcommand{\Acal}{\mathcal{A}}

\newcommand{\Lcal}{\mathcal{L}}

\newcommand{\Pcal}{\mathcal{P}}

\newcommand{\Bb}{\mathbf{B}}

\newcommand{\Q}{\mathbb{Q}}
\newcommand{\R}{\mathbb{R}}
\newcommand{\C}{\mathbb{C}}
\newcommand{\N}{\mathbb{N}}
\newcommand{\T}{\mathbb{T}}

\newcommand{\Z}{\mathbb{Z}}

\newcommand{\OC}{\bar{\mathcal{O}}}
\newcommand{\OCT}{\bar{\mathcal{O}}_T}

\def\a{\alpha}

\def\eps{{\varepsilon}}

\newcommand{\al}{\alpha}

\newcommand{\ga}{\gamma}
\newcommand{\del}{\delta}

\newcommand{\ep}{\varepsilon}
\newcommand{\sig}{\sigma}

\newcommand{\la}{\lambda}
\newcommand{\La}{\Lambda}
\newcommand{\tet}{\theta}
\newcommand{\Tet}{\Theta}
\newcommand{\om}{\omega}
\newcommand{\Om}{\Omega}

\newcommand{\OG}{\mathcal{O}_{G}}
\newcommand{\OCG}{ {\ol{\mathcal{O}}}_{G} }

\newcommand{\br}{\vspace{2 mm}}

\newcommand{\rest}{\upharpoonright}

\newcommand{\ol}{\overline}

\newcommand{\cls}{\rm{cls}}

\newcommand{\Iso}{\rm{Iso\,}}

\newcommand{\id}{\rm{id}}

\newcommand{\diam}{\rm{diam\,}}

\newcommand{\card}{\rm{card\,}}

\newcommand{\eva}{\rm{eva}}

\newcommand{\Prol}{\rm{Prol\,}}

\newcommand{\eliG}{\rm{{\ell}^{\infty}}(G)}

\newtheorem{rmk}[thm]{Remark}
\newtheorem{exa}[thm]{Example}

\begin{document}

\title
[Hereditarily non-sensitive dynamical systems] {Hereditarily
non-sensitive dynamical systems and linear representations}

%Authors
\author[E. Glasner]{E. Glasner}
\address{Department of Mathematics,
Tel-Aviv University, Ramat Aviv, Israel}
\email{glasner@math.tau.ac.il}
\urladdr{http://www.math.tau.ac.il/$^\sim$glasner}

\author[M. Megrelishvili]{M. Megrelishvili}
\address{Department of Mathematics,
Bar-Ilan University, 52900 Ramat-Gan, Israel}
\email{megereli@math.biu.ac.il}
\urladdr{http://www.math.biu.ac.il/$^\sim$megereli}

\date{November, 2005}


\begin{abstract}
For an arbitrary topological group $G$ any compact $G$-dynamical
system $(G,X)$ can be linearly $G$-represented as a
weak$^*$-compact subset of a dual Banach space $V^*$. As was shown
in \cite{M1} the Banach space $V$ can be chosen to
 be reflexive iff the metric system $(G,X)$ is weakly
almost periodic (WAP). In this paper we study the wider class of
compact $G$-systems which can be linearly represented as a
weak$^*$-compact subset of a dual Banach space with the
Radon-Nikod\'ym property. We call such a system a Radon-Nikod\'ym
system (RN). One of our main results is to show that for
metrizable compact $G$-systems the three classes: RN, HNS
(hereditarily not sensitive) and HAE (hereditarily almost
equicontinuous) coincide. We investigate these classes and their
relation to previously studied classes of $G$-systems such as WAP
and LE (locally equicontinuous). We show that the Glasner-Weiss
examples of recurrent-transitive locally equicontinuous but not
weakly almost periodic cascades are actually RN. Using
fragmentability and Namioka's theorem we give an enveloping
semigroup characterization of HNS systems and show that the enveloping
semigroup $E(X)$ of a compact metrizable HNS $G$-system is a
separable Rosenthal compact, hence of cardinality $\le
2^{\aleph_0}$. We investigate a dynamical version of the
Bourgain-Fremlin-Talagrand dichotomy and a
dynamical version of Todor\u{c}evi\'{c} dichotomy concerning
Rosenthal compacts.
\end{abstract}

\thanks{The second author thanks the Israel Science Foundation (grant
number 4699).}


\maketitle

%%%%%%%%%%%%%%%%%%%%%%%%%%%%%%%%%%%%%%%%%%%%%%%%%%%%%%%%%%%%%%%%%%%%%
%%%%%%%%%%%%%%%%%%%%%%%%  CONTENTS   %%%%%%%%%%%%%%%%%%%%%%%%%%%%%%%%
%%%%%%%%%%%%%%%%%%%%%%%%%%%%%%%%%%%%%%%%%%%%%%%%%%%%%%%%%%%%%%%%%%%%%

\tableofcontents
%\addtocontents{toc}{subsection}{\protect\hspace{10pt}}

%%%%%%%%%%%%%%%%%%%%%%%%%%%%%%%%%%%%%%%%%%%%%%%%%%%%%%%%%%%%%%%%%%%%%
%%%%%%%%%%%%%%%%%%%%%%%%%%%%  TEXT   %%%%%%%%%%%%%%%%%%%%%%%%%%%%%%%%
%%%%%%%%%%%%%%%%%%%%%%%%%%%%%%%%%%%%%%%%%%%%%%%%%%%%%%%%%%%%%%%%%%%%%

\setcounter{tocdepth}{1}

\section*{Introduction}

\sk

The main goal of this paper is to exhibit new and perhaps
unexpected connections between the (lack of) chaotic behavior of a
dynamical system and the existence of linear representations of
the system on certain Banach spaces. The property {\em sensitive
dependence on initial conditions\/} appears as a basic constituent
in several definitions of ``chaos" (see, for example,
\cite{AuYo,De,GW1,BGKM} and references thereof). In the present
paper we introduce the classes of {\em hereditarily not
sensitive\/} (HNS for short; intuitively these are the non-chaotic
systems) and {\em hereditarily almost equicontinuous systems\/}
(HAE). It turns out that these classes of dynamical systems are
well behaved with respect to the standard operations on dynamical
systems and they admit elegant characterizations in terms of
Banach space representations.

For an arbitrary topological group $G$ any compact $G$-system $X$
can be linearly $G$-represented as a weak$^*$-compact subset of a
dual Banach space $V^*$. As was shown in \cite{M1} the Banach
space $V$ can be chosen to be reflexive iff the metric $G$-system
$X$ is weakly almost periodic (WAP). We say that a dynamical
system $(G,X)$ is a Radon-Nikod\'ym system (RN) if $V^*$ can be
chosen as a Banach space with the Radon-Nikod\'ym property.
One of our main results is to show that for metrizable compact
$G$-systems the three classes of RN, HNS and HAE dynamical systems
coincide. For general compact $G$-systems $X$ we prove that $X$ is
in the class $\rm{HNS}$ iff $X$ is RN-approximable. In other
words: a compact system is non-chaotic if and only if it admits
sufficiently many $G$-representations in RN
dual Banach spaces. The link between the various topological
dynamics aspects of almost equicontinuity on the one hand and the
Banach space RN properties on the other hand is the versatile
notion of {\em fragmentability}. It played a central role in the
works on RN compacta (see e.g. Namioka \cite{N}) and their
dynamical analogues (see Megrelishvili \cite{Mefr,M,M1}). It also
serves as an important tool in the present work.

The following brief historical review will hopefully help the
reader to get a clearer perspective on the context of our results.
The theory of weakly almost periodic (WAP) functions on
topological groups was developed by W. F. Eberlein, \cite{Eb}, A.
Grothendieck, \cite{Gro} and I. Glicksberg and K. de Leeuw,
\cite{deLG}. About thirty years ago, W. A. Veech in an attempt to
unify and generalize the classical theory of weakly almost
periodic functions on a discrete group $G$, introduced a class of
functions in $\eliG$ which he denoted by $K(G)$, \cite{V}. He
showed that $K(G)$ is a uniformly closed left and right
$G$-invariant subalgebra of $\eliG$ containing the algebra of
weakly almost periodic functions $WAP(G)$ and shares with $WAP(G)$
the property that every minimal function in $K(G)$ is actually
%oct17
%deleting "Bohr" here and in many other places
%("almost periodicity" we used via "norm compactness"
%is not coincide with "Bohr almost periodicity"
almost periodic.

In \cite{S} Shtern has shown that for any compact Hausdorff
semitopological semigroup $S$ there exists a reflexive Banach
space $V$ such that $S$ is topologically isomorphic to a closed
subsemigroup of $\Bb=\{s\in \Lcal(V): \|s\|\le 1\}$. Here
$\Lcal(V)$ is the Banach space of bounded linear operators from
$V$ to itself and $\Bb$ is equipped with the weak operator
topology. Megrelishvili provided an alternative proof for this
theorem in \cite{M} and has shown in \cite{M1} that WAP dynamical
systems are characterized as those systems that have sufficiently
many linear $G$-representations on weakly compact subsets of
reflexive Banach spaces.
%aref  (1)
In particular, if $V$ is a reflexive Banach space then for every
topological subgroup $G$ of the linear isometry group $\Iso(V)$
the natural action of $G$ on the weak${}^*$ compact unit ball
$V_1^*$ of $V^*$ is WAP. Moreover, every WAP metric compact
$G$-space $X$ is a $G$-subsystem of $V_1^*$ for a suitable
reflexive Banach space $V$.
%end

A seemingly independent development is the new theory of Almost
Equicontinuous dynamical systems (AE). This was developed in a
series of papers, Glasner \& Weiss \cite{GW1}, Akin, Auslander \&
Berg \cite{AAB1,AAB2} and Glasner \& Weiss \cite{GW}. In the
latter the class of Locally Equicontinuous dynamical systems (LE)
was introduced and studied. It was shown there that the collection
$LE(G)$ of locally equicontinuous functions forms a uniformly
closed $G$-invariant subalgebra of $\ell^\infty(G)$ containing
$WAP(G)$ and having the property that each minimal function in
$LE(G)$ is almost periodic.

Of course the classical theory of WAP functions is valid for a
general topological group $G$ and it is not hard to see that the
AE theory, as well as the theory of $K(G)$-functions --- which we
call Veech functions --- extend to such groups as well.

Let $V$ be a Banach space, $V^*$ its dual. A compact dynamical
$G$-system $X$ is $V^*$-{\em representable\/} if there exist a
weakly continuous co-homomorphism $G \to{\Iso}(V)$, where
${\Iso}(V)$ is the group of linear isometries of a Banach space
$V$ onto itself, and a $G$-embedding $\phi: X \to V_1^*$, where
$V_1^*$ is the weak$^*$-compact unit ball of the dual Banach space
$V^*$ and the $G$-action is the dual action induced on $V_1^*$
from the $G$-action on $V$. An old observation (due to Teleman
\cite{Tel}) is that every compact dynamical $G$-system $X$ is
$C(X)^*$-representable.

The notion of an Eberlein compact (Eb) space in the sense of Amir
and Lindenstrauss \cite{AL} is well studied and it is known that
such spaces are characterized by being homeomorphic to a weakly
compact subset of a Banach (equivalently: reflexive Banach) space.
Later the notion of Radon-Nikod\'ym (RN) compact
%gl+
topological
spaces was introduced. These can be characterized as weak$^*$
compact sets in the duals $V^*$ with the RN property.
%gl+
A Banach space $V$ whose dual has the Radon-Nikod\'ym property is
called an {\it Asplund space} (see, for example, \cite{F, N} and
Remark \ref{r:fr}.3). We refer to the excellent 1987 paper of I.
Namioka \cite{N} where the theory of RN compacts is expounded.

One of the main objects of \cite{M1} was the investigation of RN
systems (a dynamical analog of RN compacta) and the related class
of functions called ``Asplund functions". More precisely, call a
dynamical system which is linearly representable as a
weak$^*$-compact subset of a dual Banach space with the
Radon-Nikod\'ym property
a {\em Radon-Nikod\'ym system\/} (RN for short). The class of
RN-{\em approximable\/} systems, that is the subsystems of a
product of RN systems, will be denoted by $\rm{RN_{app}}$. It was
shown in \cite{M1} that WAP $\subset \rm{RN_{app}} \subset$ LE.

Given a compact dynamical $G$-system $X$, a subgroup $H < G$ and a
function $f\in C(X)$, define a pseudometric $\rho_{H,f}$ of $X$ as
follows: $$\rho_{H,f}(x,x')=\sup_{h\in H}|f(h x) - f(h x')|.$$ We
say that $f$ is an {\em Asplund\/} function (notation: $f \in
Asp(X))$ if the pseudometric space $(X,\rho_{H,f})$ is separable
for every countable subgroup $H < G$. These are exactly the
functions which come from linear $G$-representations of $X$ on
$V^*$ with Asplund $V$. By \cite{M1}, a compact $G$-system $X$ is
$\rm{RN_{app}}$ iff $C(X)=Asp(X)$ and always $WAP(X) \subset
Asp(X)$.

\br

The first section of the paper is a brief review of some known
aspects of abstract topological dynamics which provide a
convenient framework for our results. In the Second we discuss
enveloping semigroups and semigroup compactifications. Our
treatment differs slightly from the traditional approach and
terminology and contains some new observations. For more details
refer to the books \cite{E,Gl,G2,Vr,BJM} and \cite{Au}. See also
\cite{AuHa, La, Vr-oldpaper}.

In \cite{Ko}
 K\"{o}hler shows that the well known Bourgain-Fremlin-Talagrand
dichotomy, when applied to the family $\{f^n: n\in \N\}$ of
iterates of a continuous interval map $f: I \to I$, yields a
corresponding dichotomy for the enveloping semigroups. In the
third section we generalize this and obtain a
Bourgain-Fremlin-Talagrand dichotomy
for enveloping semigroups of metric dynamical systems.

Section \ref{sect:quasi} treats the property of m-{\it
approximability}, i.e. the property of being approximable by
metric systems. For many groups $G$ every dynamical $G$-system is
m-approximable and we characterize such groups as being exactly
the {\it uniformly Lindel\"{o}f} groups.

In Section \ref{Sec-ae} we recall some important notions like
almost equicontinuity, WAP and LE and relate them to universal
systems. We also study
the related notion of {\it lightness} of a function $f \in
RUC(G)$ -- the coincidence of the pointwise and
the norm topologies on
its $G$-orbit.

Section \ref{sec:fr} is devoted to some results concerning
fragmentability. These will be crucial at many points in the rest
of the paper. In Section \ref{s:Asplund} we investigate Asplund
functions and their relations to fragmentability. In Section
\ref{s:Veech} we deal with the related class of Veech functions.
As already mentioned the latter class $K(G)$ is a generalization
of Veech's definition \cite{V}. We show that every Asplund
function is a Veech function and that for separable groups these
two classes coincide.

In Section \ref{s:main}
we introduce the dynamical properties of HAE and HNS and show that
they are intimately related to the linear representation condition
of being an RN system. In particular for metrizable compact
systems we establish the following equalities and inclusions:
$$
Eb=WAP \subset RN = HAE =HNS = RN_{app} \subset LE.
$$
Here $Eb$ stands for Eberlein systems -- a dynamical version of
Eberlein compacts (see Definition \ref{d-RN}).
Section \ref{s:ex}
is devoted to various examples and applications. We show that for
symbolic systems the RN property is equivalent to having a
countable phase space; and that any $\Z$-dynamical system $(f,
X)$, where $X$ is either the unit interval or the unit circle and
$f: X\to X$ is a homeomorphism, is an RN system.

In Section \ref{Sec-GW} we show that the Glasner-Weiss examples of
recurrent-transitive LE but not WAP metric cascades are actually
HAE. In Section \ref{s:mincenter} we investigate the mincenter of
an HAE system,
and in Section \ref{LE-HAE} we use a modified
construction to produce an example of a recurrent-transitive,
LE but not HAE system.
This example exhibits the sharp distinction
between the possible mincenters of LE and HAE systems.

%gl+

In Section \ref{envelop}, using {\it fragmented families} of functions
and {\it Namioka's joint continuity theorem}, we establish an
enveloping semigroup characterization of Asplund functions and HNS
systems. Our results imply that the Ellis semigroup $E(X)$ of a
compact metrizable HNS system $(G,X)$ is a Rosenthal compact. In
particular, by a result of Bourgain-Fremlin-Talagrand \cite{BFT},
we obtain that $E(X)$ is {\it angelic} (hence, it cannot contain a
subspace homeomorphic to $\beta \N$).
%gl+
Finally in Section \ref{Tod} we show how a theorem of
Todor\u{c}evi\'c implies that for a metric RN system $E(X)$ either
contains an uncountable discrete subspace or it admits an at most
two-to-one metric $G$-factor.

We are indebted to Stevo Todor\u{c}evi\'c for enlightening
comments. Thanks are due to Hanfeng Li for a critical reading of
the manuscript and his consequent fruitful suggestions, including
improvements in  the statement and proof of Propositions
\ref{LE-inherit} and \ref{p:AE=l-fr}.
%eli
The authors would like to thank Ethan Akin for a careful reading
of the paper and for suggesting several improvements. Finally we
thank Benjy Weiss for many helpful conversations.

 \br

\section{Topological dynamics background}
\label{Sec-TD}

 \sk

%eli
Usually all the topological spaces we deal with are assumed
to be Hausdorff and completely regular. However occasionally
we will consider a pseudometric on a space, in which case
of course the resulting topology need not be even $T_0$.
%23June
%A  $G$-{\it dynamical system\/} $(G,X)$ is a continuous action $G
%\times X \to X$ of the topological group $G$ on the
%topological space $X$. Other, synonymous, names for such
%systems will be: {\it a $G$-space or a $G$-action} or,
%when the group is understood just {\it a system or an action}.
%Most of the time we will assume that the space $X$ is compact
%\footnote{{\it compact} will mean {\it compact
%and Hausdorff}} but occasionally we deal with non
%compact systems. For instance, we can treat $X:=G$ as a
%$G$-space under the (left) regular action.
Let $G \times X \to X$ be a continuous (left) action of the
topological group $G$ on the topological space $X$. As usual we
say that $(G,X)$, or $X$ (when the group is understood), is a
$G$-\emph{space} or, a $G$-\emph{action}. Every $G$-invariant
subset $Y \subset X$ defines a $G$-\emph{subspace} of $X$. Recall
that every topological group $G$ can be treated as a $G$-space
under the left regular action of $G$ on itself. If $X$ is a
compact $G$-space then sometimes we say also a $G$-\emph{system}
or just a \emph{system}.
%3July (referee asked (see his remark (7)) to define
%"subdirect" precisely; I hope that this form is ok)
We say that a $G$-space $X$ is a {\it subdirect product} of a
class $\Gamma$ of $G$-spaces if $X$ is a $G$-subspace of a
$G$-product of some members from $\Gamma$.
%end

The notations $(X,\tau)$ and $(X,\mu)$ are
used for a topological and a  uniform space respectively. When the
acting group is the group $\Z$ of integers, we sometimes write
$(T,X)$ instead of $(\Z,X)$, where $T:X \to X$ is the
homeomorphism which corresponds to the element $1\in \Z$ (such
systems are sometimes called {\em cascades\/}). We write $g x$ for
the image of $x\in X$ under the homeomorphism $\breve{g}: X \to X$
which corresponds to $g\in G$. As usual, $Gx=\OG(x)=\{g x: g\in
G\}$ is the {\em orbit\/} of $x$ and ${\OCG}(x)={\cls} (Gx)$ is
the closure in $X$ of $\OG(x)$. If $(G,Y)$ is another $G$-system
then a surjective continuous $G$-{\it map} $\pi: X \to Y$ (that
is, $g \pi(x)=\pi (g x), \forall (g, x) \in G \times X$) is called
a {\em homomorphism\/}. We also say that $Y$ is a $G$-{\em
factor\/} of $X$. When $(G,X)$ is a dynamical system and $Y\subset
X$ is a nonempty closed $G$-invariant subset, we say that the
dynamical system $(G,Y)$, obtained by restriction to $Y$, is a
{\em subsystem\/} of $(G,X)$.

Denote by $C(X)$ the Banach algebra of all real valued bounded
functions on a topological space $X$
%23June
under the supremum norm.
Let $G$ be a topological group.
We write $RUC(G)$ for the Banach subalgebra of $C(G)$ of
{\em right uniformly
%end
continuous\/} \footnote{Some authors call these functions {\it
left uniformly continuous}} real valued bounded functions on $G$.
These are the functions which are uniformly continuous with
respect to the {\it right uniform structure} on $G$. Thus, $f\in RUC(G)$
iff for every $\ep>0$ there exists a neighborhood $V$ of the
identity element $e\in G$ such that $\sup_{g\in
G}|f(vg)-f(g)|<\ep$ for every $v \in V$. It is equivalent to say
that the orbit map $G \to C(G), \hskip 0.3cm g \mapsto {}_gf$ is
norm continuous where ${}_gf$ is the {\it left translation} of $f$
defined by ${}_gf(x)=L_g(f)(x):=f(gx)$. Analogously can be defined
the algebra $LUC(G)$ of {\it left uniformly continuous} functions
and the {\it right translations} $f_g(x)=R_g(f)(x):=f(xg)$.
It is
easy to see that $UC(G):=RUC(G) \cap LUC(G)$ is a left and right
$G$-invariant closed subalgebra of $RUC(G)$.

More generally: for a given (not necessarily compact) $G$-space
$X$ a function $f \in C(X)$ will be called {\it right uniformly
continuous} if the orbit map $G \to C(X), \ g \mapsto
{}_gf:=L_g(f)$ is norm continuous, where $L_g(f)(x):=f(gx)$.
%23June
The map $C(X) \times G \to C(X), \ (f,g) \mapsto {}_gf$ defines a
right action. The set $RUC(X)$ of all right uniformly continuous
functions on $X$ is a uniformly closed $G$-invariant subalgebra of
$C(X)$.
%end

A $G$-{\it compactification} of a $G$-space $X$ is a dense
continuous $G$-map $\nu: X \to Y$ into a compact $G$-system $Y$. A
compactification $\nu: X \to Y$ is {\it proper} when $\nu$ is a
topological embedding. We say that a $G$-compactification $\nu: G
\to S$ of $X:=G$ (the left regular action) is a {\it right
topological semigroup compactification} of $G$ if $S$ is a {\it
right topological semigroup} (that is, for every
%aref (5)
%$p \in S$
$x \in S$
the map
 $\rho_s: S \to S, \ \rho_s(x)=xs$ is continuous)
 and
$\nu$ is a homomorphism of semigroups. There exists a canonical
1-1 correspondence (see for example \cite{Vr-oldpaper}) between
the $G$-compactifications of $X$ and uniformly closed
$G$-subalgebras (``subalgebra", will always  mean a subalgebra
{\it containing the constants}) of $RUC(X)$. The $G$-compactification
$\nu: X \to Y$ induces an isometric $G$-embedding of $G$-algebras
$$
j_{\nu}: C(Y) \to RUC(X), \hskip 0.3cm \phi \mapsto \phi \circ \nu
$$
and the algebra $\Acal_{\nu}$ (corresponding to $\nu$) is defined
as the image $j_{\nu}(C(Y))$. Conversely, if $\Acal$ is a
uniformly closed $G$-subalgebra of $RUC(X)$, then its Gelfand
space $|\Acal| \subset (\Acal^*, weak^*)$ has a structure of a
dynamical system $(G,|\Acal|)$ and the map $\nu_{\Acal}: X \to
Y:=|\Acal|, \ x \mapsto eva_x$, where $eva_x(\varphi):=\varphi(x)$
is the evaluation at $x$ multiplicative functional, defines a
$G$-compactification.
%aref
%In particular, $\nu_{\Acal_{\nu}}=\nu$.
%
If $\nu_1: X \to Y_1$ and
$\nu_2: X \to Y_2$ are two $G$-compactifications then
$\Acal_{\nu_1} \subset \Acal_{\nu_2}$ iff $\nu_1= \a \circ \nu_2$
for some $G$-homomorphism $\a: Y_2 \to Y_1$. The algebra
$\Acal_{\nu}$ determines the compactification $\nu$ uniquely, up
to the equivalence of $G$-compactifications.

The $G$-algebra $RUC(X)$
% of all right uniformly continuous functions on a $G$-space $X$
defines the corresponding Gelfand
space $|RUC(X)|$ (which we denote by $\beta_G X$) and the {\it
maximal $G$-compactification} $i_{\beta}: X \to \beta_G X$. Note
that this map may not be an embedding even for Polish $X$ and $G$
(see \cite{Me-ex}); it follows that there is no proper
$G$-compactification for such $X$. If $X$ is a compact $G$-system
then $\beta_G X$ can be identified with $X$ and $C(X)=RUC(X)$.

A point $x_0\in X$ is a {\em transitive point\/}
%23June
(notation: $x_0 \in Trans(X)$) if $\OCG(x_0)=X$ and the $G$-space
$X$ is called {\em point transitive\/} (or just {\em
transitive\/}) if
%there is some transitive point in $X$.
$Trans(X) \neq \emptyset$. It is {\em topologically transitive\/}
if for every two nonempty open subsets $U,V\subset X$ there exists
$g\in G$ with $gU\cap V\ne\emptyset$. Every point transitive
$G$-space is topologically transitive.
%end
When $X$ is a metrizable system, topological transitivity is
equivalent to point transitivity and, in fact, to the existence of
a dense $G_\del$ set of transitive points.
%eli
For a $G$-space $(G,X)$ with $G$ locally compact we say that a
point $x\in X$ is a {\it recurrent point\/} if there is a net $G
\ni g_i \to \infty$ with $x=\lim_{i\to\infty} g_i x$. A system
$(G,X)$ with a recurrent transitive point is called a {\em
recurrent-transitive} system. Note that a transitive infinite
$\Z$-system is recurrent transitive iff $X$ has no isolated
points.

A system $(G,X)$ is called {\em weakly mixing\/} if the product
system $(G, X\times X)$ (where $g(x,x')=(gx,gx')$) is
topologically transitive. A system $(G,X)$ is called {\em
minimal\/} if every point of $X$ is transitive.

A triple $(G,X,x_0)$ with compact $X$ and a distinguished
transitive point $x_0$ is called a {\em pointed dynamical
system\/} (or sometimes an {\em ambit\/}). For homomorphisms $\pi:
(X,x_0) \to (Y,y_0)$ of pointed systems we require that
$\pi(x_0)=y_0$. When such a homomorphism exists it is unique. A
pointed dynamical system $(G,X,x_0)$ can be treated as a
$G$-compactification $\nu_{x_0} : G \to X, \hskip 0.3cm
\nu_{x_0}(g)=gx_0$. We associate, with every $F\in C(X)$, the
function $j_{x_0}(F)=f\in RUC(G)$ defined by $f(g)=F(g x_0)$. Then
the map $j_{x_0}$ is actually the above mentioned isometric
embedding $j_{\nu_{x_0}}: C(X) \to RUC(G)$. Let us denote its
image by $j_{x_0}(C(X))=\Acal(X,x_0)$. We have ${}_g
f={}_g(j_{x_0}(F))=j_{x_0}(F\circ g)$. The Gelfand space
$|\Acal(X,x_0)|$ of the algebra $\Acal(X,x_0)$ is naturally
identified with $X$ and in particular the multiplicative
functional ${\eva}_e:f\mapsto f(e)$, is identified with the point
$x_0$. Moreover the action of $G$ on $\Acal(X,x_0)$ by left
translations induces an action of $G$ on $|\Acal(X,x_0)|$ and
under this identification the pointed systems $(X,x_0)$ and
$(|\Acal(X,x_0)|,{\eva}_e)$ are isomorphic.

Conversely, if $\Acal$ is a $G$-invariant uniformly closed
subalgebra of $RUC(G)$ (here and in the sequel when we say that a
subalgebra of $RUC(G)$ is $G$-invariant we mean left
$G$-invariant;
%23June
that is invariant with respect to the action
$\Acal \times G \to \Acal, (f,g) \mapsto {}_gf)$,
%end
then its Gelfand space $|\Acal|$ has a structure
of a pointed dynamical system $(G,|\Acal|,{\eva}_e)$.
In particular, we have, corresponding to the algebra $RUC(G)$, the
{\em universal ambit\/} $(G,G^{R},{\eva}_e)$ where we denote the
Gelfand space
$|RUC(G)|=\beta_G G$ by $G^{ R}$ (See for example \cite{E} or
\cite{Vr} for more details).

It is easy to check that for any collection
$\{(G,X_\tet,x_\tet):\tet\in \Tet\}$ of pointed systems we have
$$
\Acal\left(\bigvee \{(X_\tet,x_\tet):\tet\in \Tet\}\right)
=\bigvee \{\Acal(X_\tet,x_\tet):\tet\in \Tet\},
$$
where $\bigvee \{(X_\tet,x_\tet):\tet\in \Tet\}$ is the orbit
closure of the point $x$ in the product space $\prod_{\tet\in
\Tet}X_\tet$ whose $\tet$ coordinate is $x_\tet$, and the algebra
on the right hand side is the closed subalgebra of  $RUC(G)$
generated by the union of the subalgebras $\Acal(X_\tet,x_\tet)$.

\begin{defn} \label{d:coming}
\ben
\item
We say that a function $f\in C(X)$ on a $G$-space $X$ \emph{comes
from} a $G$-system $Y$ if there exist a $G$-compactification $\nu:
X \to Y$ (so, $\nu$ is onto if $X$ is compact) and a function
$F\in C(Y)$ such that $f=\nu \circ F$ (equivalently, $f \in
\Acal_{\nu}$). Then necessarily $f \in RUC(X)$. Only the maximal
$G$-compactification $i_{\beta}: X \to \beta_G X$ has the property
that every $f \in RUC(X)$ comes from $i_{\beta}$.
\item
A function $f\in C(G)$ {\em comes from\/} a pointed system
$(Y,y_0)$
%aref
(and then necessarily $f \in RUC(G)$)
%end
if for some continuous
function $F\in C(Y)$ we have $f(g)=F(g y_0),\ \forall g\in G$;
i.e.$f=j_{y_0}(F)$
%aref
(equivalently, if $f \in \Acal(Y,y_0)$).
%end
Defining $\nu: X=G \to Y$ by $\nu(g)=gy_0$ observe that this is
indeed a particular case of \ref{d:coming}.1.
\item
A function $f\in RUC(X)$ is called {\em minimal\/} if it comes
from a minimal system.
\een
\end{defn}

\br

\section{The enveloping semigroup}\label{es}

\sk

The {\em enveloping (or Ellis) semigroup\/} $E=E(G,X)=E(X)$ of a
dynamical system $(G,X)$ is defined as the closure in $X^X$ (with
its compact, usually non-metrizable, pointwise convergence
topology) of the set $\breve{G}=\{\breve{g}: X \to X\}_{g \in G}$
considered as a subset of $X^X$. With the operation of composition
of maps this is a right topological semigroup.
Moreover, the map $i: G \to E(X), g \mapsto \breve{g}$ is a right
topological semigroup compactification of $G$.

\begin{prop} \label{envel}
The enveloping semigroup of a dynamical system $(G,X)$ is
isomorphic (as a dynamical system) to the pointed product
$$
(E',\om_0)=\bigvee \{({\OCG}(x),x): x\in X\}\subset X^X.
$$
\end{prop}
\begin{proof}
It is easy to see that the map $p \mapsto p\om_0,\ (G,E,i(e)) \to
(G,E',\om_0)$ is an isomorphism of pointed systems.
\end{proof}

Let $X$ be a (not necessarily compact) $G$-space.
%aref  I guess this construction easily goes also in the case of
%$f \in RUC(X,K)$ where
%$K$ is a compact space which replaces our $I$. Then $\Omega:=K^G$ ...
Given $f\in RUC(X)$ let $I=[-\|f\|,\|f\|]\subset \R$ and
$\Om=I^G$, the product space equipped with the compact product
topology. We let $G$ act on $\Om$ by $g\om(h)=\om(hg),\ g,h\in G$.

Define the continuous map
$$
f_{\sharp}: X \to \Om, \hskip 0.3cm f_{\sharp}(x)(g)=f(gx)
$$
and the closure $X_f:= {\cls} (f_{\sharp}(X))$ in $\Om$.
Note that $X_f
= f_{\sharp}(X)$ whenever $X$ is compact.

Denoting the unique continuous
extension of $f$ to $\beta_G X$ by $\tilde f$ we now define a map
$$
\psi: \beta_G X \to X_f \quad \text{by} \quad \psi(y)(g)= \tilde
f(gy), \qquad y\in \beta_G X, g\in G.
$$
Let
%3July replacing \pi_e by pr_e (too many \pi_s !)
$pr_e: \Om \to \R$ denote the  projection of $\Om=I^G$ onto the
$e$-coordinate and let $F_e:=pr_e \rest_{X_f}: X_f \to \R$ be its
%end
restriction to $X_f$. Thus, $F_e(\om):=\om(e)$ for every $\om \in
X_f$.

%23June
For every $f \in RUC(X)$ denote by $\Acal_f$ the smallest closed
$G$-invariant subalgebra of $RUC(X)$ which contains $f$.
There is then a naturally defined $G$-action on the Gelfand space
$|\Acal_f|$ and a $G$-compactification
(homomorphism of dynamical systems if $X$ is
compact) $\pi_f: X \to |\Acal_f|$.
Next consider the map $\pi: \beta_G X \to |\Acal_f|$,
the canonical extension of
 $\pi_f$.
 %: X \to |\Acal_f|$
%end

%aref (8)
The action of $G$ on $\Omega$ is not in general continuous.
However, the restricted action on $X_f$ is continuous for every $f
\in RUC(X)$.
%eli
This follows from the second assertion of the
next proposition.

\begin{prop} \label{p:generalX_f}
\begin{enumerate}
\item
Each $\om \in X_f$
is an element of $RUC(G)$.
\item
The map $\psi: \beta_G X \to X_f$ is a continuous homomorphism of
$G$-systems.
The dynamical system $(G,|\Acal_f|)$ is isomorphic to $(G,X_f)$
and the diagram
%gl-
\begin{equation*}
\xymatrix {
X \ar[d]_{\pi_f} \ar[r]^{i_{\beta}} \ar[dr]^-{\pi}  &
\beta_G X \ar[dl]^{f_{\sharp}} \ar[d]^{\psi} \ar[dr]^{\tilde{f}} & \\
|\Acal_f| \ar[r]  &  X_f \ar[l] \ar[r]^{F_e} & \R & }
\end{equation*}
commutes.
\item $f=F_e \circ f_{\sharp}$.
Thus every $f \in RUC(X)$
%aref
%is coming
comes from the system $X_f$. Moreover, if $f$
%is coming
comes from a system $Y$ and a $G$-compactification $\nu: X \to Y$
 then there exists a homomorphism $\a: Y \to
X_f$ such that $f_{\sharp}=\a \circ \nu$.
%aref
In particular, $f \in \Acal_f \subset \Acal_{\nu}$.
\end{enumerate}
\end{prop}

\begin{proof}
1.\ $f \in RUC(X)$ implies that $f_{\sharp}(X)$ is a uniformly
equicontinuous subset of $I^G$ (endowing $G$ with its right
uniform structure). Thus, the pointwise closure ${\cls} (f_{\sharp}(X)) =
X_f$ is also uniformly equicontinuous. In particular, for every
$\omega \in X_f$ the function $\omega: G \to I$ is right uniformly
continuous.

2.\ Suppose $i_{\beta}(x_\nu) \in i_{\beta}(X)$ is a net
converging to $y \in \beta_G X$. Then $\psi(y)(g)={\tilde
f}(gy)=\lim_\nu f(g x_\nu) =\lim_\nu f_{\sharp}(x_\nu)(g)$. Thus
$\psi(y)=\lim_\nu f_{\sharp}(x_\nu)$ is indeed an element of $X_f$
and it is easy to see that $\psi$ is a continuous $G$-homomorphism.
In particular, we see that $X_f$,
being a $G$-factor of $\beta_G X$, is indeed a $G$-system (i.e.
the $G$-action on $X_f$ is jointly continuous).

%14June
Now we use the map
$\pi: \beta_G X \to |\Acal_f|$.
By definition, the elements of $\beta_G X$ are continuous
multiplicative linear functionals on the algebra $RUC(X)$, and for
$y\in \beta_G X$ its value $\pi(y) \in |\Acal_f|$ is the
restriction $y \rest_{\Acal_f}$ to the subalgebra $\Acal_f \subset
RUC(X)$. For $g\in G$, as above, let $_g f\in \Acal_f \subset
RUC(X)$ be defined by $_g f(x)=f(gx)$. Then $\pi(y_1)=\pi(y_2)$
implies $y_1(_g f)=\tilde f(gy_1)=\tilde f(gy_2)=y_2(_g f)$ for
every $g\in G$.

Conversely, assuming, $\tilde f(gy_1)=\tilde f(gy_2)$ for every
$g\in G$, we observe that, as $y_1$ and $y_2$ are multiplicative
functionals, we also have $y_1(h)=y_2(h)$ for every $h$ in the
subalgebra $\Acal_0$ generated by the family $\{_g f: g\in G\}$.
Since $\Acal_0$ is dense in $\Acal_f$ and as $y_1$ and $y_2$ are
continuous we deduce that $\pi(y_1)=y_1 \rest_{\Acal_f}=y_2
\rest_{\Acal_f}=\pi(y_2)$.

We clearly have $\psi(y_1)=\psi(y_2) \iff \tilde f(gy_1)= \tilde
f(gy_2)$ for every $g\in G$. Thus for $y_1,y_2 \in \beta_G X$ we
have $\pi(y_1)=\pi(y_2) \iff \psi(y_1)=\psi(y_2) \iff \tilde
f(gy_1)= \tilde f(gy_2)$ for every $g\in G$, and we find that
indeed $|\Acal_f|$ and $X_f$ are isomorphic $G$-systems.

%gl-
The verification of the commutativity of the diagram is straightforward.

3.\ Clearly, $F_e ( f_{\sharp}(x))=f_{\sharp}(x)(e)=f(ex)=f(x)$
for every $x \in X$. For the rest use the $G$-isomorphism
$|\Acal_f| \leftrightarrow X_f$ (assertion
%me  "3" before
2). If $f=F \circ \nu$ for some $F \in C(Y)$ then $f \in
\Acal_{\nu}$. This implies the inclusion of $G$-subalgebras
$\Acal_f \subset \Acal_{\nu}$ which leads to the desired
$G$-homomorphism $\a: Y \to X_f$.
\end{proof}

\begin{rmk} \label{r:X_f}
\ben
\item
Below we use the map $f_{\sharp}: X \to X_f$ and Proposition
\ref{p:generalX_f} in two particular cases. First, for a compact
$G$-space $X$ when clearly $\beta_G X$ can be replaced by $X$. We
frequently consider also the case of left regular action of $G$ on
itself $X:=G$ (see Proposition \ref{X_f}). Here the canonical
maximal $G$-compactification $i_{\beta}: X \to \beta_G X$ is
actually the compactification $G \to G^R$ and the orbit
%14June
$Gf= \{R_g(f)\}_{g \in G}=f_{\sharp}(G)$ of $f \in RUC(G)$
is pointwise dense in $X_f={\cls} (f_{\sharp}(G)) \subset
\Om=I^G$.
\item
$\beta_G X$ is a subdirect product
%(i.e. a $G$-subspace of a
%$G$-product)
of the $G$-systems $X_f$ where $f \in RUC(X)$. This
follows easily from Proposition \ref{p:generalX_f} and the fact
that elements of $C(\beta_G X)=\{{\tilde f}: f \in RUC(X) \}$
separate points and closed subsets of $\beta_G X$.
\item
Proposition \ref{p:generalX_f}.3 actually says that the
compactification $f_{\sharp}: X \to X_f$ is {\it minimal}
(in fact, {\it the smallest}) among all the $G$-compactifications
$\nu: X \to Y$ such that $f \in RUC(X)$
%aref
%is coming
comes from $\nu$. The {\it maximal} compactification in the same
setting is clearly $i_{\beta}: X \to \beta_G X$.
 \een
\end{rmk}

\begin{prop} \label{X_f}
\ben
\item
Consider
%oct21
the
left regular action of $G$ on itself $X:=G$.
For every $f \in RUC(G)$ we have $Gf \subset X_f= \OCG(f) \subset
\Om$, $f_{\sharp}(e)=f$ and
$F_e(g f)=f(g)$ for every $g \in G$.
\item
The pointed $G$-system $(|\Acal_f|,{\eva}_e)$ is isomorphic to
$(X_f, f)$
%3July
(hence $\Acal_f=\Acal(X_f,f)$).
%end
\item
$f=F_e \circ f_{\sharp}$. Thus every $f \in RUC(G)$
%aref
%is coming
comes from the pointed system $(X_f,f)$. Moreover, if $f$
%aref
%is coming
comes from a pointed system $(Y,y_0)$ and $\nu: (G,e) \to (Y,y_0)$
is the corresponding $G$-compactification then there exists a
homomorphism $\a: (Y,y_0) \to (X_f,f)$ such that $f_{\sharp}=\a
\circ \nu$. In particular, $f \in \Acal_f \subset \Acal (Y,y_0)$.
\item
Denote by $X_f^H\subset I^H$ the dynamical system constructed for
the subgroup $H < G$ and the restriction $f \rest_H$ (e.g.,
$X_f^G=X_f$). If $H< G$ is a dense subgroup then, for every $f\in
RUC(G)$, the dynamical systems $(H,X_f) $ and $(H,X^H_f)$ are
canonically isomorphic.
\een
\end{prop}
\begin{proof} For the assertions 1, 2 and 3
use Proposition \ref{p:generalX_f} and Remark \ref{r:X_f}.1.

4.\ Let $j: X_f \to X_f^H$ be the restriction of the natural
projection $I^G \to I^H$. Clearly, $j:(H,X_f)\to (H,X_f^H)$ is a
surjective homomorphism. If $j(\om)=j(\om')$ then $\om(h)=\om'(h)$
for every $h\in H$. Since, by Proposition \ref{p:generalX_f}.1
every $\om \in X_f$ is a
continuous function on $G$ and since we assume that $H$ is dense
in $G$ we conclude that $\om=\om'$ so that $j$ is an isomorphism.
\end{proof}

\begin{defn} \label{d:p-univ}
We say that a pointed dynamical system $(G,X,x_0)$ is {\em
point-universal\/} if it has the property that for every $x\in X$
there is a homomorphism $\pi_x:(X,x_0)\to (\OCG(x),x)$. A closed
$G$-invariant subalgebra $\Acal\subset RUC(G)$ is called {\em
point-universal\/} if the corresponding Gelfand system
$(G,|\Acal|,\eva_e)$ is point-universal.
\end{defn}

\begin{prop}\label{univ}
The following conditions on the pointed
dynamical system $(G,X,x_0)$ are equivalent:
\begin{enumerate}
\item
$(X,x_0)$ is point-universal.
\item
$\Acal(X,x_0) = \bigcup_{x \in X} \Acal (\OCG(x), x)$.
\item
$(X,x_0)$ is isomorphic to its enveloping semigroup $(E(X),i(e))$.
\end{enumerate}
\end{prop}

\begin{proof}
$1 \Rightarrow 2:$\ Clearly, $\Acal(X,x_0) = \Acal(\OCG(x_0), x_0)
\subset \bigcup_{x \in X} \Acal (\OCG(x), x)$. Suppose
$f(g)=F(gx)$ $\forall g \in G$ for some $F\in C(\OCG(x))$ and
$x\in X$. Since $(X,x_0)$ is point-universal there exists a
homomorphism $ \pi_x: (X,x_0) \to (\OCG(x),x)$. Hence
$f(g)=F(gx)=F(g \pi_x(x_0))=F(\pi_x (gx_0)) =
 (F\circ \pi_{x})(gx_0)=j_{x_0}(F\circ \pi_x)(g)$ and we conclude that
$f =j_{x_0}(F\circ \pi_x) \in \Acal(X,x_0)$.
%$\Acal(X,x_0)= \bigcup_{x \in X} \Acal (\OCG(x), x)$.

$2 \Rightarrow 3:$\ Proposition \ref{envel} guarantees the
existence of a pointed isomorphism between the systems $(E(X),i(e))$
and $\bigvee_{x\in X}(\OCG(x),x)$. Now, using our assumption we
have:
$$
\Acal(E(X),i(e))=\Acal\left(\bigvee_{x\in X}(\OCG(x),x)\right)=
\bigvee_{x\in X}\Acal(\OCG(x),x)= \Acal(X,x_0),
$$
whence the isomorphism of $(X,x_0)$ and $(E(X),i(e))$.

$3 \Rightarrow 1:$\ For any fixed $x\in X$ the map $\pi_x:E(X)\to
X$, defined by $\pi_x(p)=px$, is a $G$-homomorphism with
$\pi_x(i(e))=x$. Our assumption that  $(X,x_0)$ and $(E(X),i(e))$ are
isomorphic now implies the point-universality of $(X,x_0)$.
\end{proof}

%3July (replacing Corollaries by Propositions)
%We now formulate three important corollaries.

\begin{prop}\label{Cor-univ1}
A transitive
%aref
%pointed (mentioning "pointed" is necessary or not ?)
%Or maybe "transitive" can be deleted ?
system $(G,X,x_0)$ is point-universal iff the map $G \to
X,\ g\mapsto gx_0$ is a right topological semigroup
compactification of $G$.
\end{prop}
\begin{proof}
The necessity of the condition follows directly from Proposition
\ref{univ}. Suppose now that the map $G \to X,\ g\mapsto gx_0$ is
a right topological semigroup compactification of $G$. Given $x\in
X$ we observe that the map $\rho_x: (X,x_0) \to (X,x)$,
$\rho_x(z)=zx$ is a homomorphism of pointed systems, so
that $(G,X,x_0)$ is point-universal.
\end{proof}

In particular, for every $G$-system $X$ the enveloping semigroup
$(E(X), i(e))$, as a pointed $G$-system, is point-universal. Here,
%aref
as before,
 $i: G \to E(X), \hskip 0.2cm g \mapsto \breve{g}$ is
the canonical enveloping semigroup compactification.

%3July
\begin{prop}\label{Cor-univ2}
Let $(G,X,x_0)$ be a pointed compact system and
$\Acal=\Acal(X,x_0)$ the corresponding (always left $G$-invariant)
subalgebra of $RUC(G)$. The following conditions are equivalent:
\begin{enumerate}
\item
$(G,X,x_0)$ is point-universal.
%3July (perhaps now is more clear why this condition $X_f \subset \Acal$
%appears naturally in the paper of Veech and in our second paper
%end
\item
$X_f \subset \Acal$ for every $f \in \Acal$ (in particular,
$\Acal$ is also right $G$-invariant).
\end{enumerate}
\end{prop}
\begin{proof}
$1 \Rightarrow 2$:
\ Let $f: G \to \R$
%gl++
belong to $\Acal$.
Consider the $G$-compactification $f_{\sharp}: G \to X_f:={\cls}
(Gf)$ as defined by Proposition \ref{X_f}. We have to show that
$\varphi \in \Acal$ for every $\varphi \in X_f$. Consider the
orbit closure $X_{\varphi}= {\cls} (G\varphi)$ in $X_f$. By
Definition \ref{d:coming}.2 there exists a continuous function
%gl++
$F: X \to \R$ such that $f(g)=F(gx_0)$ for every $g \in G$.
That is, $f$ comes from the pointed system $(X,x_0)$.
For some net $g_i\in G$ we have $\varphi(g)=
\lim_i f(gg_i)$ for every $g \in G$ and with no loss in generality
we have $x_1=\lim_i g_ix_0\in X$.
Then
$$
\varphi(g)= \lim_i f(gg_i) = \lim_i F(gg_i x_0) =
F(gx_1).
$$
Thus $\varphi$ comes from the pointed system $(\OCG(x_1), x_1)$
and in view of Proposition \ref{univ} we conclude that
indeed $\varphi\in \Acal$.
%
%
%
%
%$f^*: X \to \R$ such that $f(g)=f^*(gx_0)$ for every $g \in G$.
%That is, $f$ comes from the pointed system $(X,x_0)$. By
%Proposition \ref{X_f}.3 there exists a homomorphism $\a: (X,x_0)
%\to (X_f,f)$ such that $f_{\sharp}=\a \circ \nu$ (where
%$\nu(g)=gx_0$). One can choose a point $z_0 \in X$ such that
%$\a(\OCG(z_0))=X_{\varphi}$ and $\a(z_0)=\varphi$. Let $F_e:
%X_{\varphi} \to \R, \ \ F_e(\omega)=\omega(e)$ be the map defined
%as in Proposition \ref{X_f}.3 (for $X_{\varphi}$). Then we get
%$$F_e(\a(gz_0))=\a(gz_0)(e)=(g\a(z_0))(e)=\a(z_0)(g)=\varphi(g).$$
%Therefore, $\varphi$ comes
%from the pointed system $(\OCG(z_0),z_0)$.  This means that
%$\varphi \in \Acal(\OCG(z_0),z_0)$. Since $(X,x_0)$ is point
%universal, Proposition \ref{univ} guarantees that
%$\Acal(\OCG(z_0),z_0) \subset \Acal(X,x_0)$. Thus, $\varphi \in
%\Acal(X,x_0)$.

$2 \Rightarrow 1$: \ Define the $G$-ambit
$$
(Y,y_0):=\bigvee\{(X_f, f): f\in \Acal\}.
$$
First we show that $\Acal(X,x_0) = \Acal(Y,y_0)$. Indeed, as we
know
$$
\Acal(Y,y_0)=\bigvee\{\Acal(X_f, f): f\in \Acal\}.
$$
Proposition \ref{X_f} implies that $f \in \Acal_f=\Acal(X_f,f)$
for every $f \in \Acal(X,x_0)$. Thus we get
$$
f \in \Acal_f=\Acal(X_f,f) \subset \Acal(Y,y_0) \ \ \ \forall f
\in \Acal(X,x_0).
$$
Therefore, $\Acal(X,x_0) \subset \Acal(Y,y_0)$. On the other hand,
$\Acal_f=\Acal(X_f,f) \subset \Acal(X,x_0)$ (for every $f \in
\Acal(X,x_0)$) because $\Acal(X,x_0)$ is left $G$-invariant and
$\Acal_f$ is the smallest closed left $G$-invariant subalgebra of
$RUC(G)$ which contains $f$. This implies that $\Acal(Y,y_0)
\subset \Acal(X,x_0)$. Thus, $\Acal(X,x_0) = \Acal(Y,y_0)$. Denote
this algebra simply by $\Acal$.

%by construction, $\Acal \subset \Acal(Y,y_0)$. On the other hand
%the algebra $\Acal$ is generated (as a closed invariant algebra)
%by the functions of $K(G)$ and by part 1 it follows that
%$K(G)=\Acal$.

Suppose $py_0=qy_0$ for $p,q\in E(Y)$ (the enveloping semigroup of
$(G,Y)$). By our assumption, $X_f \subset \Acal$ for every $f \in
\Acal$. Then every $y\in Y$, considered as an element of the
product space $\prod_{f \in \Acal} X_f$, has the property that its
$f$-coordinate, say $y_f$ is again an element of $\Acal$ and it
follows that $y_f$ appears as a coordinate of $y_0$ as well.
Therefore also $py_f=qy_f$ and it follows that $py=qy$. Thus the
map $p \mapsto py_0$ from $(E(Y),i(e))$ to $(Y,y_0)$ is an
isomorphism. By Proposition \ref{univ}, $(Y,y_0)$ (and hence also
$(X,x_0)$) is point-universal.

%gl+++
(Observe that $Gf=\{R_g(f)\}_{g \in G} \subset X_f:={\cls}(Gf)$.
Therefore, the condition $X_f \subset \Acal$, $\forall f \in
\Acal$ trivially implies that $\Acal$ is right invariant.)
\end{proof}

\begin{prop}\label{Cor-univ3}
Let $P$ be a property of compact $G$-dynamical systems which is
preserved by products, subsystems and $G$-isomorphisms.
%such that the one point
%trivial $G$-system has the property $P$.
\begin{enumerate}
\item
Let $X$ be a (not necessarily compact) $G$-space and let $\Pcal_X
\subset C(X)$ be the collection of functions coming from systems
having property $P$. Then there exists a maximal
$G$-compactification $X^{\Pcal}$ of $X$ with property $P$.
Moreover, $j(C(X^{\Pcal}))=\Pcal_X$. In particular, $\Pcal_X$ is a
uniformly closed, $G$-invariant subalgebra of $RUC(X)$.
%If $X$ is compact then $(G,X^{\Pcal})$ is the maximum factor of
%$(G,X)$ with property $P$.
\item
Let $\Pcal \subset C(G)$ be the set of functions coming from
systems with property $P$. Then $(G^{\Pcal}, \eva_e)$ is the
universal point transitive compact $G$-system having property $P$.
Moreover $\Pcal$ is a point-universal  subalgebra of $RUC(G)$.
(Thus, $\Pcal$ is uniformly closed, right and left $G$-invariant
and $X_f \subset \Pcal$ for every $f \in \Pcal$.)
% and $\Pcal$ defines a factor $G^{\Pcal}$ of $G^R$ which is the maximum
% factor of $G^R$ with property $P$.

%; i.e. if $(G,Z)$ is any dynamical system with property $P$ and
%$z_0\in Z$ is a transitive point then there exists a homomorphism
%$\phi: G^{\Pcal}\to Z$ mapping the distinguished point
%${\eva}_e\in G^{\Pcal}$ onto $z_0$.
\item
If in addition $P$ is preserved by factors then $f \in \Pcal$ iff
$X_f$ has property $P$.
\end{enumerate}
\end{prop}

\begin{proof}
1.\ We only give an outline of the rather standard procedure.
There is a complete \emph{set} $\{\nu_i: X \to Y_i \}_{i \in I}$
of equivalence classes of $G$-compactifications of $X$ such that
each $Y_i$ has the property $P$. Define the desired
compactification $\nu: X \to Y=\cls (\nu (X)) \subset \prod_{i \in
I}Y_i$ via the diagonal product. Then we get the suprema of our
class of $G$-compactifications.  In fact, $Y$ has the property $P$
because the given class is closed under subdirect products. $f \in
\Pcal$ means that it comes from some $Y_i$ via the
compactification $\nu_i: X \to Y_i$. Denote $Y$ by $X^{\Pcal}$.
Now using the natural projection of $Y$ on $Y_i$ it follows that
$f$ comes from $Y=X^{\Pcal}$. This implies the coincidence
$j(C(X^{\Pcal}))=\Pcal_X$.

2.\ The construction of the maximal ambit $(G^\Pcal, \eva_e)$ with
the property $P$ is similar. In fact it is a particular case of
the first assertion identifying $G$-ambits $(Y,y_0)$ and
$G$-compactifications $\nu_{y_0} : G \to Y, \hskip 0.2cm
\nu_{y_0}(g)=gy_0$ of $X:=G$. As to the point-transitivity of
$\Pcal$ note that according to the definition the uniformly closed
subalgebra $\Pcal \subset RUC(G)$ is the set of functions coming
from systems with property $P$. Every subsystem of $G^{\Pcal}$ has
the property $P$. In particular, $(\OCG(x), x)$ has the property
$P$. Therefore, $\Pcal$ contains the algebra $\Acal(\OCG(x), x)$
for every $x \in X$. By Proposition \ref{univ} it follows that
$\Pcal$ is point-universal. Thus Proposition \ref{Cor-univ2}
guarantees that $X_f \subset \Pcal$ for every $f \in \Pcal$ (and
that $\Pcal$ is right and left $G$-invariant).

3.\ Use Proposition \ref{p:generalX_f}.3.
\end{proof}

\br

\section{A dynamical version of the
Bourgain-Fremlin-Talagrand theorem}\label{Sec-BFT}

\sk

Let $E=E(X)$ be the enveloping semigroup of a $G$-system $X$. For
every $f \in C(X)$ define
$$
E^f:=\{p_f: X \to \R \}_{p \in E}= \{f\circ p: p\in E\}, \quad
p_f(x)=f(px).
$$
Then $E^f$ is a pointwise compact subset of $\R^X$, being a
continuous image of $E$ under the map $q_f: E \to E^f, \hskip
0.2cm p \mapsto p_f$.

%aref (11)
%we do not need at all this right action so why we should define it ?
%We have a natural separately continuous {\it right} action
%$$E^f
%\times G \to E^f, \hskip 0.2cm (p_f, g) \mapsto (pg)_f.$$

%aref [I think that these 4 lines are superfluous or even misleading]
%When $G$ is a commutative group then $E^f$ can be treated also as
%a left $G$-system. In this case the action $G \times E^f \to E^f$
%is jointly continuous; so that moreover, the map $E \to E^f,
%\hskip 0.2cm p \mapsto p_f$ is a homomorphism of $G$-systems.
%

\sk

Recall that a topological space $K$ is {\it Rosenthal compact}
\cite{Godefroy} if it is homeomorphic to a pointwise compact
subset of the space $B_1(X)$ of functions of the first Baire class
on a Polish space $X$. All metric compact spaces are Rosenthal. An
example of a separable non-metrizable Rosenthal compact is the
{\it Helly compact} of all (not only strictly) increasing selfmaps
of $[0,1]$ in the pointwise topology. Another is the ``two arrows"
space of Alexandroff and Urysohn (see Example \ref{two-arr}). A
topological space $K$ is {\em angelic} if the closure of every
subset $A \subset K$ is the set of limits of sequences from $A$
and every relatively countably compact set in $K$ is relatively
compact.
%aref (12) (this change here is not exactly what
%referee wants - let's discuss it in our conversation
Note that the second condition is superfluous if $K$ is compact.
Clearly, $\beta \N$ the Stone-\v{C}ech
compactification of the natural numbers $\N$,
%aref
is not angelic, and hence it
cannot be embedded into a Rosenthal compact space.

The following theorem is due to Bourgain-Fremlin-Talagrand
\cite[Theorem 3F]{BFT}, generalizing a result of Rosenthal. The
second assertion (BFT dichotomy) is presented as in the book of
Todor\u{c}evi\'{c} \cite{T-b} (see Proposition 1 of Section 13).

\begin{thm} \label{BFT}
\begin{enumerate}
\item
Every Rosenthal compact space $K$ is angelic.
\item (BFT dichotomy)
Let $X$ be a Polish space and let $\{f_n\}_{n=1}^\infty \subset
C(X)$ be a sequence of real valued functions which is {\em
pointwise bounded}  (i.e. for each $x\in X$ the sequence
$\{f_n(x)\}_{n=1}^\infty$ is bounded in $\R$).
Let $K$ be the pointwise closure of $\{f_n\}_{n=1}^\infty$
in $\R^X$. Then either $K \subset B_1(X)$ (i.e. $K$ is
Rosenthal compact) or $K$ contains a homeomorphic copy of $\beta\N$.
\end{enumerate}
\end{thm}

Next we will show how the BFT dichotomy leads to a corresponding
dynamical dichotomy (see also \cite{Ko}). In the proof we will use
the following observation. Let $G$ be an arbitrary topological
group. For every compact $G$-space $X$, denote by $j: G \to
{\Homeo}(X), g \mapsto \breve{g}$ the associated (always {\em continuous\/})
homomorphism into the group of all selfhomeomorphisms of $X$. Then
the topological group $\breve{G}=j(G)$ (we will call it the {\em
natural restriction\/}) naturally acts on $X$. If $X$ is a compact
metric space then ${\Homeo}(X)$, equipped with the topology of
uniform convergence, is a Polish group. Hence, the subgroup
$\breve{G}=j(G)$ is second countable. In particular one can always
find a countable dense subgroup $G_0$ of $\breve{G}$.

\begin{thm}[A dynamical BFT dichotomy]\label{D-BFT}
Let $(G,X)$ be a metric dynamical system and let $E=E(X)$ be its
enveloping semigroup. We have the following alternative. Either
\begin{enumerate}
\item
$E$ is a separable Rosenthal compact (hence ${\card}{E} \leq
2^{\aleph_0}$); or
\item
the compact space $E$ contains a homeomorphic
copy of $\beta\N$, hence ${\card}{E} = 2^{2^{\aleph_0}}$.
\end{enumerate}
The first possibility holds iff $E^f$ is a Rosenthal compact for
every $f\in C(X)$.
\end{thm}

\begin{proof}
Since $X$ is compact and metrizable one can choose a
sequence $\{f_n \}_{n \in \N}$ in $C(X)$ which separates
the points of $X$. For every pair
%aref
%%p, q$ of distinct elements of $E$ there exist a point $x_0 \in X$
%and a function $f_{n_0}$ from our sequence such that
%$f_{n_0}(px_0) \neq f_{n_0}(qx_0)$. It follows that the continuous
$s, t$ of distinct elements of $E$ there exist a point $x_0 \in X$
and a function $f_{n_0}$ from our sequence such that
$f_{n_0}(sx_0) \neq f_{n_0}(tx_0)$. It follows that the continuous
diagonal map
$$
\Phi: E \to \prod_{n\in \N} E^{f_n},\qquad
p\mapsto (f_1\circ p, f_2\circ p,
\dots )
$$
separates the points of $E$ and hence is a topological embedding.
%aref
%(For $G$ commutative, clearly $\Phi$ is an isomorphism of
%dynamical systems.)

Now if for each $n$ the space $E^{f_n}$ is a Rosenthal compact
then so is $E\cong \Phi(E)\subset \prod_{n=1}^\infty E^{f_n}$,
because the class of Rosenthal compacts is closed under countable
products and closed subspaces. On the other hand
%aref
% -- since for each $f\in C(X)$
the map $q_f: E \to E^f, p \mapsto f \circ p$ is a continuous
surjection for each $f\in C(X)$. Therefore,
$E^f = {\cls}(q_f(G_0))={\cls}\{f \circ g : g\in G_0\}$, where
$G_0$ is a countable dense subgroup of $\breve{G}$.
%aref
%--- by the BFT theorem, Theorem \ref{BFT}, if at least one
By Theorem \ref{BFT} (BFT dichotomy), if at least one $E^{f_n}$ is
not Rosenthal then it contains a homeomorphic copy of $\beta\N$
and it is easy to see that so does its preimage $E$. (In fact if
$\beta\N \cong Z\subset E^{f_n}$ then any closed subset $Y$ of $E$
which projects onto $Z$ and is minimal with respect to these
properties is also homeomorphic to $\beta\N$.)

Again an application of the BFT
%aref
%theorem
dichotomy yields the fact that in the first case $E$ is angelic.
Clearly, the cardinality of every separable angelic space is at
most $2^{\aleph_0}$. Now in order to complete the proof observe
that for every compact metric $G$-system $X$ the space $E$, being
the pointwise closure of $\breve{G}$ in $X^X$, is separable, hence
${\card} E \le 2^{2^{\aleph_0}}$.

The last assertion clearly follows from the above proof.
\end{proof}

\br

\section{Metric approximation of dynamical systems} \label{sect:quasi}

\sk

Let $(X,\mu)$ be a uniform space and let $\ep \in \mu$. We say
that $X$ is $\ep$-{\em Lindel\"of{\/}} if the uniform cover $\{\ep
(x) \bigm| x\in X\},$ where $\ep(x) = \{y\in X\bigm| (x,y) \in \ep
\}$, has a countable subcover. If $X$ is $\ep$-Lindel\"of for each
$\ep \in \mu,$ then it is called {\em{uniformly Lindel\"of\/}}
\cite{Mefr}. We note that $(X,\mu)$ is uniformly Lindel\"of iff it
is $\aleph_0$-{\em {precompact\/}} in the sense of Isbell
\cite{Isb}. If $X$, as a topological space, is either separable,
Lindel\"of or {\em {ccc}\/} (see \cite[p. 24]{Isb}), then
$(X,\mu)$ is uniformly Lindel\"of. For a metrizable uniform
structure $\mu, (X,\mu)$ is uniformly Lindel\"of iff $X$ is
separable. Uniformly continuous maps send uniformly Lindel\"of
subspaces onto uniformly Lindel\"of subspaces.

A topological group $G$ is $\aleph_0$-{\it bounded} (in the sense
if Guran \cite{Gu}) if for every neighborhood $U$ of $e$ there
exists a countable subset $C \subset G$ such that $G=CU$. Clearly,
$G$ is $\aleph_0$-bounded means exactly that $G$ is uniformly
Lind\"elof with respect to its right (or, left) uniform structure. By
\cite{Gu} a group $G$ is $\aleph_0$-bounded iff $G$ is a
topological subgroup of a product of second countable topological
groups. If $G$ is either separable or Lindel\"of
($\sigma$-compact, for instance) then $G$ is uniformly Lindel\"of.

Recall our notation for the ``natural restriction"
$\breve{G}=j(G)$,
where for a compact $G$-system $(G,X)$, the
map $j: G \to {\Homeo}(X)$ is
the associated continuous homomorphism of $G$
into the group of all selfhomeomorphisms of $X$
(see Section \ref{Sec-BFT}).

We say that a compact $G$-system $X$ is {\em m-approximable\/} if
it is a subdirect product of {\em metric\/} compact $G$-systems
(see also the notion of {\it quasi-separablity} in the sense of
\cite{key,Vr}). By Keynes \cite{key}, every transitive system $X$
with $\sigma$-compact acting group $G$ is m-approximable. The
following generalization provides a simple criterion for
m-approximability.

\begin{prop} \label{prop:quasisep}
Let $X$ be a compact $G$-system. The following conditions are
equivalent:
\begin{enumerate}

\item $X$ is an inverse limit of metrizable compact $G$-systems (of
dimension $\leq dimX$);

\item $(G,X)$ is m-approximable;

\item $\breve{G}$ is uniformly Lindel\"of.
\end{enumerate}
\end{prop}
\begin{proof}
1 $\Rightarrow$ 2 is trivial.

For 2 $\Rightarrow$ 3 observe that for every metric compact
$G$-factor $X_i$ of $X$ the corresponding natural restriction $G_i
\subset {\Homeo}(X_i)$ of $G$ is second countable with respect to
the compact open topology. By our assumption it follows that the
group $\breve{G} \subset {\Homeo}(X)$ can be topologically
embedded into the product $\prod_i G_i$ of second countable
groups. Hence $\breve{G}$ is uniformly Lindel\"of by the theorem
of Guran mentioned above.

The implication 3 $\Rightarrow$ 1 has been proved (one can assume
that $G=\breve{G}$) in \cite[p. 82]{Me-diss}, and \cite[Theorem
2.19]{Me-sing} (see also \cite[Lemma 10]{usp}).
\end{proof}

\begin{prop} \label{prop:blocks}
Let $G$ be a topological group. The following conditions are
equivalent.
\begin{enumerate}
\item
$G$ is uniformly Lindel\"of.
\item
The greatest ambit $G^{ R}$ is m-approximable.
 \item
 Every compact $G$-system is m-approximable.
 \item
 For every $G$-space $X$ and each $f \in RUC(X)$
 the $G$-system $X_f$ is metrizable.
\end{enumerate}
\end{prop}
\begin{proof}
$1 \Rightarrow 4$:\ Given $f \in RUC(X)$ the orbit map $ G \to
RUC(X), \hskip 0.3cm g \mapsto {}_gf$ is uniformly continuous,
where $G$ is endowed with its right uniform structure. Since $G$ is
uniformly Lindel\"of the orbit $fG=\{{}_gf\}_{g \in G}$ is also
uniformly Lindel\"of, hence separable in the Banach space $RUC(X)$
(inspired by \cite[Lemma 10]{usp}). It follows that the Banach
$G$-algebra $\Acal_f$ generated by $fG$ is also separable. By
Proposition \ref{p:generalX_f}.2, $X_f$ is metrizable.

$4 \Rightarrow 2$:\ Consider the $G$-space $X:=G$. Assuming that
each $X_f$ is metrizable, we see, by Remark \ref{r:X_f}.2, that
$G^{ R}=\beta_G X$ is an m-approximable $G$-system.

$2 \Rightarrow 1$:\ Since $G$ naturally embeds as an orbit into
$G^{ R}$, we get that the map $j: G \to \breve{G} \subset
\Homeo(G^R)$ is a homeomorphism. If $G^R$ is m-approximable then
by Proposition \ref{prop:quasisep}, $\breve{G}$ (and hence $G$) is
uniformly Lindel\"of.

$1 \Rightarrow 3$:\ Immediately follows by Proposition \ref
{prop:quasisep}.

$3 \Rightarrow 2$:\ Trivial.
\end{proof}

\br

\section{Almost equicontinuity, local equicontinuity and
variations}\label{Sec-ae}

\sk

By a {\em uniform\/} $G$-space $(X,\mu)$ we mean a $G$-space
$(X,\tau)$ where $\tau$ is a
%aref  I think we really do not need non-Hausdorff spaces
%(at least treating them as a G-space)
%(completely regular not necessarily Hausdorff)
(completely regular Hausdorff) topology, with a {\it compatible}
%aref
%(not necessarily separated)
uniform structure $\mu$, so that the topology $top(\mu)$ defined
by $\mu$ is $\tau$.

\begin{defn} \label{d:un-def}
Let $(X, \mu)$ be a uniform $G$-space.
\begin{enumerate}
\item
A point $x_0\in X$ is a {\em point of equicontinuity\/} (notation:
$x_0 \in Eq(X)$) if for every entourage $\eps\in \mu$, there is a
neighborhood $U$ of $x_0$ such that $(g x_0,g x)\in \eps$ for
every $x\in U$ and $g\in G$. The $G$-space $X$ is {\it
equicontinuous} if $Eq(X)=X$. As usual, $X$ is {\it uniformly
equicontinuous} if for every $\eps\in \mu$ there is $\delta \in
\mu$ such that $(gx,gy) \in \eps$ for every $g \in G$ and $(x,y)
\in \delta$. For compact $X$, equicontinuity and uniform
equicontinuity coincide.
\item
The $G$-space $X$ is {\em almost equicontinuous\/} (AE for short)
if $Eq(X)$ is dense in $X$.
\item
We say that the $G$-space $X$ is {\em hereditarily almost
equicontinuous\/} (HAE for short) if every closed uniform
$G$-subspace of $X$ is AE.
\end{enumerate}
\end{defn}

%eli
The following fact is well known
%23June
at least for metric compact $G$-spaces. See for example
\cite[Proposition 3.4]{AAB2}. Note that neither metrizability nor
compactness of $(X, \mu)$ are needed in the proof.

\begin{lem} \label{l:eq=tr}
If $(X,\mu)$ is a point transitive \footnote{By Lemma \ref{l:NS}.5
one can assume that $X$ is only topologically transitive} uniform
$G$-space and $Eq(X)$ is not empty then $Eq(X)=Trans(X)$.
%6July
%(in particular, $X$ is AE).
%end
\end{lem}
%end

Let $\pi: G \times X \to X$ be a separately continuous (at least)
action on a uniform space $(X,\mu)$. Following \cite[ch. 4]{AG}
define the injective map
$$
\pi_{\sharp}: X \to C(G,X), \hskip 0.3cm \pi_{\sharp}(x)(g)=gx,
$$
where $C(G,X)$ is the collection of continuous  maps from $G$ into
$X$. Given a subgroup $H < G$ endow $C(H,X)$ with the uniform
structure of \emph{uniform convergence} whose basis consists of
sets of the form
$$
\tilde {\ep}=\{(f,f')\in C(H,X): (f(h),f'(h))\in \ep \hskip 0.2cm
\text{for all} \hskip 0.2cm h \in H \}.
$$

We use the map $\pi_{\sharp} : X \to C(H,X)$ to define a uniform
structure $\mu_H$ on $X$, as follows. For $\ep \in \mu$ set
$$
[\eps]_H:=\{(x,y) \in X \times X: (hx,hy) \in \eps \hskip 0.2cm
\text{for all} \hskip 0.2cm h \in H\} \ \ \  (\eps \in \mu).
$$
The collection $\{[\eps]_H: \ep\in \mu\}$ is a basis for $\mu_H$.

Always $\mu \subset \mu_H$ and equality occurs iff the
action of $H$ on $(X, \mu)$ is uniformly equicontinuous. If $(X,
\mu)$ is metrizable and $d$ denotes some compatible metric on $X$,
then the corresponding $\mu_H$ is uniformly equivalent to the
following metric
$$
d_H(x,x')={\sup}_{g\in H}d(g x,g x').
$$

\begin{remarks} \label{r:mu-on-X_f}
\ben
\item
It is easy to characterize $\mu_G$ for
$G$-subsets of $RUC(G)$ (e.g., for $X_f={\cls} (f_{\sharp}(X)) \subset
RUC(G)$), where $\mu$ is the pointwise uniform structure on $RUC(G)$.
The corresponding $\mu_G$ is the metric
uniform structure inherited from the norm of $RUC(G)$.
\item
The arguments of \cite[Theorem 2.6]{AAB1} show
%aref (as I understand the proof in the original work [AAB1]
%the assumptions about
%metrizability of X (and the action of Z_+) are not essential)
that the uniform space $(X, \mu_G)$ is complete
for every compact
%23June
(not necessarily metric)
%end
$G$-system $(X,\mu)$.
 \een
\end{remarks}

\begin{lem}\label{l-equic-p}
Let $(X, \mu)$ be a uniform $G$-space.
The following conditions are equivalent:
\begin{enumerate}
\item
$x_0$ is a point of equicontinuity of the $G$-space $(X, \mu)$.
\item
$x_0$ is a point of continuity of the map $\pi_{\sharp}: X \to
C(G,X)$.
\item
$x_0$ is a point of continuity of the map
${\id}_X: (X, \mu) \to (X,\mu_G)$.
\end{enumerate}
\end{lem}
\begin{proof}
Straightforward.
\end{proof}

\begin{cor} \label{c-equic}
Given a compact system $(G,(X, \mu))$ (with the unique compatible
uniform structure $\mu$) the following conditions are equivalent:
\begin{enumerate}
  \item  $(G,(X,\mu))$ is (uniformly) equicontinuous.
  \item  $\mu_G=\mu$.
  \item  $\pi_{\sharp}: X \to C(G,X)$ is continuous.
  \item  $\mu_G$ is precompact.
\end{enumerate}
\end{cor}
\begin{proof}
 By Remark \ref{r:mu-on-X_f}.2 the uniform space $(X, \mu_G)$ is complete.
 Thus precompact implies compact. This establishes 4 $\Rightarrow$ 1.

%oct17
 The implications
 1 $\Rightarrow$ 2 $\Rightarrow$ 3 $\Rightarrow$ 4 are trivial
 taking into account Lemma \ref{l-equic-p}.
\end{proof}

\begin{lem} \label{l-unif-st}
The uniform structure $\mu_G$ defined above is compatible with subdirect
products. More precisely:
\begin{enumerate}
\item
Let $G$ act on the uniform space $(X, \mu)$ and let $Y$ be a
$G$-invariant subset. Then $(\mu_G)\rest_Y = (\mu\rest_Y)_G$.
\item
Let $\{(X_i, \mu_i): i \in I \}$ be a family of uniform
$G$-spaces. Then $(\prod_i \mu_i)_G = \prod_i (\mu_i)_G$.
\end{enumerate}
\end{lem}
\begin{proof} Straightforward.
\end{proof}

\begin{defn} \label{d:LE=hl}
\ben
\item
Let us say that a subset $K$ of a uniform $G$-space $(X,\mu)$ is
{\it light} if the topologies induced by the uniformities $\mu$
and $\mu_G$ coincide on $K$. We say that $X$ is {\it orbitwise
light} if all orbits are light in $X$.
\item
 $(X,\mu)$ is said to be {\em locally equicontinuous\/} (LE for
short) if every point $x_0 \in X$ is a point of equicontinuity of the
uniform $G$-subspace ${\cls} (Gx_0)$. That is, for every $x_0 \in
X$ and every element $\varepsilon$ of the uniform structure $\mu$
there exists a neighborhood $O$ of $x_0$ in $X$ such that
$(gx,gx_0) \in \varepsilon$ for every $g\in G$ and every $x \in O
\cap {\cls}(Gx_0)$  (see \cite{GW}).
%Since $O \cap \cls (Gx_0) \subset \cls(O \cap Gx_0)$
It is easy to see that the
latter condition, equivalently, can be replaced by the weaker
condition: $x \in O \cap Gx_0$ (this explains Lemma
\ref{l:pr-light}.1 below).
%aref
It follows by Lemma \ref{l:eq=tr} that $X$ is LE iff every
%6July
point
%end
transitive closed $G$-subspace of $X$ is AE.
 \een
\end{defn}

\begin{lem} \label{l:pr-light}
\ben
\item
$x_0 \in X$ is a point of equicontinuity of ${\cls} (Gx_0)$ iff
$Gx_0$ is light in $X$.
\item
$X$ is LE iff $X$ is orbitwise light.
\item
A pointed system $(X,x_0)$ is AE iff the orbit $Gx_0$ is light in
$X$.
\item
Let $f \in RUC(X)$. A subset $K \subset X_f={\cls} (f_{\sharp}(X))$
is light iff the pointwise and norm topologies coincide on $K
\subset RUC(G)$. \een
\end{lem}
\begin{proof}
1. \ Straightforward.

2. \ Follows directly from Claim 1.

3. \ $X$ is point transitive and AE. Therefore the nonempty set
$Eq(X)$ coincides with the set of transitive points (Lemma
\ref{l:eq=tr}). In particular, $x_0 \in Eq(X)$. Thus, $Gx_0$ is
light in $X={\cls} (Gx_0)$ by Claim 1.

Conversely, let $Gx_0$ be a light subset and $x_0$ be a transitive
point. Then again by the first assertion $x_0 \in Eq(X)$. Hence $Eq(X)$
(containing $Gx_0$) is dense in $X$.

4. \ For the last assertion see Remark \ref{r:mu-on-X_f}.1.
\end{proof}

%aref (10)
Given a $G$-space $X$ the collection $AP(X)$ of functions in
$RUC(X)$ coming from equicontinuous systems is the $G$-invariant
uniformly closed algebra of almost periodic functions, where a
function $f \in C(X)$ is {\em almost periodic\/}
%aref  if
iff the set of translates $\{L_g (f):g\in G\}$, where $L_g
(f)(x)=f(g x)$, forms a precompact subset of the Banach space
$C(X)$.
This happens iff $X_f$ is norm compact iff $(G,X_f)$ is an AP
system.
%end

A function $f \in C(X)$ is called {\em weakly almost periodic\/}
(WAP for short, notation: $f \in WAP(X)$) if the set of translates
$\{L_g (f):g\in G\}$ forms a weakly precompact subset of $C(X)$.
We say that a dynamical system $(G,X)$ is {\em weakly almost
periodic\/} if $C(X)=WAP(X)$. The classical theory shows that
$WAP(G)$ is a left and right $G$-invariant, uniformly closed,
point-universal algebra containing $AP(G)$ and that every minimal
function in $WAP(G)$ is in $AP(G)$.
%aref
In fact $f \in WAP(X)$ iff $X_f$ is weakly compact iff $(G,X_f)$
is a WAP system.
%end

The following characterization of WAP dynamical systems is due to Ellis
\cite{E0} (see also Ellis and Nerurkar \cite{EN})
and is based on a result of Grothendieck \cite{Gro}
(namely: pointwise compact bounded subsets in $C(X)$ are weakly
compact for every compact $X$).

\begin{thm}\label{semi}
Let $(G,X)$ be a dynamical system.
The following conditions are equivalent.
\begin{enumerate}
\item
$(G,X)$ is WAP.
\item
The enveloping semigroup
%me
$E(X)$ consists of continuous maps.
\end{enumerate}
\end{thm}

\begin{rmk}
When $(G,X)$ is WAP the enveloping semigroup $E(X)$ is a
semitopological semigroup; i.e. for each $p\in E$ both
$\rho_p: q\mapsto qp$ and $\la_p: q\mapsto pq$ are continuous maps. The
converse holds if in addition we assume that $(G,X)$ is point
transitive. As one can verify the enveloping semigroup of the
dynamical system described in Example \ref{exp} below is
isomorphic to the Bohr compactification of the integers (use
Proposition \ref{envel}). In particular it is a topological group;
however the original system is not even AE and therefore not WAP as we
will shortly see.
\end{rmk}

The next characterization, of AE metric systems, is due to Akin,
Auslander and Berg \cite{AAB2}.

\begin{thm}\label{aees}
Let $(G,X)$ be a compact metrizable system. The following
conditions are equivalent.
\begin{enumerate}
\item
$(G,X)$ is almost equicontinuous.
\item
There exists a dense $G_\del$ subset $X_0\subset X$ such that
every member of the enveloping semigroup $E$ is continuous
on $X_0$.
\end{enumerate}
\end{thm}
Combining these results Akin, Auslander and Berg deduce that every
compact metric WAP system is AE, \cite{AAB2}. Since every
subsystem of a WAP system is WAP it follows from Theorems
\ref{semi} and \ref{aees} that every metrizable WAP system is both
AE and LE. This result is retrieved, and generalized, in \cite{M1}
for all compact $\rm RN_{app}$ $G$-systems using linear
representation methods.

Note that a
%aref
point transitive LE system is of course AE but there are
nontransitive LE systems which are not AE (e.g., see Remark
\ref{exp+}.1 below). It was shown in \cite{GW} that the LE
property is preserved under products under passage to a subsystem
and under factors $X\to Y$ provided that $X$ is metrizable (for
arbitrary systems $X$ see Proposition \ref{LE-inherit} below).

Let $LE(X)$ be the set of functions on a $G$-space $X$ coming from
LE dynamical systems. It then follows from Proposition
\ref{Cor-univ3} that $LE(G)$ is a uniformly closed point-universal
left and right $G$-invariant subalgebra of $RUC(G)$ and that
$LE(X)$, for compact $X$, is the $G$-subalgebra of $C(X)$ that
corresponds to the unique maximal LE factor of $(G,X)$. The
results and methods of \cite{GW} show that $WAP(X) \subset LE(X)$
and that a minimal function in $LE(X)$ is almost periodic (see
also Corollary \ref{cor:LE}.2 below).

\begin{rmk}
In contrast to the well behaved classes of WAP and LE systems, it
is well known that the class of AE systems is closed neither under
passage to subsystems nor under taking factors, see \cite{GW1,
AAB1} and Remark \ref{exp+}.1 below.
\end{rmk}

 By Proposition \ref{Cor-univ3} we see that for every
 $G$-space $X$ the
classes $AP(X), WAP(X)$, $LE(X)$ form $G$-invariant Banach
subalgebras of $RUC(X)$. Recall that for a topological group $G$ we
denote the greatest ambit of $G$ by $G^{ RUC(G)} =G^{R}=|RUC(G)|$.
It is well known that the maximal compactification $u_R: G \to
G^R$ is a right topological semigroup compactification of $G$.
We adopt the following notation. For a $G$-invariant closed
subalgebra $\Acal$ of $RUC(G)$ let $G^{\Acal}$ denote the
corresponding factor $G^{R} \to G^{ \Acal}$ and for a
$G$-space $X$ and a closed $G$-subalgebra $\Acal\subset RUC(X)$, let
$X^{ \Acal}=|\Acal|$ denote the corresponding factor
$\beta_GX \to X^{ \Acal}$.

In the next proposition we sum up some old and new observations
concerning some subalgebras of $RUC(X)$ and $RUC(G)$.

\begin{prop} \label{sp-inclusions}
Let $G$ be a topological group.
\begin{enumerate}
\item
For every $G$-space $X$ we have the following inclusions
$$
RUC(X) \supset LE(X) \supset Asp(X) \supset WAP(X) \supset AP(X),
$$
and the corresponding $G$-factors
$$
\beta_G X  \to X^{ LE} \to X^{Asp} \to X^{WAP} \to X^{ AP}.
$$
\item
For every topological group $G$
we have the following inclusions
$$
RUC(G) \supset UC(G) \supset LE(G) \supset Asp(G)
 \supset WAP(G) \supset AP(G),
$$
and the corresponding $G$-factors
$$
G^{ R}  \to G^{UC} \to G^{ LE} \to G^{Asp} \to G^{ WAP} \to
G^{AP}.
$$
\item
The compactifications $G^{ AP}$ and  $G^{WAP}$ of $G$ are
respectively: a {\em topological group} and a {\em semitopological
semigroup};  $G^R$ and $G^{Asp}$ are {\em right topological
semigroup compactifications\/} of $G$.
\end{enumerate}
\end{prop}
\begin{proof}
For the properties of $Asp(X)$ we refer to section
\ref{s:Asplund}, Theorem \ref{thm1}.6 and Lemma \ref{l-RNis}.2.

In order to show that $UC(G) \supset LE(G)$ we have only to check that
$LUC(G) \supset LE(G)$. Let $f \in LE(G)$. By the definition $f$
%aref (9)
%is coming
comes from a point transitive LE system $(X, x_0)$. Therefore for
some continuous function $F: X \to \R$ we have $f(g)=F(gx_0)$. Let
$\mu$ be the natural uniform structure on $X$. For a given $\ep >0$
choose an entourage $\delta \in \mu$ such that $|F(x)-F(y)| <
\eps$ for every $(x,y) \in \delta$. Since $x_0$ is a point of
equicontinuity we can choose a neighborhood $O$ of $x_0$ such that
$(gx, gx_0) \in \delta$ for every $(g,x) \in G \times O$. Now pick
a neighborhood $U$ of $e \in G$ such that $Ux_0 \subset O$. Then
clearly $|F(gux_0) - F(gx_0)| < \eps$ for every $(g,u) \in G
\times U$; or, equivalently $|f(gu)-f(g)| < \eps$. This means that
$f \in LUC(G)$.
\end{proof}

Now we show the inheritance of LE under factors.

\begin{prop}\label{LE-inherit}
Let  $X$ be a compact LE $G$-system.
If $\pi: X \to Y$ is a $G$-homomorphism then $(G,Y)$ is LE.
\end{prop}

\begin{proof}
We have to show that each point $y_0$ in the space $Y$ is an
equicontinuity point of the subsystem $\OCG(y_0)$. Fix $y_0\in Y$
and assume, with no loss in generality, that $\OCG(y_0)=Y$.
Furthermore, since by Zorn's lemma there is a subsystem of $X$
which is minimal with the property that it projects onto $Y$, we
may and will assume that $X$ itself is minimal with respect to
this property. Denoting by $Y_0$ the subset of transitive points
in $Y$ it then follows that the set $X_0=\pi^{-1}(Y_0)$ coincides
with the set of transitive points in $X$. Let $\ep$ be an element
of the uniform structure of Y (i.e. a neighborhood of the identity
in $Y\times Y$). Then the preimage $\delta:=\pi^{-1}(\ep)$ is an
element of the uniform structure of $X$. Let $q$ be a preimage of
$y_0$. Then $q \in Eq(X)$ since $q$ is transitive and $X$ is LE
(see Lemma \ref{l:eq=tr}). Thus there exists an open neighborhood
$U_q$ of $q$ such that $(gx, gq) \in \delta$ for all $g \in G$ and
$x \in U_q$. Let $V$ be the union of all such $U_q$'s for $q$
running over the preimages of $y_0$. Then $V$ is an open
neighborhood of
%the preimage of $y_0$.
$\pi^{-1}(y_0)$. Set $W$ to be $Y \setminus \pi(X \setminus V)$.
Then $W$ is an open neighborhood of $y_0$
%23June
and $W \subset \pi(V)$.
%because $Y \setminus \pi(V) \subset \pi(X \setminus V) \subset Y$. Hence
%$W=Y \setminus \pi(X \setminus V) \subset \pi(V)$.
%end
For any $y \in W$ we
can find some preimage
$q$ of $y_0$ and some point $x \in U_q$ such
that $\pi(x)=y$. Then $(gx, gq) \in \delta$ for all $g \in G$,
which means that $(gy, gy_0) \in \ep$ for all $g \in G$. Therefore
$y_0 \in Eq(Y)$.
\end{proof}

\begin{cor}\label{cor:LE}
Let $G$ be a topological group, $X$ a $G$-space and $f \in
RUC(X)$. Then
\begin{enumerate}
\item
$f \in LE(X)$ $\Leftrightarrow$ $X_f$ is LE.
\item
If $f \in LE(X)$ is a minimal function then $f \in AP(X)$.
\end{enumerate}
\end{cor}
\begin{proof}
1.\ Use Propositions \ref{LE-inherit} and \ref{Cor-univ3}.3.

2.\
Observe that every minimal LE system is AP.
\end{proof}

Our next result is an intrinsic characterization of the LE
property of a function.

First recall that for the left regular action of $G$ on $X:=G$,
the space $X_f$ can be defined as the pointwise closure of the
orbit $Gf$ (Remark \ref{r:X_f}.1) in $RUC(G)$.

\begin{defn} \label{def:LE}
We say that a function $f \in RUC(G)$ is
\begin{enumerate}
\item
{\em light} (notation: $f \in light(G)$), if the pointwise and
norm topologies coincide on the orbit
$Gf=\{R_g(f)\}_{g \in G}=\{f_g\}_{g \in G}
\subset X_f$ (with $X:=G$) as a subset of $RUC(G)$.
\item
{\em hereditarily light} (notation: $f \in hlight(G)$), if the
pointwise and norm topologies coincide on the orbit $G h$ for
every $h \in X_f$.
\end{enumerate}
By Lemma \ref{l:pr-light}.4 and Definition \ref{d:LE=hl}.1, $f \in
light(G)$ ($f \in hlight(G)$) iff $Gf$ is a {\it light subset}
of the $G$-system $X_f$
(resp.: iff $X_f$ is {\it orbitwise light}).
\end{defn}

\begin{prop} \label{light}
For every topological group $G$ and $f \in RUC(G)$ we have:
\begin{enumerate}
\item
$UC(G) \supset light(G)$.
\item
$f \in light(G)$ $\Leftrightarrow$ $X_f$ is AE.
\item
$f \in hlight(G)$ $\Leftrightarrow X_f$ is LE.
\end{enumerate}
\end{prop}
\begin{proof}
1.\ $f \in light(G)$ means that the pointwise and norm topologies
coincide on $Gf$. It follows that the orbit map $G \to RUC(G),
\hskip 0.2cm g \mapsto f_g$ is norm continuous. This means that
$f$ is also left uniformly continuous.

2. \ Since $f$ is a transitive point of $X_f={\cls} (Gf)$ we can use
Lemma \ref{l:pr-light}.3.

3. \ Use Lemma \ref{l:pr-light}.2.
\end{proof}

\begin{thm} \label{c:Hlight=LE}
   $LE(G) = hlight(G)$ for every topological group $G$.
\end{thm}
\begin{proof}
 Follows from Proposition \ref{light}.3 and Corollary \ref{cor:LE}.1.
\end{proof}

%23June
%Q: let pointwise and norm topologies coincide on $Gf$
%(that is $f \in light(G)$)
%Is it true that the pointwise and norm topologies coincide on $fG$  ??
%Probably one more question which suits to "second" paper

\begin{remarks}
\begin{enumerate}
\item
%23June
By \cite[Theorem 8.5]{M1}, for every topological group $G$ and
every $f \in WAP(G)$ the pointwise and norm topologies coincide on
$fG=\{L_g(f)\}_{g \in G}=\{{}_gf \}_{g \in G}$. Using the involution
$$
UC(G) \to UC(G), \hskip 0.3cm f \mapsto f^* \hskip 0.6cm
(f^*(g):=f(g^{-1}))
$$
(observe that $Gf^*=(fG)^*$) we get the coincidence of the above
mentioned topologies also on $Gf^*$. Since $(WAP(G))^*=WAP(G)$ we
can conclude that
 $WAP(G) \subset light(G)$ for every
topological group $G$. Theorem \ref{c:Hlight=LE} provides a
stronger inclusion $ LE(G) \subset light(G)$ (since $WAP(G)
\subset LE(G)$ by Proposition \ref{sp-inclusions}.2).
%end
\item
In view of Proposition \ref{light}.2 a minimal function is light
iff it is AP.
%oct17 (we can assume wrg that the function $f$ comes
%from $X_f$ and moreover, that $X_f$ is minimal) RIGHT ??
Thus, for example, the function $f(n)=\cos(n^2)$ on the integers,
which comes from a minimal distal but not equicontinuous
$\Z$-system on the two torus, is not light.
\end{enumerate}
\end{remarks}

\br

\section{Fragmented maps and families} \label{sec:fr}

\sk

The following definition is a generalized version of
{\it fragmentability} (implicitly
it appears in a paper of Namioka and Phelps \cite{NP}) in
the sense of Jayne and Rogers \cite{JR}.

\begin{defn} \label{def:fr}
\cite{Mefr} Let $(X,\tau)$ be a topological space and $(Y,\mu)$ a
uniform space.
\begin{enumerate}
\item
We say that $X$ is {\em $(\tau, \mu)$-fragmented\/} by a (not
necessarily continuous) function $f: X \to Y$ if for every
nonempty subset $A$ of $X$ and every $\eps \in \mu$ there exists
an open subset $O$ of $X$ such that $O \cap A$ is nonempty and the
set
%23June
%oscillation $osc(f, O \cap A)$ ?
%end
$f(O \cap A)$ is $\eps$-small in $Y$. We also say in that case
that the function $f$ is {\em fragmented\/}. Note that it is
enough to check the condition above only for closed subsets $A
\subset X$ and
%use the fact that if $O_2 \cap \cls(O_1 \cap A)$ is nonempty
%then $(O_2 \cap O_1) \cap A$ is also nonempty ...
for $\ep \in \mu$ from a {\it subbase} $\gamma$ of $\mu$ (that is,
the finite intersections of the elements of $\gamma$ form a base
of the uniform structure $\mu$).
\item
If the condition holds only for every nonempty {\em open subset\/}
$A$ of $X$ then we say that $f$ is {\em locally fragmented\/}.
%7July
%\item
%end
\item
%gl+++
When the inclusion map $i:X \subset Y$ is (locally) fragmented
we say that $X$ is
{\em (locally) $(\tau, \mu)$-fragmented\/}, or more simply,
{\em (locally) $\mu$-fragmented}.
%If $X \subset Y$ and $f$ is the inclusion map then we say that $X$
%is (locally) $(\tau, \mu)$-fragmented, or even simpler, (resp.:
%locally) $\mu$-fragmented.
\end{enumerate}
\end{defn}

\begin{rmk} \label{r:fr}
\ben
\item
Note that in Definition \ref{def:fr}.1 when $Y=X, f={\id}_X$ and
$\mu$ is a metric uniform structure, we get the usual definition
of fragmentability \cite{JR}. For the case of functions see also
\cite{JOPV}.
\item
Namioka's {\it joint continuity theorem} \cite{N-jct} (see also
Theorem \ref{t:N-jct} below) implies that every weakly compact
subset $K$ of a Banach space is (weak,norm)-fragmented (that is,
$id_K: (K,weak) \to (K,norm)$ is fragmented).
\item
Recall that a Banach space $V$ is an {\em Asplund\/} space if the
dual of every separable Banach subspace is separable, iff every
bounded subset $A$ of the dual $V^*$ is
(weak${}^*$,norm)-fragmented, iff $V^*$ has the Radon-Nikod\'ym
property. Reflexive spaces and spaces of the type $c_0(\Gamma)$
are Asplund. For more details cf. \cite{Bo, F, N}.
\item
A topological space $(X,\tau)$ is {\em scattered\/} (i.e., every
nonempty subspace has an isolated point) iff $X$ is $(\tau,
\rho)$-fragmented, where $\rho(x,y)=1$ iff $x \neq y$.
 \een
\end{rmk}

Following \cite{MN} we say that $f: X \to Y$ is {\em barely
continuous\/} if for every nonempty closed subset $A \subset X$,
the restricted map $f\rest_A$ has at least one point of
continuity.

\begin{lem} \label{simple-fr}
\begin{enumerate}
\item
If $f$ is $(\tau,\mu)$-continuous then $X$ is
$(\tau,\mu)$-fragmented by $f$.
\item
Suppose that there exists a dense subset of
$(\tau,\mu)$-continuity points of $f$. Then $X$ is locally
$(\tau,\mu)$-fragmented by $f$.
\item
$X$ is $(\tau,\mu)$-fragmented by $f$ iff $X$ is {\em hereditarily
locally fragmented\/} by $f$ (that is, for every closed subset $A
\subset X$ the restricted function $f\rest_A$ is (relatively)
locally $(\tau,\mu)$-fragmented).
\item
Every barely continuous $f$ is fragmented.
\item
Fragmentability is preserved under products. More precisely, if
$f_i: (X_i, \tau) \to (Y_i, \mu_i)$ is fragmented for every $i \in
I$ then the product map
$$
f:=\prod_{i \in I} f_i: \prod_{i \in I} X_i \to \prod_{i \in I}
Y_i
$$
is $(\tau, \mu)$-fragmented with respect to the product topology
$\tau$ and the product uniform structure $\mu$.
%The same is true also for diagonal products.
%oct17
\item
Let $\al: X \to Y$ be a continuous map.
If $f: Y \to (Z,\mu)$ is a fragmented map then
the composition $f \circ \al: X \to (Z,\mu)$ is also fragmented.
\end{enumerate}
\end{lem}
\begin{proof} The assertions 1, 2 and 6 are straightforward.
%(6): Let $A$ be a nonempty subset of $X$ and let $\eps \in \mu$.
%Choose an open subset $O$ in $Y$ such that $\a(A) \cap O$ is
%nonempty and $f(\a(A) \cap O)$ is $\eps$-small.
%Since $\al(A \cap \al^{-1}(O))=\al(A) \cap O$ we get that
%$A \cap \al^{-1}(O)$ is nonempty and
%$(f \circ \a) (A \cap \al^{-1}(O))$ is $\eps$-small in $Y$. This completes
%the proof because $\al^{-1}(O)$ is open in $X$ by the continuity of $\a$.

For 3 and 4 use the fact that it is enough to check the
fragmentability condition only for closed subsets $A \subset X$.

The verification of 5 is straightforward taking into account that
 it is enough to check the
fragmentability (see Definition \ref{def:fr}.1) for $\eps \in
\gamma$, where $\gamma$ is a subbase of $\mu$.
\end{proof}

Fragmentability has good stability properties being closed under
passage to subspaces (trivial), products (Lemma \ref{simple-fr}.5)
and quotients. Here we include the details for quotients. The
following lemma is a generalized version of \cite[Lemma 4.8]{Mefr}
which in turn was inspired by Lemma 2.1 of Namioka's paper
\cite{N}.

\begin{lem} \label{l-quot-fr}
Let $(X_1,\tau_1)$ and $(X_2,\tau_2)$ be compact (Hausdorff)
spaces, and let $(Y_1, \mu_1)$ and $(Y_2, \mu_2)$ be uniform
spaces. Suppose that: $F: X_1\to X_2$ is a continuous surjection,
$f: (Y_1, \mu_1) \to (Y_2, \mu_2)$ is uniformly continuous,
$\phi_1: X_1 \to Y_1$ and $\phi_2: X_2 \to Y_2$ are maps such that
the following diagram

\begin{equation*}
\xymatrix {
(X_1, \tau_1) \ar[d]_{F} \ar[r]^{\phi_1} & (Y_1, \mu_1) \ar[d]^{f} \\
(X_2, \tau_2) \ar[r]^{\phi_2} & (Y_2, \mu_2) }
\end{equation*}
commutes. If $X_1$ is fragmented by $\phi_1$ then
$X_2$ is fragmented by $\phi_2$.
\end{lem}
\begin{proof}
We modify the proof of \cite[Lemma
4.8]{Mefr}. In the definition of fragmentability it suffices to
check the condition for {\em closed subsets\/}. So, let $\ep \in
\mu_2$ and let $A$ be a non-empty closed and hence, compact subset
of $X_2$.
 Choose $\delta \in\mu_1$ such that $(f\times
f)\ (\delta)\subset \ep.$ By Zorn's Lemma, there exists a minimal
compact subset $M$ of $X_1$ such that $F(M) = A.$ Since $X_1$ is
fragmented by $\phi_1$, there exists $V \in \tau_1$ such that
$V\cap M\neq \emptyset$ and $\phi_1(V\cap M)$ is $\delta$-small.
Then the set $f \phi_1(V\cap M)$ is $\ep$-small. Consider the set
$W:=A \setminus F(M \setminus V).$ Then
\begin{itemize}
\item [(a)] $\phi_2(W)$ is $\ep$-small,
being a subset of $f \phi_1(V\cap M)= \phi_2 F( V \cap M)$.
%$A \setminus F(V \cap M) \subset F(M \setminus V) \subset A$. Hence
%$W=A \setminus F(M \setminus V) \subset F(V \cap M)$
\item [(b)] $W$ is relatively open in $A$;
\item [(c)] $W$ is non-empty
\newline (otherwise, $M \setminus V$ is a {\em {proper}\/} compact
subset of $M$ such that $F(M \setminus V) = A$).
\end{itemize}
\end{proof}

The next lemma provides a key to understanding the connection
between fragmentability and separability properties.

\begin{lem} \label{fr-sep-0}
Let $(X,\tau)$ be a separable metrizable space and $(Y,\rho)$ a
pseudometric space. Suppose that $X$ is $(\tau, \rho)$-fragmented
by a surjective map $f: X \to Y$. Then $Y$ is separable.
\end{lem}

\begin{proof}
Assume (to the contrary) that the pseudometric space $(Y,\rho)$ is not
separable. Then there exist an $\eps >0$ and an uncountable subset
$H$ of $Y$ such that
$\rho (h_1,h_2) > \eps$ for all distinct $h_1,h_2 \in H$. Choose a
subset $A$ of $X$ such that $f(A)=H$ and $f$ is bijective on $A$.
Since $X$ is second countable the uncountable subspace $A$ of $X$
(in its relative topology) is a disjoint union of a countable set
%of isolated points
and a nonempty closed perfect set $M$ comprising the condensation
points of $A$ (this follows from the proof of the Cantor-Bendixon
theorem; see e.g. \cite{Kech}). By fragmentability there exists an
open subset $O$ of $X$ such that $O \cap M$ is nonempty and $f(O
\cap M)$ is $\eps$-small. By the property of $H$ the intersection
$O \cap M$ must be a singleton, contradicting the fact that no
point of $M$ is isolated.
\end{proof}

\begin{prop} \label{fr-sep-1}
If $X$ is locally fragmented by $f:X \to Y$, where $(X,\tau)$ is a
Baire space and $(Y,\rho)$ is a pseudometric space then $f$ is
continuous at the points of a dense $G_{\delta}$ subset of $X$.
\end{prop}
\begin{proof}
For a fixed $\ep > 0$ consider
$$
O_{\ep}:= \{ \text{union of all $\tau$-open subsets $O$ of $X$
with} \hskip 0.3cm {\diam}_{\rho} f(O) \leq \ep \}.
$$
The local fragmentability implies that $O_{\ep}$ is dense in $X$.
Clearly, $\bigcap \{O_{\frac{1}{n}} :   n \in \N \}$ is the
required dense $G_{\delta}$ subset of $X$.
\end{proof}

A topological space $X$ is {\it hereditarily Baire} if every
closed subspace of $X$ is a Baire space. Recall that for
metrizable spaces $X$ and $Y$ a function $f: X \to Y$ is of {\em
Baire class 1\/} if $f^{-1}(U)\subset X$ is an $F_\sig$ subset for
every open $U\subset Y$.
%23June
If $X$ is separable then a real valued function $f: X \to \R$ is
of Baire class 1 iff $f$ is the pointwise limit of a sequence of
continuous functions
%end
(see e.g. \cite{Kech}).

\begin{prop} \label{Baire-1}
Let $(X,\tau)$ be a hereditarily Baire (e.g., Polish, or compact)
space, $(Y,\rho)$ a pseudometric space. Consider the following
assertions:
\begin{enumerate}
\item [(a)]
$X$ is $(\tau,\rho)$-fragmented by $f: X \to Y$;
\item [(b)]
$f$ is barely continuous;
\item [(c)]
$f$ is of Baire class 1.
\end{enumerate}
\ben
\item
Then $(a) \Leftrightarrow (b)$.
%If $X$ is Polish then (a)
%$\Leftrightarrow$ (b) $\Rightarrow$ (c).
\item
If $X$ is Polish and $Y$ is a separable metric space then $(a)
\Leftrightarrow (b) \Leftrightarrow (c)$.
\een
\end{prop}
\begin{proof} For (a) $\Leftrightarrow$ (b) combine
Lemma \ref{simple-fr} and Proposition \ref{fr-sep-1}.

The equivalence
(b) $\Leftrightarrow$ (c) for Polish $X$ and separable $Y$ is well
known (see \cite[Theorem 24.15]{Kech}) and actually goes back to
Baire.
%Observe that if $f$ is barely continuous then $f$ is
%fragmented (Lemma \ref{simple-fr}.4). Hence by Lemma
%\ref{fr-sep-0} we can assume in the proof of the implication (b) $\Rightarrow$ (c)
%that $f(X)$ is separable.
\end{proof}

\sk

The following new definition will play a crucial role in Section
\ref{envelop}.

\begin{defn} \label{d:fr-family}
\ben
\item
We say that a {\it family of functions} $\mathcal{F}=\{f: (X,\tau)
\to (Y,\mu) \}$ is {\it fragmented} if the condition of Definition
\ref{def:fr}.1 holds simultaneously for all $f \in \mathcal F$.
That is, $f(O \cap A)$ is $\eps$-small for every $f \in \mathcal
F$. It is equivalent to say that the mapping
$$
\pi_{\sharp}: X \to Y^{\mathcal{F}}, \hskip 0.4cm
\pi_{\sharp}(x)(f)=f(x)
$$
is $(\tau, \mu_U)$-fragmented, where $\mu_U$ is the uniform structure of
uniform convergence on the set $Y^{\mathcal{F}}$ of all mappings
from $\mathcal{F}$ into $(Y, \mu)$.
\item
Analogously one can define the notions of {\it a locally
fragmented family} and {\it a barely continuous family}. The
latter means that every closed nonempty subset $A\subset X$ contains a
point $a\in A$ such that ${\mathcal F}_A =\{f
\upharpoonright_A: f\in {\mathcal F} \}$ is equicontinuous at $a$.
%me!
If $\mu$ is pseudometrizable then so is $\mu_U$. Therefore if in
addition $(X, \tau)$ is hereditarily Baire then it follows by
Proposition \ref{Baire-1}.1 that $\mathcal{F}$ is {\it fragmented}
iff $\mathcal{F}$ is barely continuous. \een
\end{defn}

Fragmented families, like equicontinuous families, are stable
under pointwise closures as the following lemma shows.

\begin{lem} \label{l:fr-f-cls}
Let $\mathcal{F}=\{f: (X,\tau) \to (Y,\mu) \}$ be a fragmented
family of functions. Then the pointwise closure
$\bar{\mathcal{F}}$ of $\mathcal{F}$ in $Y^X$ is also a $(\tau,
\mu)$-fragmented family.
\end{lem}
\begin{proof}
Use straightforward ``$3\eps$-trick" argument.
\end{proof}

\br

\section{Asplund functions and RN systems} \label{s:Asplund}

\sk

Let $H$ be a subgroup of $G$. Recall that we denote by $\mu_H$ the
uniform structure on the uniform $G$-space $(X, \mu)$ inherited by
the inclusion $\pi_{\sharp}: X \to C(H,X)$. Precisely, $\mu_H$ is
generated by the basis $\{[\eps]_H: \eps \in \mu \}$, where
$$
[\eps]_H:=\{(x,y) \in X \times X: (hx,hy) \in \eps \hskip 0.2cm
\text{for all} \hskip 0.2cm h \in H\}.
$$
For every $f \in C(X)$ and $H < G$ denote by $\rho_{H,f}$ the
pseudometric on $X$ defined by
$$
\rho_{H,f}(x,y)=\sup_{h \in H} |f(hx)-f(hy)|
$$
Then
$\mu_{{\cls}(H)} = \mu_H$ and $\rho_{{\cls}(H),f} = \rho_{H,f}$.

\begin{defn}\label{d:asp}
\begin{enumerate}
\item
A continuous function $f: X \to \R$ on the compact $G$-space $X$
is an {\em Asplund function\/} \cite{M1}
if
%fG is an Asplund set
for every countable subgroup $H \subset G$ the
pseudometric space $(X,\rho_{H,f})$ is separable. It is an {\em
s-Asplund function\/} (notation: $f \in Asp_s(X)$) when
$(X,\rho_{G,f})$ is separable. A pseudometric $d$ on a set $X$ is
called {\em Asplund\/} (respectively, {\em s-Asplund\/}) if for
every countable subgroup $H < G$ (respectively, for $H=G$) the
pseudometric space $(X,d_{H})$ is separable, where
$$
d_{H}(x, y)= {\sup}_{h \in H}  d(hx , hy).
$$
\item
More generally, we say that a function $f \in RUC(X)$ on a (not
necessarily compact) $G$-space $X$ is an {\it Asplund function}
(notation: $f \in Asp(X)$) if $f$ is coming (in the sense of
Definition \ref{d:coming}) from an Asplund function $F$ on a
$G$-system $Y$ and a $G$-compactification $\nu: X \to Y$. By
Remark \ref{r:asp-def}.2 below, equivalently, one can take each of
the following $G$-compactifications (see Remark \ref{r:X_f}.3)
$f_{\sharp}: X \to X_f$ (minimally possible) or $i_{\beta}: X \to
\beta_G X$ (maximal).
%oct17
Analogously we define
the class $Asp_s(X)$ of s-Asplund functions on a $G$-space $X$.
\item
In particular, a function $f\in RUC(G)$ is an {\em Asplund
function\/} ({\em s-Asplund function\/}) if it is Asplund
(s-Asplund) for the  $G$-space $X:=G$
%aref
with respect to the regular left action.
Notation: $f \in Asp(G)$ (resp.: $f \in Asp_s(G)$).
\end{enumerate}
\end{defn}

\begin{remarks} \label{r:asp-def}
\ben
\item
 Note that in the definition of Asplund functions $F: X \to
\R$, equivalently, we can run over all uniformly Lindel\"of
subgroups $H < G$. Indeed, as in the proof of Proposition
\ref{prop:blocks}, the orbit
%23June
$FH=\{{}_hF\}_{h \in H}$
%$HF$
is norm separable. Let $K < G$
be a countable subgroup of $H$ such that $FK$ is dense in $FH$.
Then $\rho_{H,F} = \rho_{K,F}$.
\item
Let $q: Y_1 \to Y_2$ be a $G$-homomorphism of compact $G$-spaces.
It is straightforward to show that a continuous bounded function
$F: Y_2 \to \R$ is Asplund
(s-Asplund) iff the function $f=F \circ q: Y_1 \to \R$ is Asplund
(resp.: s-Asplund).
\item
Of course every s-Asplund function is Asplund. If $G$, or, the
natural restriction $\breve{G}$, is uniformly Lindel\"of (e.g.
$\breve{G}$ is second countable if $X$ is compact and metrizable)
then clearly the converse is also true. Thus in this case
$Asp(X)=Asp_s(X)$.
\item
Let $(G,X)$ be a dynamical system and $d$ a pseudometric on $X$.
Suppose $F: X\to \R$ is $d$-uniformly continuous. If $d$ is
Asplund or s-Asplund then so is $F$.
\een
\end{remarks}

Let $X$ be a $G$-space. By Proposition \ref{p:generalX_f}.1,
$X_f:=\cls f_{\sharp}(X)$ is a subset of $RUC(G)$ for every $f \in
RUC(X)$. Let $r_G: X_f \hookrightarrow RUC(G)$ be the inclusion
map. For every subgroup $H < G$ we can define the natural
restriction operator $q_H: RUC(G) \to RUC(H)$.
%23June
Denote by $r_H:= q_H \circ r_G: X_f \to RUC(H)$ the composition
and let $\xi_{H,f}$ be the corresponding pseudometric induced on
$X_f$ by the norm of $RUC(H)$. Precisely,
$$
\xi_{H,f}(\om,\om')= \sup_{h\in H} |\om(h)-\om'(h)|.
$$
Finally define the composition $f^H_{\sharp}:= r_H \circ
f_{\sharp}: X \to RUC(H)$. The corresponding pseudometric induced
by $f^H_{\sharp}$ on $X$ is just $\rho_{H,f}$.
%end

\begin{lem} \label{l:asp-def}
Let $X$ be a $G$-space and $f \in RUC(X)$.
%oct17
Let $F_e: X_f \to \R$ be the map
%oct21
$F_e(\omega)=\omega(e)$
(defined before
Proposition \ref{p:generalX_f}). The following are
equivalent:
 \ben
\item $f \in Asp(X)$.
\item $F_e \in Asp(X_f)$.
\item $(X_f, \xi_{H,f})$ is separable for every
countable (uniformly Lindel\"of) subgroup $H < G$.
\item $r_H(X_f)$ is norm separable in $RUC(H)$ for every
countable (uniformly Lindel\"of) subgroup $H < G$.
 \een
\end{lem}
\begin{proof}
1 $\Leftrightarrow$ 2 follows by Definition \ref{d:asp}.2,
Remark \ref{r:asp-def}.2 and Proposition \ref{p:generalX_f}.3.

3 $\Leftrightarrow$ 4 is clear by the definitions of $\xi_{H,f}$
and $r_H$.

2 $\Leftrightarrow$ 3: $F_e \in Asp(X_f)$ means, by Definition
\ref{d:asp}.1, that for every countable (uniformly Lindel\"of)
subgroup $H < G$ the pseudometric space $(X_f,\rho_{H,F_e})$ is
separable, where
$$
\rho_{H,F_e}(\om, \om')=\sup_{h\in H}|F_e(h \om) - F_e(h \om')|.
$$

 Recall that by the definition
of $F_e: X_f \to \R$ we have $F_e(h \om)=(h\om)(e)=\om(h)$. Hence

$$
\xi_{H,f}(\om,\om')=\sup_{h\in H}|\om(h) - \om'(h)|=\sup_{h\in
H}|F_e(h \om) - F_e(h \om')|=\rho_{H,F_e}(\om, \om')
$$

Therefore we obtain that the pseudometrics $\xi_{H,f}$ and
$\rho_{H,F_e}$ coincide on $X_f$. This clearly completes the
proof.
\end{proof}

\begin{cor} \label{c:s-asp}
Let $X$ be a $G$-space and $f \in RUC(X)$. The following are
equivalent:
%oct17
 \ben
\item $f \in Asp_s(X)$.
\item $F_e \in Asp_s(X_f)$.
\item $X_f$ is norm separable in $RUC(G)$.
\een
\end{cor}
\begin{proof}
The proof of Lemma \ref{l:asp-def} shows that in fact $\xi_{H,f}$
and $\rho_{H,F_e}$ coincide on $X_f$ for every $H<G$. Consider the
particular case of $H:=G$ taking into account that $r_G(X_f)=X_f$.
%and use Remark \ref{r:asp-def}.2.
\end{proof}

The following definition of RN dynamical systems (a natural
generalization of RN compacta in the sense of Namioka \cite{N})
and {\it Eberlein systems} (a natural generalization of Eberlein
compacta in the sense of Amir-Lindenstrauss \cite{AL}) were
introduced in \cite{M1}. About the definition and properties of
Asplund spaces see Remark \ref{r:fr}.3 and \cite{Bo,N,F}.

\begin{defn} \label{d-RN}
Let $(G,X)$ be a compact dynamical system.
\begin{enumerate}
%aref   "representation" should be defined in general
%(not necessarily embeddings) - Eyal's correction (March, 2005)
%$(G,X)$ is a {\em Radon-Nikod\'ym system\/} (RN for short) if
%there exists a representation $(h,\al): G\times X\to
%{\Iso}(V)\times V^*$, where $V$ is an Asplund Banach space, $h: G
%\to {\Iso}(V)$ is a strongly continuous co-homomorphism and $\alpha
%: X \to (V^*, w^*)$ is a bounded weak$^*$ $G$-embedding (with respect to
%the {\em dual action\/} $G \times V^* \to V^*,
%(g\varphi)(v):=\varphi(h(g)(v))$). If we can choose $V$ to be
%reflexive, then $(G,X)$ is called an {\em Eberlein system}.
%Notation: RN and Eb respectively.
\item A continuous (\emph{proper}) \emph{representation}
of $(G,X)$ on a Banach space $V$ is a pair $(h,\al)$,
%23June
%$(h,\al): G\times X\to {\Iso}(V)\times V^*$,
%end
where $h: G \to Iso(V)$ is a strongly
continuous co-homomorphism of topological groups and $\a: X \to
V^*$ is a weak star continuous bounded $G$-mapping (resp.
\emph{embedding}) (with respect to the {\em dual action\/} $G
\times V^* \to V^*, (g\varphi)(v):=\varphi(h(g)(v))$).
\item
$(G,X)$ is a {\em Radon-Nikod\'ym system\/} (RN for short) if
there exists a proper representation
%23June
%$(h,\al): G\times X\to {\Iso}(V)\times V^*$, where $V$ is
of $(G,X)$ on an Asplund Banach space $V$. If we can choose $V$ to
be reflexive, then $(G,X)$ is called an {\em Eberlein system}. The
classes of Radon-Nikod\'ym and Eberlein compact systems will be
defined by RN and Eb respectively.
\item
$(G,X)$ is called an RN-{\em approximable\/} system
($\rm{RN_{app}}$) if it can be represented as a subdirect product
%gl+++
(or equivalently, as an inverse limit) of RN systems.
%(that is, $(G,X)$ is isomorphic to a subsystem of a
%product of RN systems).
\end{enumerate}
\end{defn}

Note that compact spaces which are not RN necessarily are
non-metrizable, while there are many natural {\it metric} compact
$G$-systems which are not RN.

The next theorem collects some useful properties which were
obtained recently in \cite{M1}.

\begin{thm} \label{thm1}
Let $(G,X)$ be a compact $G$-system.
\begin{enumerate}
\item
$X$ is WAP iff $X$ is a subdirect product of Eberlein $G$-systems.
A metric system $X$ is WAP iff $X$ is Eberlein.
\item
The system $(G,X)$ is RN iff there exists a representation $(h,
\alpha)$ of $(G,X)$ on a Banach space $V$ such that: $h: G \to
{\Iso}(V)$ is a co-homomorphism (no continuity assumptions on $h$),
$\alpha : X \to V^*$ is a bounded weak$^*$ $G$-embedding and
$\alpha(X)$ is (weak$^*$, norm)-fragmented.
\item
$f: X \to \R$ is an Asplund function iff $f$ arises from an
\emph{Asplund representation} (that is, there exists a continuous
representation $(h,\al)$ of $(G,X)$ on an Asplund space $V$, such
that $f(x)= \al(x)(v)$ for some $v \in V$). Equivalently, iff $f$
comes from an RN (or, $\rm{RN_{app}}$) $G$-factor $Y$ of $X$.
\item
The system $(G,X)$ is $RN_{app}$ iff $Asp(X)=C(X)$.
\item
RN is closed under countable products and
$\rm{RN_{app}}$ is closed under quotients.
 For metric compact systems $\rm{RN_{app}=RN}$ holds.
\item
$Asp(X)$ is a closed $G$-invariant subalgebra of $C(X)$ containing
$WAP(X)$.
%oct17
The canonical compactification $u_A: G \to G^{Asp}$ is
the universal RN$_{app}$ compactification of $G$.
Moreover, $u_A$ is a right topological semigroup
compactification of $G$.
\item
$(G,X)$ is RN iff $(G,(C(X)^*_1,w^*))$ is RN iff $(G,P(X)$) is RN,
where $P(X)$ denotes the space of all probability measures on $X$
(with the induced action of $G$).
\end{enumerate}
\end{thm}

The proofs of the assertions 1, 2 and 3 use several ideas from
Banach space theory; mainly the notion of {\em Asplund sets\/} and
Stegal's generalization of a factorization construction by Davis,
Figiel, Johnson and Pe{\l}czy\'nski \cite{D,Bo,N,St,F}.

\begin{prop} \label{r:ASP}
Let $G$ be an arbitrary topological group. Then $(G^{Asp},u_A(e))$
is point-universal (hence $X_f \subset Asp(G)$ for every $f \in
Asp(G)$).
\end{prop}
\begin{proof}
%1.\
$\Pcal:=Asp(G)$ is an algebra of functions coming from RN$_{app}$
systems. Since the class RN$_{app}$ is preserved by products and
subsystems we can apply Proposition \ref{Cor-univ3}.2.
\end{proof}

Let $(X, \tau)$ be a topological space. As usual, a metric $\rho$
on the set $X$ is said to be {\em lower semi-continuous\/} if the
set $\{ (x,y) : \rho(x,y) \leq t \} $ is closed in $X \times X$
for each $t > 0$. A typical example is any subset $X \subset V^*$
of a dual Banach space equipped with the weak${}^*$ topology and
the norm metric. It turns out that every lower semi-continuous
metric on a compact Hausdorff space $X$ arises in this way (Lemma
\ref{JNR}.1). This important result has been established in
\cite{JNR} using ideas of Ghoussoub and Maurey.

\begin{lem} \label{JNR}
\begin{enumerate}
\item \cite{JNR}
Let $(X,\tau)$ be a compact space and let $\rho \leq 1$ be a lower
semi-continuous metric on $(X,\tau)$. Then there is a dual Banach
space $V^*$ and a homeomorphic embedding  $\alpha: (X,\tau) \to
(V_1^*, w^*)$ such that
$$
|| \a (x)-\a(y)||=\rho(x,y)
$$
for all $x,y \in X$.
\item
If in addition $X$ is a $G$-space and $\rho$ is $G$-invariant,
then the Claim 1 admits a $G$-generalization. More precisely,
there is a linear isometric (not necessarily jointly continuous)
right action $V \times G \to V$ such that $\a: X \to V_1^*$ is a
$G$-map.
\end{enumerate}
\end{lem}
\begin{proof} 2:
As in the proof of \cite[Theorem 2.1]{JNR} the required Banach
space $V$ is defined as the space of all continuous real-valued
functions $f$ on $(X,\tau)$ which satisfy a uniform Lipschitz
condition of order $1$ with respect to $\rho$, endowed with the
norm
$$
p(f)={\max}\{||f||_{Lip}, ||f|| \},
$$
where $||f||={\sup}\{|f(x)|: x\in X\}$ and the seminorm
$||f||_{Lip}$ is defined to be the least constant $K$ such that
$|f(x_1)-f(x_2)| \leq K \rho(x_1,x_2)$ for all $x_1,x_2 \in X$.
Then $\alpha: (X,\tau) \to (V_1^*, w^*)$ is defined by
$\alpha(x)(f)=f(x)$.

%23June
Define now the natural right action $\pi: V \times G \to V$ by
$\pi(f,g)=fg={}_gf$, where ${}_gf(x):=f(gx)$. Then clearly
$p(fg)=p(f)$ and $\a: X \to V_1^*$ is a $G$-map.
%end
\end{proof}

\begin{thm} \label{frag=rn}
Let $(G,X)$ be a compact dynamical system. The following
conditions are equivalent:
\begin{enumerate}
\item
$(G,X)$ is RN.
\item
$X$ is fragmented with respect to some bounded lower
semi-continuous $G$-invariant metric $ \rho$.
\end{enumerate}
\end{thm}

\begin{proof}
1 $\Rightarrow$ 2: Our $G$-system $X$, being RN, is a $G$-subsystem
of the ball $V^*_1=(V^*_1,w^*)$ for some Asplund space $V$. By a
well known characterization of Asplund spaces, $V^*_1$ is $(w^*,
norm)$-fragmented. Hence, $X$ is also fragmented by the lower
semi-continuous $G$-invariant metric on $X$,
$\rho(x_1,x_2)=||x_1-x_2||$ inherited by the norm of $V^*$.

2 $\Rightarrow$ 1 :
%23June
We can suppose that $\rho \leq 1$.
%end
Using Lemma \ref{JNR}.1 we can find a Banach space $V$ and a
weak${}^*$ embedding $\alpha: (X,\tau) \to V^*_1$ such that
$\alpha$ is $(\rho, norm)$-isometric. Since $X$ is
$(\tau,\rho)$-fragmented we get that $\alpha (X) \subset V^*_1$ is
$(w^*,norm)$-fragmented. Moreover, by Claim 2 of the same lemma,
there exists a co-homomorphism (without continuity assumptions)
$h: G \to \Iso (V)$ (the right action $V \times G \to V$ leads to
the co-homomorphism $h$) such that the map $\alpha: X \to V^*_1$
is $G$-equivariant with respect to the dual action of $G$ on $V^*$
defined by $(g\varphi)(v):=\varphi(h(g)(v))$. Therefore we get a
representation $(h,\alpha)$ of $(G,X)$ on $V$ such that $\alpha
(X) \subset V^*_1$ is $(w^*,norm)$-fragmented. By Theorem
\ref{thm1}.2 we deduce that the $G$-system $(X, \tau)$ is RN.
\end{proof}

\br

\section{Veech functions} \label{s:Veech}

\sk

The algebra $K(G)$ was defined
by Veech in \cite{V} --- for a discrete group $G$
--- as the algebra of functions
$f\in \ell^\infty(G)$ such that for every countable subgroup $H <
G$ the collection $X_{f \rest_H}=\OC_H(\eta_0)\subset \Om_H=
[-\|f\|,\|f\|]^H$, with $\eta_0=f \rest_H$, considered as a
subspace of the Banach space $\ell^\infty(H)$, is norm separable.
Replacing $\ell^\infty(G)$ and $\ell^\infty(H)$ by $RUC(G)$ and
$RUC(H)$, respectively, we define --- for any topological group
$G$ --- the algebra $K(G)\subset RUC(G)$ as follows.

\begin{defn}
Let $G$ be a topological group. We say that a function $f\in
RUC(G)$ is a {\em Veech function\/} if for every countable
(equivalently: separable) subgroup $H < G$ the corresponding
$H$-dynamical system $(H,X_{f \rest_H},\eta_0)$, when considered
%23June
%I guess X_{f \rest_H} is exactly $X_f^H$ from 2.4.4
%end
as a subspace of the Banach space $RUC(H)$ (see Proposition
\ref{X_f}.4), is norm separable (that is, $r_H(X_{f \rest_H})
\subset RUC(H)$ is separable;
%oct21
see the definitions before Lemma \ref{l:asp-def}). We denote by
$K(G)$ the collection of Veech functions in $RUC(G)$.
\end{defn}

\begin{thm}\label{Veech}
For any topological group $G$ we have:
\begin{enumerate}
\item
$K(G)$ is a closed left $G$-invariant
subalgebra of $RUC(G)$.
\item
The algebra  $K(G)$ is point-universal.
%3July
%In particular it is also right $G$-invariant.
%end
\item
$ Asp(G) \subset K(G) $.
\item
$K(G)=Asp(G)=Asp_s(G)$ for every separable $G$.
\end{enumerate}
\end{thm}
\begin{proof}
1.\
For every $f\in K(G)$ let $(G,X_f, f)$ be the corresponding
pointed dynamical system as constructed in Proposition \ref{X_f}.
If $f_i,\ i=1,2$ are in $K(G)$ and $H < G$ is a
countable subgroup then the subsets $X_{f_i \rest_H}, i=1,2$ are
norm separable in $RUC(H)$ and therefore so is $X=\{\om + \eta:
\om \in X_{f_1 \rest_H},\ \eta\in X_{f_2 \rest_H}\}$. Since
$X_{(f_1+f_2) \rest_H}\subset X$ it follows that $f_1+f_2\in
K(G)$. Likewise $f_1\cdot f_2\in K(G)$, and we conclude that
$K(G)$ is a subalgebra. Uniformly convergent countable sums are
treated similarly and it follows that $K(G)$ is uniformly closed.
The left $G$-invariance is clear.

%23June
%I guess the present proof of (2) can be generalized to something
%like this:
%if the given algebra $\Acal$ is left and right G-invariant
%and $X_f \subset \Acal$
%(what else ?) then $\Acal$ is point-transitive
%end

 2.\ Given $f\in K(G)$ one shows, as in
\cite[Lemma 3.4]{V}, that every element $\om\in X_f$ is also in
$K(G)$.
%3July
Now use Proposition \ref{Cor-univ2}.
%end

3.\
By Lemma \ref{l:asp-def}, a function $f\in RUC(G)$ is Asplund
iff $r_H(X_f)$ is norm separable in $RUC(H)$ for every countable
subgroup $H < G$. Consider ${\cls} (Hf)$ the $H$-orbit closure
%oct17
in $X_f$ (for $f \in X_f={\cls} (Gf)$).
Then $r_H({\cls} (Hf))$ is also separable in
$RUC(H)$. On the other hand, it is easy to see
that the set $r_H(X_{f \rest_H})$ coincides with $r_H({\cls} (Hf))$.
Hence, $r_H(X_{f \rest_H})$ is also separable in $RUC(H)$. This
exactly means that $f \in K(G)$.

4.\ Let $f \in K(G)$. Then the collection $X_{f \rest_H}$ is norm
separable for every separable subgroup $H < G$. In particular,
$X_f$ (for $H:=G$) is norm separable. Now by Corollary \ref{c:s-asp}
we can conclude that $f \in Asp_s(G)$.
\end{proof}

\br

\section{Hereditary AE and NS systems} \label{s:main}

\sk

We begin with a generalized version of sensitivity. The {\it
functional version} (Definition \ref{d:sens-f}.3) will be
convenient in the proof of Theorem \ref{t:env-asp}.

\begin{defn} \label{d:sens-f}
\begin{enumerate}
\item The uniform $G$-space $(X,\mu)$ has
{\em sensitive dependence on initial conditions\/} (or, simply is
{\it sensitive}) if there exists an $\ep \in \mu$ such that for
every $x\in X$ and any neighborhood $U$ of $x$ there exists $y\in
U$ and $g\in G$ such that $(gx,gy) \notin \ep$ (for metric
cascades see for example \cite{AuYo, De, GW1}).

Thus a (metric) $G$-space $(X,\mu)$ is {\em not sensitive\/}, NS
for short, if for every ($\ep>0$) $\ep \in \mu$ there exists an
open nonempty subset $O$ of $X$ such that $gO$ is $\eps$-small in
$(X, \mu)$ for all $g \in G$, or, equivalently, $O$ is
$[\eps]_G$-small in $(X, \mu_G)$ (respectively: whose
$d_G$-diameter is less than $\eps$, where $d$ is the metric on $X$
and as usual $d_G(x,x')=\sup_{g\in G}d(gx,gx')$).
\item We say that $(G, X)$ is {\em
hereditarily not sensitive\/} (HNS for short) if every nonempty
closed $G$-subspace $A$ of $X$ is not sensitive.
\item
More generally, we say that a
%23June
map $f: (X, \tau) \to (Y,\mu)$ is {\it not sensitive} if there
exists an open nonempty subset $O$ of $X$ such that $f(gO)$ is
$\eps$-small in $(Y,\mu)$ for every $g \in G$.
%gl
The function $f$ is
{\em hereditarily not sensitive\/} if for every
closed $G$-subspace $A$ of $X$ the restricted
function $f \rest_A: A \to (Y,\mu)$ is not sensitive.
Using these notions we can define
the classes of NS and HNS functions. Observe that $(X,\mu)$ is NS
iff the map $id_X: (X,top(\mu)) \to (X,\mu)$ is NS.
%end
\end{enumerate}
\end{defn}

%23June inspired by the referee's suggestions
Let $(X, \mu)$ be a uniform $G$-space and $\eps \in \mu$. Define
$Eq_{\eps}$ as the union of all nonempty $top(\mu)$-open
$[\eps]_G$-small subsets in $X$. More precisely
$$
Eq_{\eps}:=\cup\{U \in top(\mu): \hskip 0.2cm (gx,gx') \in \eps
\quad \forall \hskip 0.2cm (x,x',g) \in U \times U \times G\}.
$$
Then $Eq_{\eps}$ is an open $G$-invariant subset of $X$ and
$Eq(X)=\cap\{Eq_{\eps}: \hskip 0.2cm \eps \in \mu \}$.

\begin{lem} \label{l:NS}
Let $(X,\mu)$ be a uniform $G$-space.
 \ben
 \item
 $X$ is NS if and only if $Eq_{\eps} \neq \emptyset$ for every $\eps \in
 \mu$. Therefore, if $Eq(X) \neq \emptyset$ then $(X, \mu)$ is NS.
\item $X$ is locally $\mu_G$-fragmented iff $Eq_{\eps}$ is dense
 in $X$ for every $\eps \in
 \mu$. Thus, if $X$ is locally $\mu_G$-fragmented then $X$ is NS.
\item
If $X$ is NS then $Eq(X) \supset Trans(X)$.
\item
%gl-
If $X$ is NS and topologically transitive then $Eq(X) = Trans(X)$
and so  $X$ is point transitive iff $Eq(X) \not= \emptyset$.
%If $X$ is NS and topologically transitive then $Eq(X)=Trans(X)$.
\item
If $Eq(X) \neq \emptyset$ and $X$ is topologically transitive then
%$X$ is point transitive, AE and
$Eq(X)=Trans(X)$.
 \een
\end{lem}
\begin{proof}
%gl
The first two assertions are trivial.

3. If $X$ is NS then $Eq_{\eps}$ is not empty for every $\eps \in
\mu$. Any transitive point is contained in any nonempty invariant
open subset of $X$. In particular, $Trans(X) \subset Eq_{\eps}$.
Hence, $Trans(X) \subset \cap\{Eq_{\eps}: \hskip 0.2cm \eps \in
\mu \}=Eq(X)$.

4. By assertion 3 it suffices now to show that if $X$ is
topologically transitive then $Eq(X) \subset Trans(X)$. Let $x_0
\in Eq(X)$, $y \in X$ and let $\eps \in \mu$. We have to show that
the orbit $Gx_0$ intersects the $\eps$-neighborhood $\eps(y):=\{x
\in X: (x,y) \in \eps \}$ of $y$. Choose
%gl
$\delta \in \mu$
such that $\delta \circ \delta \subset \eps$. Since $x_0 \in Eq(X)$
there exists a neighborhood $U$ of $x_0$ such that $(gx_0,gx) \in
\delta$ for every $(x,g) \in U \times G$. Since $X$ is
topologically transitive we can choose $g_0 \in G$ such that $g_0U
\cap \delta(y) \neq \emptyset$. This implies that $(g_0x,y) \in
\delta$ for some $x \in U$. Then $(g_0x_0,y) \in \delta \circ
\delta \subset \eps$.

5. Combine the assertions 1 and 4.
\end{proof}
%end

%Note also that if a $G$-space $X$ has an isolated point
%then trivially every map $f: X \to (Y, \mu)$ is NS.

%gl

%3July maybe it is better to call it "Proposition" and not "Corollary" ?
%end
\begin{cor} \label{wm}
A weakly mixing NS system is trivial.
\end{cor}

\begin{proof}
Let $(G,X)$ be a weakly mixing NS system.
Let $\ep$ be a neighborhood of the diagonal and
%3July
%chose
choose
%end
a symmetric neighborhood of the diagonal $\del$
with $\del\circ \del \circ \del \subset \ep$. By the NS property
%3July
and Lemma \ref{l:NS}.1
%end
$Eq_\del$ is nonempty. Thus there exists a nonempty open subset
$U\subset X$ such that $W=\cup_{g\in G} gU \times g U\subset \del$.
By weak mixing the open invariant set $W$ is dense in $X\times X$
and hence $X \times X \subset \ep$.
Since $\ep$ is arbitrary we conclude that $X$ is trivial.
\end{proof}

 \sk

Next we provide some useful results which link our dynamical and
topological definitions (and involve fragmentability
and sensitivity).

\begin{lem}\label{l-=}
\begin{enumerate}
\item Let $f: X \to Y$ be a $G$-map from a topological $G$-space $(X,\tau)$
into a uniform $G$-space $(Y,\mu)$. Then the following are
equivalent:
\begin{itemize}
\item [(a)]
$f: (X, \tau) \to (Y,\mu)$ is HNS.
\item [(b)]
$f: (X, \tau) \to (Y,\mu_G)$
is fragmented.
\item [(c)]
$f: (A,\tau \rest_A) \to (Y,\mu_G)$ is locally fragmented
for every closed nonempty $G$-subset $A$ of $X$.
\end{itemize}
\item
$(X,\mu)$ is HNS iff ${\id}_X: (X,\tau) \to (X, \mu_G)$ is
fragmented.
\item HAE $\subset$ HNS.
\end{enumerate}
\end{lem}
\begin{proof}
1.\ (a) $\Rightarrow$ (b): Suppose that $f: (X, \tau) \to (Y,\mu)$
is HNS. We have to show that $f$ is $(\tau, \mu_G)$-fragmented.
Let $A$ be a nonempty subset of $X$ and $[\ep]_G \in \mu_G$.
Consider the closed $G$-subspace $Z:={\cls}(GA)$ of $X$. Then by
our assumption $f \rest_Z: Z \to (Y,\mu)$ is NS. Hence there
exists a relatively open nonempty subset $W \subset Z$ such that
$(f(gx),f(gy))=(gf(x),gf(y)) \in \eps$ for every $(g,x,y) \in G
\times W \times W$. Therefore, $f(W)$ is $[\ep]_G$-small. Since
$GA$ is dense in $Z$, the intersection $W \cap GA$ is nonempty.
There exists $g_0 \in G$ such that $g_0^{-1}W \cap A \neq
\emptyset$. On the other hand, clearly, $f(g_0^{-1}W)$ is also
$[\eps]_G$-small. Thus the same is true for $f(g_0^{-1}W \cap A)$.

(b) $\Rightarrow$ (c): Is trivial by Definition \ref{def:fr}.

(c) $\Rightarrow$ (a):
%23June
Let $A$ be a closed nonempty $G$-subspace of $X$ and $\eps \in
\mu$. Take a nonempty open subset $O$ of the space $A$ (say,
$O=A$). Since $f: A \to (Y,\mu_G)$ is locally fragmented one can
choose a nonempty open subset $U \subset O$  such that $f(U)$ is
$[\eps]_G$-small in $Y$. This means in particular, that $f
\rest_A: A \to (Y, \mu)$ is NS for every closed $G$-subspace $A$.
Hence, $f$ is HNS.
%end

2:\ This is a particular case of the first assertion for $f=
{\id}_X: (X, \mu) \to (X, \mu)$.

3:\  Let $(G,X)$ be HAE. For every closed nonempty $G$-subsystem
$A$ there exists a point of equicontinuity of $(G,A)$. By Lemma
\ref{l:NS}.1, $(G,A)$ is NS. Therefore, $(G,X)$ is HNS.
\end{proof}

%In particular, this proposition (apply it for the trivial group)
%implies also that a function $f: (X,\tau) \to (X,\mu)$ is
%fragmented iff $f$ is HNS.

\begin{prop} \label{p:AE=l-fr}
Let $X$ be a compact $G$-system with its unique uniform structure $\mu$.
Consider the following conditions:
 \ben
\item [(a)]
$X$ is AE.
\item [(b)]
$X$ is locally $\mu_G$-fragmented.
\item [(c)]
$X$ is NS.
 \een
Then we have:
 \ben
\item
 Always, $(a) \Rightarrow (b) \Rightarrow (c)$.
\item If $\mu_G$ is metrizable (e.g., if $\mu$ is metrizable)
then $(a) \Leftrightarrow (b) \Rightarrow (c)$.
\item
If $X$ is
%aref (15)  (in some cases I add "point")
point transitive then $(a) \Leftrightarrow (b) \Leftrightarrow
(c)$.
%23June
\item
If $X$ is topologically transitive then (a) $\Rightarrow$ (b)
$\Leftrightarrow$ (c).
 \een
\end{prop}
\begin{proof}
1.\ (a) $\Rightarrow$ (b):\ Let $U$ be a nonempty open subset of
$X$ and $\eps \in \mu$. Since $X$ is AE we can choose a point $x_0
\in Eq(X) \cap U$. Now we can pick an open neighborhood $O \subset
U$ of $x_0$ such that $(gx,gx') \in \eps$ for every $g \in G$ and
$x,x' \in O$. Therefore, $(x,x') \in [\eps]_G$. This proves that
$X$ is locally $\mu_G$-fragmented.

(b) $\Rightarrow$ (c):\ Trivial by Lemma \ref{l:NS}.2.

2.\ (a) $\Leftarrow$ (b): If $\mu_G$ is metrizable then
Proposition \ref{fr-sep-1} guarantees that $id_X: (X,\mu) \to
(X,\mu_G)$ is continuous at the points of a dense $G_{\delta}$
subset (say, $Y$) of $X$. By Lemma \ref{l-equic-p}, $Y \subset
Eq(X)$. Hence, $Eq(X)$ is also dense in $X$. Therefore, $X$ is AE.

3.\ (c) $\Rightarrow$ (a): Observe that $Trans(X) \subset Eq(X)$
by Lemma \ref{l:NS}.3.
% Let $\eps \in
%\mu$. Since $X$ is NS one can choose a nonempty open subset U of X
%being $[\eps]_G$-small. Since $z$ is transitive, we can find $g
%\in G$ such that $gz \in U$. Then $g^{-1}U$ is an open
%neighborhood of $z$ being $[\eps]_G$-small. Thus, $z \in Eq(X)$.

4.\ (b) $\Leftarrow$ (c): Since $X$ is NS the subset $Eq_{\eps}$
is nonempty for every $\eps \in \mu$ (Lemma \ref{l:NS}.1). Since
the open set $Eq_{\eps}$ is invariant and $X$ is topologically
transitive we obtain that $Eq_{\eps}$ is dense for every $\eps \in
\mu$. By Lemma \ref{l:NS}.2 this means that $X$ is locally
$\mu_G$-fragmented.
\end{proof}

The equivalence of AE and NS for transitive metric systems is
shown in \cite{GW1,AAB1}.
%gl
The referee proposed the following problem.
Does there exist a topologically transitive
NS system which is not point transitive?
That is, can it happen for a topologically transitive system
that every $Eq_\ep$ is dense but the intersection $Eq$ is empty?

\begin{cor} \label{c:X_f-lofr}
For every topological group $G$ and $f \in RUC(G)$ the following
are equivalent:
\ben
\item
$f \in light(G)$.
\item
$X_f$ is AE.
\item
$X_f$ is locally norm-fragmented (with respect to the norm of
$RUC(G)$).
\item
$X_f$ is NS.
\een
\end{cor}
\begin{proof}
Use Propositions \ref{p:AE=l-fr}.3 and \ref{light}.2.
%me
It should be noted here that if $\mu$ is the natural pointwise
uniform structure on $X_f={\cls}(Gf) \subset RUC(G)$ then the norm of
$RUC(G)$ induces on $X_f$ the uniform structure $\mu_G$ (Remark
\ref{r:mu-on-X_f}.1).
\end{proof}

\begin{lem} \label{l-quot}
HNS is closed under quotients of compact $G$-systems.
\end{lem}
\begin{proof}
Let $f: X \to Y$ be a $G$-quotient. Denote by $\mu_X$ and $\mu_Y$
the original uniform structures on $X$ and $Y$ respectively.
Assume that $X$ is HNS. Or, equivalently (see Lemma \ref{l-=}.2),
that $X$ is $(\mu_X)_G$-fragmented. Since $f: (X, \mu_X) \to (Y,
\mu_Y)$ is uniformly continuous, it is easy to see that the
$G$-map $f: (X, (\mu_X)_G) \to (Y, (\mu_Y)_G)$ is also uniformly
continuous. We can now apply Lemma \ref{l-quot-fr}. It follows
that $Y$ is
%me
$(\mu_Y)_G$-fragmented.
Hence, $Y$ is HNS (use again Lemma
\ref{l-=}.2).
\end{proof}

\sk

Note that the class NS is not closed under quotients (see \cite{GW1}).

\begin{lem} \label{l-RNis}
\ben
\item
 Every RN compact $G$-system $X$ is HAE. In particular, such a
system is always LE and HNS.
\item
$Asp(X) \subset LE(X)$ for every $G$-space $X$.
 \een
\end{lem}
\begin{proof} 1.
By Definition \ref{d-RN} there exists a representation $(h,\al)$
of $(G,X)$
%aref "into" or "on " ??
%eli
on an Asplund space $V$ such that $h:G \to {\Iso}(V)$
is a co-homomorphism and $\al: (X,\tau) \to (V^*,w^*)$ is a bounded weak$^*$
$G$-embedding. Since $V$ is Asplund, it follows that $\al(X)$ is
(weak$^*$, norm)-fragmented. The map ${\id}_X: (X,\tau) \to (X, norm)$
has a dense subset of points of continuity by Proposition \ref{fr-sep-1}.
 The norm induces on $X$ the metric uniform
structure which majorizes the original uniform structure $\mu$ on
$X$.
%oct17
On the other hand the norm is $G$-invariant. It follows that every
point of continuity of ${\id}_X: (X, \mu) \to (X,norm)$ is a point
of equicontinuity for the system $(G,X)$. Clearly, the same is
true for every restriction on a closed $G$-invariant nonempty
subset $Y$ of $X$. Hence $X$ is HAE. Then clearly $X$ is LE (see
Definition \ref{d:LE=hl}.2). Lemma \ref{l-=}.3 implies that $X$ is
also HNS.

2. Use the first assertion and Theorem \ref{thm1}.3 (taking into
account Definition \ref{d:asp}.2).
\end{proof}

\sk

In the following theorem we show that the classes HNS and
$\rm{RN_{app}}$ coincide. Loosely speaking we can rephrase this by
saying that a compact $G$-system $X$ admits sufficiently many good
(namely: Asplund) representations if and only if $X$ is
``non-chaotic''.

\begin{thm} \label{RN-app=HNS}
For a compact $G$-space $X$ (with its unique compatible
uniform structure $\mu$)
the following are equivalent:
\begin{enumerate}
\item
$X$ is $RN_{app}$.
\item
$X$ is HNS.
\item
$\pi_{\sharp}: X \to C(G,X)$ is a fragmented map.
\item
$\breve{G}=\{\breve{g}: X \to X \}_{g \in G}$ is a fragmented
family.
\item
$(X, \mu_H)$ is uniformly Lindel\"of for every countable
(equivalently, uniformly Lindel\"of) subgroup $H < G$.
\end{enumerate}
\end{thm}
\begin{proof}
1 $\Rightarrow$ 2 :\  The first assertion means that $(X,\mu)$ is
a subdirect product of a collection $X_i$ of RN $G$-systems (with
the uniform structure $\mu_i$). By Lemma \ref{l-RNis}.1 every $X_i$ is
HNS. Lemma \ref{l-=}.2 guarantees that each $X_i$ is
$(\mu_i)_G$-fragmented. Then $X$ is $\mu_G$-fragmented. Indeed,
this follows by Lemma \ref{l-unif-st} and the fact that
fragmentability is closed under passage to products (Lemma
\ref{simple-fr}.5) and subspaces. Now, by Lemma \ref{l-=}.2, $X$
is HNS.

2 $\Leftrightarrow$ 3 :\ $\pi_{\sharp}: X \to C(G,X)$ is
fragmented iff $X$ is $\mu_G$-fragmented. Hence, we can use Lemma
\ref{l-=}.2.

3 $\Leftrightarrow$ 4: See Definition \ref{d:fr-family}.1.

2 $\Rightarrow$ 5 :\ Let $X \in HNS$ and $H < G$ be a uniformly
Lindel\"of subgroup. We have to show that $(X, \mu_H)$ is
uniformly Lindel\"of. The system $(H,X)$ (being m-approximable by
virtue of Proposition \ref{prop:quasisep}) is a subdirect product
of a family of compact metric $H$-systems $\{X_i: i \in I \}$.
Uniform product of uniformly Lindel\"of spaces is uniformly
Lindel\"of. Therefore by Lemma \ref{l-unif-st} it suffices to
establish that every $(X_i, (\mu_i)_H)$ is uniformly Lindel\"of.
Since $\mu_i$ and $(\mu_i)_H$ are metrizable, this is equivalent
to showing that $(\mu_i)_H$ is separable. Since $(H,X)$ is HNS
then, by Lemma \ref{l-quot}, the $H$-quotient $(H, X_i)$ is also
HNS. Hence, $id_{X_i}: (X_i, \mu_i) \to (X_i,(\mu_i)_H)$ is fragmented
by virtue of Lemma \ref{l-=}.2. Now, Lemma \ref{fr-sep-0}
guarantees that $(X_i,(\mu_i)_H)$ is separable.

5 $\Rightarrow$ 1 :\ We have to show that $X$ is $RN_{app}$.
Equivalently, by Theorem \ref{thm1}.4 we need to check that
$C(X)=Asp(X)$. Let $F \in C(X)$ and $H < G$ be a countable
subgroup. By our assumption, $(X, \mu_H)$ is uniformly Lindel\"of.
Since $F: (X, \mu) \to \R$ is uniformly continuous then ${\id}_X:
(X, \mu_H) \to (X,\rho_{H,F})$ is also uniformly continuous.
Therefore, $(X, \rho_{H,F})$ is
uniformly Lindel\"of, too. Since $\rho_{H,F} $ is a pseudometric,
we conclude that $(X, \rho_{H,F})$ is separable. This proves that
$F \in Asp(X)$.
\end{proof}

\begin{remarks} \label{rem:cor}
\begin{enumerate}
\item Every precompact uniform space is uniformly Lindel\"of.
Note here that $(X, \mu_G)$ is precompact iff $(G,X)$ is
equicontinuous (cf. Corollary \ref{c-equic}). Therefore,
$\rm{RN_{app}}$, and its equivalent concept HNS, can be viewed as
a natural generalization of equicontinuity.
\item
Theorem \ref{RN-app=HNS} implies that $\rm{RN_{app}}$ (or, HNS) is
``countably-determined". That is, $(G,X)$ is $\rm{RN_{app}}$ iff
$(H,X)$ is $\rm{RN_{app}}$ for every {\em countable\/} subgroup $H
< G$.
\item Let $H < G$ be a {\em syndetic\/} subgroup
(that is, there exists a compact subset $K \subset G$ such that
$G=KH$) of a uniformly Lindel\"of group $G$. Then a system $(G,X)$
is RN$_{app}$ iff the system $(H,X)$ is RN$_{app}$. Indeed, $K$
acts $\mu$-uniformly equicontinuously on $X$. Thus if $(X, \mu_H)$
is uniformly Lindel\"of then $(X, \mu_{KH})$ is uniformly
Lindel\"{o}f, too.
%\item We get also that always, HAE
%$\subset$ $\rm{RN_{app}}$ (because HAE $\subset$ HNS =
%$\rm{RN_{app}}$ ).
\item
RN$_{app} \subset LE$ by Lemma
\ref{l-RNis}.1 (or, by \cite[Theorem 6.10]{M1}).
\end{enumerate}
\end{remarks}

We now have the following diagram for compact $G$-systems:
\begin{equation*}
\xymatrix {
 Eb \ar[r] \ar[drr] & RN  \ar[r]  & HAE  \ar[r] & HNS  =  RN_{app}
\ar[r] & LE  \\
& & {WAP}\ar[ur] & & }
\end{equation*}

\begin{remark}
\begin{enumerate}
\item
We do not know (even for cascades) if  HAE $ \neq $ HNS for
nonmetrizable systems. All other implications, in general, are
proper:
\item
RN $\neq$ HAE, Eb $\neq$ WAP. Indeed, take a system $(G,X)$ with
trivial $G$ and a compact $X$ which is not RN in the sense of
Namioka, and hence not Eberlein, as a compact space (e.g.
$X:=\beta \N$). Such a $G$-system however trivially is WAP and also
HAE.
\item
Eb $\neq$ RN. Take a trivial action on a compact RN space which is
not Eberlein.
\item
RN$_{app} \neq$ LE even for transitive metric systems (Remark
\ref{exp+}.1 and Theorem \ref{t:GW}).
\item
WAP $\neq$ HNS. See Theorem \ref{t:GW}.
\end{enumerate}
\end{remark}

%oct24
%By Ellis-Nerurkar \cite[Proposition II.8.2]{EN}
%every wap distal $G$-system is equicontinuous.

%\begin{prop} \label{p:distal}
%Let $X$ be an LE compact $G$-system.
%Denote by $E:=E(X)$ the corresponding enveloping semigroup.
%Then the map $E \times E \to E$ is jointly continuous
%at every point $(p,z)$ for every transitive point $z$ of the system
%$(G,E)$ and every $p \in E$.
%\end{prop}
%\begin{proof} Since the class LE is closed under subdirect products
%we get, in particular, that
%the enveloping semigroup $E$, as a $G$-space,
%is also LE. Then $(G,E)$ is AE and every transitive point
%$z \in E$ is a point of equicontinuity with respect to
%the action of $G$ on $E$. Then the same is true
%for the action of the monoid $E$ on $E$. Using the continuity of
%$\rho_z: E \to E, \rho_z(s)=sz$
%it follows that $E \times X \to X$ is jointly continuous at every
%$(p,z)$, where $p \in E$.
%\end{proof}
%%%%%%%%%%%

\begin{thm} \label{thm-faces}
For a compact $G$-system $X$
the following are equivalent:
\begin{enumerate}
  \item  $f \in Asp(X)$.
  \item $f^G_{\sharp}: X \to RUC(G)$ is fragmented.
  %23June
  \item $f_{\sharp}: X \to X_f$ is HNS.
  \item $f: X \to \R$ is HNS.
  %end
  \item $\breve{G}^f:=\{\breve{g}_f: X \to \R \}_{g \in G}$
  (where $\breve{g}_f(x)=f(gx)$) is a fragmented family.
  \item  $X_f \subset RUC(G)$ is norm fragmented.
  \item  The $G$-system $X_f$ is RN.
\end{enumerate}
\end{thm}
\begin{proof} 1 $\Rightarrow$ 2: By Theorem \ref{thm1}.3 there
exist a $G$-quotient $\a: (X, \mu_X) \to (Y, \mu_Y)$ with $Y \in
\rm{RN}$ and $F \in C(Y)$ such that $f= F \circ \a$. Then
$f_{\sharp} =F_{\sharp} \circ \a$. Therefore,
%oct21
by Lemma \ref{simple-fr}.6,
it is enough to show
that $F_{\sharp}: Y \to RUC(G)$ is fragmented, or equivalently,
that $Y$ is $\rho_{G, F}$-fragmented
(see remark before Lemma \ref{l:asp-def}).
By our assumption $(Y,\mu_Y)$ is RN. Therefore,
Theorem \ref{RN-app=HNS} guarantees that
$Y$ is ($\mu_Y)_G$-fragmented. Since ${\id}_Y: (Y, (\mu_Y)_G) \to
(Y, \rho_{G, F})$ is uniformly continuous, it follows that $Y$ is
$\rho_{G, F}$-fragmented, as required.

%3July
2 $\Leftrightarrow$ 3: Use Lemma \ref{l-=}.1 taking into account
Remark \ref{r:mu-on-X_f}.1.
%end

3 $\Leftrightarrow$ 4: Let $f_{\sharp}: X \to X_f$ be HNS. Then
$f_{\sharp}\rest_A: A \to X_f$ is NS for every nonempty invariant
closed subset
of $A \subset X$. Therefore by Definition \ref{d:sens-f}
(observe that the
%gl
uniform structure of $X_f \subset \R^G$ is the pointwise
uniform structure inherited from $\R^G$)
for every $\eps >0$ and every finite subset $S \subset G$
there exists
%6July
%an open nonempty subset $O \subset X$
a relatively open nonempty subset $O \subset A$ such that
$$|f_{\sharp}(gx)(s)-f_{\sharp}(gx')(s)| < \eps
\quad \forall \hskip 0.2cm (s,g) \in S \times G, \quad \forall
\hskip 0.2cm (x,x') \in O \times O.$$
 Now since
$|f_{\sharp}(gx)(s)-f_{\sharp}(gx')(s)|=|f(sgx)-f(sgx')|$ and $g$
runs over all elements of
$G$ our condition is equivalent to the inequality
%6June
$$
|f(gx)-f(gx')| <\eps  \quad \forall \ g \in G.
$$
%$$
%|f(gx)-f(gx') <\eps|  \quad \forall \hskip 0.2cm g \in G.
%$$
%end

%On the other hand the latter means that $f: X \to \R$ is HNS.
The latter means that $f(gO)$ is $\eps$-small for every $g \in G$.
Equivalently, $f: X \to \R$ is HNS.

2 $\Leftrightarrow$ 5:
See Definition \ref{d:fr-family}.1.

2 $\Rightarrow$ 6: Let $f_{\sharp}: X \to X_f$ be the canonical
$G$-quotient. Then by Lemma \ref{l-quot-fr} (with
$Y_1=Y_2=RUC(G)$) the fragmentability of $f^G_{\sharp}: X \to
RUC(G)$ guarantees the fragmentability of $r_G: X_f \to RUC(G)$.
This means that $X_f$ is norm fragmented.

6 $\Rightarrow$ 7: The norm on $RUC(G)$ is lower semicontinuous
with respect to the pointwise topology. Hence, Theorem
\ref{frag=rn} ensures that the $G$-system $X_f$ is RN.

7 $\Rightarrow$ 1: Since $X_f$ is RN, by Theorem \ref{thm1}.3 and
Proposition \ref{p:generalX_f}.3 we get that $f \in Asp(X)$.
\end{proof}

%6July
%Note that $X_f$ is an Eberlein system if (and only if) $f
%\in WAP(X)$. Indeed we can apply \cite[Theorem 4.11]{M1} making
%use of the {\it right strict $G$-duality} (in the sense of
%\cite{M1}) $Y \times X_f \to \R$ where $Y:={\cls}_{w}(fG)$ is the
%weak closure of $fG$ in $WAP(X)$.
%end

\begin{rmk}
Note in the following list how, for a $G$-space $X$,
topological properties of $X_f$
correspond to dynamical properties of $f \in
RUC(X)$ and provide an interesting dynamical hierarchy.
\begin{align*}
&\text{$X_f$ is norm compact $\leftrightarrow$ $f$ is AP}\\
&\text{$X_f$ is weakly compact $\leftrightarrow$ $f$ is WAP}\\
&\text{$X_f$ is norm fragmented $\leftrightarrow$ $f$ is Asplund
$\quad$ }\\
%me interchanging last two lines
&\text{$X_f$ is orbitwise light $\leftrightarrow$ $f$ is LE
$\quad$}\
%oct11  deleting the last line (of the "hierarchy") - pahot mesurbal !
% &\text{$X_f$ is locally fragmented $\leftrightarrow$
%$X_f$ is AE$ \leftrightarrow$ $f$ is light $\quad$ (here $X:=G$)}\\
\end{align*}
\end{rmk}
In the domain of compact metric systems NS and AE are distinct
properties. In contrast to this fact, if these conditions hold
hereditarily then they are equivalent.

\begin{thm} \label{hae=rn=v}
Let $(X,d)$ be a compact metric $G$-space.
The following properties are equivalent.
\begin{enumerate}
\item
X is RN.
\item
X is HAE.
\item
Every closed $G$-subsystem $Y$ of $X$ has a point of
equicontinuity.
\item
$X$ is HNS.
\item
X is
$d_G$-fragmented (recall that $d_G(x,x')={\sup}_{g\in G}d(g x,g
x') $).
\item
$(X,d_G)$ is separable (that is, $d$ is an s-Asplund metric).
\item
Every continuous function $F: X \to \R$ is s-Asplund.
\end{enumerate}
\end{thm}
\begin{proof} Since $X$ is metric,
$\breve{G} \subset Homeo(X)$ is second countable.
So we can and do assume, for simplicity, that $G$ is second countable.

By Theorem \ref{thm1}.5, $RN=RN_{app}$ in the domain of compact
metric systems. Hence, it follows by our diagram above that 1
$\Leftrightarrow$ 2 $\Leftrightarrow$ 4.

2 $\Rightarrow$ 3: Is trivial.

3 $\Rightarrow$ 4: By the assumption $Eq(Y) \neq \emptyset$ for
every subsystem $(G,Y)$. Thus, $Y$ is NS by Lemma \ref{l:NS}.1. It
follows that $X$ is HNS.

4 $\Leftrightarrow$ 5: By Lemma \ref{l-=}.2.

5 $\Rightarrow$ 6: Apply Lemma \ref{fr-sep-0} to the map $id_X:
(X,d) \to (X,d_G)$.

6 $\Rightarrow$ 7: By our assumption $(X,d_G)$ is separable. Since
${\id}_X: (X, d_G) \to (X, \rho_{G,F})$ is uniformly continuous,
we obtain that $(X, \rho_{G,F})$ is also separable. Hence, $f \in
Asp_s(X)$.

7 $\Rightarrow$ 1: Every s-Asplund function is Asplund. Hence,
$C(X)=Asp(X)$. By assertions 4 and 5 of Theorem \ref{thm1} we can
conclude that $X$ is RN.
\end{proof}

Summing up we have the following simple diagram (with two proper
inclusions) for {\em metric\/} compact systems :
$$Eb=WAP\to RN = HAE =HNS = RN_{app} \to LE$$

\br

%aref
%I think "semigroup case" should be checked more carefully.
%so I suggest not include it at all in this paper
%\begin{remark}(Semigroup case) \label{semigroups}
%For a semigroup $S$ acting jointly continuously on a compact space
%$X$ the main results of this section (Theorems \ref{RN-app=HNS},
%\ref{thm-faces}, \ref{hae=rn=v}
%(with the exception of \ref{hae=rn=v}.3)) are still valid after minor
%natural modifications of the definitions of HAE and HNS.
%%oct17
%More precisely, we require that the corresponding conditions hold
%for every closed nonempty, not necessarily
%$S$-invariant, subset $A$ of $X$.
%%we say that an $S$-system $(S,X)$ is HNS (HAE) if the map $(A,\mu) \to
%%(X,\mu_G)$ is locally fragmented (respectively: has a dense subset
%%of continuity points)
%\end{remark}

\br

\section{Some Examples} \label{s:ex}

\sk

\begin{cor}\label{cor-hae-fac-pro}
The class of compact metrizable HNS (hence also RN, HAE) systems
is closed under factors and countable products.
\end{cor}
\begin{proof}  RN=HAE=HNS by Theorem \ref{hae=rn=v}. Now use Lemma
\ref{l-quot} and Theorem \ref{thm1}.5.
\end{proof}

\begin{cor} \label{count}
%oct17
 Every {\em scattered\/} (e.g., countable)
compact $G$-space $X$ is $RN$ (see also \cite{M1}).
\end{cor}
\begin{proof}
%nov3
Apply Theorem \ref{frag=rn} using Remark \ref{r:fr}.4.
\end{proof}

A metric $G$-space $(X,d)$ is called {\em expansive\/} if there
exists a constant $c > 0$ such that
$d_G(x, y):=sup_{g \in G} d(gx,gy)
> c$ for every distinct $x,y \in X$.

\begin{cor} \label{expansive}
An expansive compact metric $G$-space $(X,d)$ is RN iff $X$ is
countable.
\end{cor}
\begin{proof} If $X$ is RN then by Theorem
\ref{hae=rn=v}, $(X,d_G)$ is separable. On the other hand,
$(X,d_G)$ is discrete for every expansive system $(X,d)$. Thus,
$X$ is countable.
\end{proof}

For a countable discrete group $G$ and a finite alphabet $S$
%oct17
the compact space $S^G$ is a $G$-space under left
translations $(g \omega)(h)=\omega (g^{-1}h),\ \omega \in S^G$,
$g,h \in G$. A closed invariant subset $X \subset S^G$ defines a
subsystem $(G,X)$. Such systems are called {\em subshifts\/} or
{\em symbolic dynamical systems\/}.

\begin{cor} \label{shift}
For a countable discrete group $G$ and a finite alphabet
$S$ let $X \subset S^G$ be a subshift.
The following properties are equivalent:
\begin{enumerate}
\item
$X$ is RN.
\item
$X$ is countable.
\end{enumerate}
Moreover if  $X \subset S^G$ is an RN subshift and
$x \in X$ is a
%23June!
recurrent
%agrees with the corrected def of "recurrent" ?
point then it is periodic (i.e. $Gx$ is a finite set).
\end{cor}
\begin{proof}
It is easy to see (and well known) that every subshift is expansive.

For the last assertion recall that if $x$ is a recurrent point
with an infinite orbit then its orbit closure contains a
homeomorphic copy of the Cantor set.
\end{proof}

For some (one-dimensional) compact spaces every selfhomeomorphism
will produce an RN system.

\begin{prop}\label{interval}
\begin{enumerate}
\item
For each element $f\in \Homeo(I)$, the homeomorphism
group of the unit interval $I=[0,1]$, the
corresponding dynamical system $(f,I)$ is HNS.
\item
For each element $f\in \Homeo(S^1)$, the homeomorphism
group of the circle $S^1=\{z\in \C: |z|=1\}$, the
corresponding dynamical system $(f,S^1)$ is HNS.
\end{enumerate}
\end{prop}

\begin{proof}
1.\
Fix an element $f\in \Homeo(I)$, which with no loss
of generality we assume is orientation preserving.
Consider the dynamical system $(f,I)$
and for a set $A\subset I$ denote
$O_f(A)=\cup_{n\in \Z}f^{n}(A)$.
Let us note first that for every $x\in [0,1]$
the sequence $ \dots f^{-2}(x),f^{-1}(x),x,f(x),
\allowbreak f^{2}(x),
\dots$ is monotone increasing hence
the orbit closure of $x$ is just the orbit together with
the points $\lim_{n\to \infty} f^{-n}(x)$ and
$\lim_{n\to \infty} f^{n}(x)$. In particular
the dynamical system $(f,I)$ is LE.

Next we will show that $(f,I)$ is NS. If this is
not the case then there exists an $\epsilon>0$ such that
for every non empty open set $U\subset I$ there exists
$n\in \Z$ such that ${\diam}(f^{n}U) \ge \epsilon$.
Let $(a,b) \subset I$ be an open interval and let
$\{U_k\}_{k\in \N}$ be a countable basis for open sets in $(a,b)$.
If for every $k$ the set $(a,b) \cap O_f(U_k)$
is dense in $(a,b)$ then the orbit of any point $x \in
(a,b) \cap\bigcap_{k=1}^\infty  O_f(U_k)$
will be dense in $(a,b)$ which is impossible.

We conclude that for every interval $(a,b)$ and
every proper subinterval $J_1$ there is another
subinterval $J_2\subset
(a,b)$ which is disjoint from $O_f(J_1)$.
By induction we can find an infinite sequence of disjoint
intervals $J_j$ in $(a,b)$ such that for every $j$ the set
$J_{j+1}$, and hence also $O_f(J_{j+1})$,
is disjoint from $\cup_{i\le j}O_f(J_i)$.
Since for each $j$ the set $O_f(J_j)$ contains an interval
of length at least $\epsilon$ we arrive at a contradiction.
This concludes the proof that $(f,I)$ is NS.

Next consider any nonempty closed invariant subset $Y\subset I$.
If $Y$ contains an isolated point then clearly the system
$(f,Y)$ is NS. Thus we now assume that $Y$ is a perfect set.
We can then repeat the argument that showed that $(f,I)$ is NS
for the system $(f,Y)$ and arrive at the same kind of contradiction
since again an orbit of a single point in $Y$ can not be everywhere
dense in a nonempty set of the form $(a,b)\cap Y$.

\br

2.\
We will use Poincar\'e's classification of the systems $(S^1,f)$
whose nature is well understood (see for example \cite{KH},
Section 11.2).
Again we can assume with no loss of generality that
our homeomorphism $f$ preserves the orientation on $S^1$.
Let $r(f)\in \R$ denote the rotation number of $f$.
If $r(f)$ is rational then some power of
$f$ has a fixed point and we are reduced to the case
of a homeomorphism of $I=[0,1]$. Thus we can assume that
$r(f)$ is irrational. There are two cases to consider.
The first is when the system $(S^1,f)$ is minimal, in which
case $f$ is conjugate to an irrational rotation and is
therefore equicontinuous.

In the second case, when  $(S^1,f)$ is not minimal,
there exists a unique minimal subset $K\subset S^1$
with $K$ a Cantor set and there are wandering intervals $J\subset S^1$.
For such an interval, given an $\ep>0$ there exists an $N$
such that for every $n\in \Z$ with $|n|\ge N$, ${\diam}(f^n(J))< \ep$;
whence the NS property of $(S^1,f)$.

For the HNS property consider an arbitrary subsystem $(Y,f)$
with $Y \subset S^1$. Again distinguish between the cases
when $Y$ has an isolated point and when it is a perfect set.
The presence of an isolated point ensures
NS. Finally when $Y$ is
perfect it is either equal to $K$, hence equicontinuous,
or we can still use the existence
of the wandering intervals in $(S^1,f)$ to obtain a nonempty set
$J\cap Y$ with the property that the diameter
of its images under the iterates of $f$ tends to zero.
\end{proof}

%\begin{prop}\label{wm+ae}

%23June

%gl (see Corollary \ref{wm})
%Let $X$ be
%a weakly mixing compact $G$-system.
%Then every continuous NS function $f: X \to \R$
%is necessarily trivial. In particular,
%$Asp(X)=\{\text{constants}\}$ and every weakly mixing RN$_{app}$
%system is trivial.
%\end{prop}
%\begin{proof}
%Let $f: X \to \R$ be a continuous NS function.
%Assume in contrary that $f(a) \neq f(b)$ for some
%$a,b \in X$. Denote by $\eps_0$ the positive number $\eps_0:=|f(a)-f(b)|$.
%Since $f$ is NS one can choose a nonempty open subset
%$O \subset X$ such that
%$|f(gx)-f(gy)|< \frac{\eps_0}{3}$ for every $x,y \in O$.
%By weak mixing there is a transitive
%point $(x_0,y_0) \in O \times O$. Taking into account the continuity
%of the function $f$
%one can choose $g \in G$ such that $|f(a)-f(gx_0)|< \frac{\eps_0}{3}$ and
%$|f(gy_0)-f(b)| <  \frac{\eps_0}{3}$.
%Then $|f(a) - f(b)| < \eps_0$, a contradiction.
%\end{proof}
%end

%Fix a metric $d$ on $X$ and let $0 < \ep < \frac{1}{2} {\diam}
%(X)$. Let $x_0\in X$ be an equicontinuity point of $X$. Let $\del
%>0$ satisfy: $x\in B_\del(x_0) \Rightarrow d(gx_0,gx)< \ep  \
%\text{for all} \ g\in G$. By weak mixing there is a transitive
%point $(z,w) \in B_\del(x_0) \times B_\del(x_0)$
%%me
%in $X \times X$. Since for every $g\in G$
%$$
%d(gz,gw)\le d(gz,gx_0) + d(gx_0,gw) < 2\ep,
%$$
%we conclude that $\diam(X)\le 2\ep$, contradicting
%our choice of $\ep$.
%\end{proof}

\begin{examples} \label{ex-simple}
Of course it is easy to find non-RN metric
 systems. Here are some
``random" examples.
\begin{enumerate}
\item
The cascades on the torus $\T^2$ defined by a hyperbolic
automorphism, or the horocycle flows,
%3July
being weakly mixing (see Corollary \ref{wm}), are not RN. Likewise
Anosov diffeomorphisms on a compact manifold, being expansive (see
\cite{AH}), are not RN by Corollary \ref{expansive}.
\item
Systems which contain non-equicontinuous minimal subsystems
fail to be RN.
\item
Let $X$ be compact metric and uncountable and set $G={\Homeo}(X)$.
Then in many cases (like $X=[0,1]$) the action is expansive, hence
not RN
%me
(Corollary \ref{expansive}).
\item
As we have seen, any uncountable subshift is not RN.
Thus, for example the well known ``generator of the Morse cascade"
$$
w=\ldots 0110100110010110 \dot{ 0}110100110010110 \ldots
$$
considered as a function $w:\Z\to \R$
is not an Asplund function
on the group $\Z$.
%aref
%\item
% In Proposition \ref{interval}.2 it is essential that $f: S^1
%\to S^1$ is a homeomorphism. Indeed, defining $f(z)=z^2$, we get a
%sensitive action of the semigroup $\N$ on $S^1$. Thus, the
%%me
%corresponding semigroup action $(\N, S^1)$ is not RN.
\end{enumerate}
\end{examples}

\br

A
%aref (17)
point
transitive LE system is, by definition, AE but there are
nontransitive LE systems which are not AE.

\begin{exa}\label{exp}
As can be easily seen
the $\Z$ system $(T,D)$, where $D=\{z\in \C:
|z|\le 1\}$ is the unit disk in the complex plain and $T: D\to D$
is the homeomorphism given by the formula $Tz=z\exp(2\pi i |z|)$,
is an LE system which is not AE.
\end{exa}

There exist many compact metrizable
transitive AE systems which fail to be HAE. This follows, for
example, from the following lemma. We will use the following
construction which is due to Takens.
%23June
%any reference ?
%end
For a metric cascade $(T,X)$ define an {\it asymptotic
pseudo-orbit\/} to be a bi-infinite sequence $\{x_n\}$ such that
$\lim_{|n|\to\infty}d(Tx_n, x_{n+1})=0$. Note that $(T,X)$ is
chain transitive iff it admits an asymptotic pseudo-orbit with
alpha and omega limit point sets the whole space.

%eli

\begin{lem}\label{trans}
Let $(T,X)$ be a metric cascade.
\begin{enumerate}
\item
If $(T,X)$ is chain recurrent $\Z$-space then $X$ is isomorphic
to a subsystem
of a compact metric transitive AE cascade $(T,Y)$.
\item
If $(T,X)$ is transitive-recurrent then $X$ is
also a retract of the ambient transitive AE system $(T,Y)$.
\end{enumerate}
\end{lem}
\begin{proof}
Let $\{t_n\}$ be a bi-infinite monotonic sequence in $(0,1)$
with $\lim_{n\to\infty}t_n =1,\
\lim_{n\to\infty} t_{-n}\allowbreak=0$. Let $S$ be the circle
represented as the interval $[0,1]$ with
$0$ identified with $1$.
%23June
%What is x_n here ? I guess the asymptotic pseudo-orbit ...
%gl
Let $\{x_n\}$ be an asymptotic
pseudo-orbit in $X$. Identify $X$ with the subset $X \times \{0\}
\subset X \times S$ and let $Y= X \cup \{(x_n,t_n): n\in \Z\}$.
Extend $T$ to $Y$ by $T(x_n,t_n)=(x_{n+1},t_{n+1})$.
This completes the proof of part 1.
For part 2 note that if the pseudo orbit is actually an
orbit then the first coordinate projection from $Y$ to $X$
is a $\Z$-retraction.
\end{proof}

\begin{remarks}\label{exp+}
\ben
\item
 If we apply the construction of Lemma \ref{trans} to the
(clearly chain recurrent) system $(T,X)=(T,D)$ of Example
\ref{exp} we obtain a transitive (but not recurrent-transitive)
metric LE system $(T,Y)$ which is not HAE (or, RN$\rm{_{app}}$).
Applying Lemma
%eli
\ref{trans} to a transitive non AE system $(T,X)$ (e.g. a
minimal weakly mixing system), we obtain an example of an AE
system with both a subsystem and a factor which are not AE (see
\cite{GW1}).
\item
As noted above, HAE is preserved under both passage to subsystems
and the operation of taking factors. In the next section we will
show that the Glasner-Weiss family of recurrent-transitive LE but
not WAP systems consists, in fact, of HAE systems. On the other
hand, in Section \ref{LE-HAE} we will modify these examples so
that the resulting dynamical system will still be
recurrent-transitive, LE, but no longer HAE. Thus even among
metric recurrent transitive $\Z$-systems we have the proper
inclusions
$$
WAP \subset HAE \subset LE.
$$
Then we can conclude that the following inclusions are also proper
$$
WAP(\Z) \subset Asp(\Z) \subset LE(\Z).
$$
%23June  The following remarks move here
%I deleted the first remark - I think that it is not interesting enough
%\item
%Note that a subshift $X \subset S^G$ is equicontinuous iff $X$ is finite.
%Once
%again (Remark \ref{rem:cor}.1) we see the analogy between
%``equicontinuous --- finite" and ``RN --- countable".
\item
It is interesting to compare some of the current definitions of
chaos and the corresponding classes of dynamical systems (see, for
example, \cite{De, GW1, BGKM}) with the class of $G$-systems $X$
such that $Asp(X)=\{\text{constants}\}$. The latter are the
systems which admit only trivial representations on Asplund Banach
spaces. Every weakly mixing compact system belongs to this class
%3July
because by
%gl
Corollary \ref{wm} every Asplund function (in fact, every
continuous NS function) on such a system is constant.
%end
\item
By Theorem 1.3 of \cite{GW1} and the variational principle an LE
(e.g., RN) cascade has topological entropy zero. This probably
holds for a much broader class of acting groups but we have not
investigated this direction.
%aug18 Discussion: is it true that every Li-Yorke chaotic
%system is sensitive ?
% (as an alternative, possibly, a sharper proof)
%aug18 Probably factor of Li-Yorke system is not Li-Yorke - right ?
%what is known about FACTORS of "chaotic" (in one of known senses) ?
%I suggest the following
%DEFINITION. Let's say that the system $X$ is "chaotic"
%(or, "GM-chaotic"  :-)
%if Asp(X) = {constants}
%I believe that this definition of chaos is "categorically correct"
%(in contrast to existing, in fact, empiric definitions, which mostly are
%inspired by very important but exceptional case of 1-dimensional systems)
%1. The class of GM-chaotic systems is closed with respect to factors.
%2. What can be said about GM-chaotic groups ?
%(That is, the case when the greatest ambit G^{RUC} is GM-chaotic)
% We know that there exists a non-trivial GM-chaotic group
%(I mean of course, G:=H_+[0,1])
\een
\end{remarks}
%end

\br

\section{The G-W examples are HAE}
\label{Sec-GW}

\sk

In this section we assume that the reader is familiar
with the details of the paper \cite{GW}. In particular
we use the notations of that paper with no further comments.

\begin{thm} \label{t:GW}
The G-W examples of recurrent-transitive LE but not WAP systems
are actually HAE.
\end{thm}

\begin{proof}
Recall, $\Omega$ is the space of continuous maps $x:\R \to 2^I$,
where $I=[0,1]$ and $2^I$ is the compact metric space of closed
subsets of $I$ equipped with the Hausdorff metric $d$. (In fact,
the values $x$ assumes are either intervals or points.) The
topology on $\Omega$ is that of uniform convergence on compact
sets: $x_n\to x$ if for every $\ep>0$ and every $M>0$ there exists
$N>0$ such that for all $n>N$,\ $\sup_{|t|\le
M}d(x_n(t),x(t))<\ep$.
%3July is not "compact"
%This topology makes $\Omega$ a compact metrizable space.
On $\Omega$ there is a natural $\R$-action defined by
translations: $(T^tx)(s)=x(s+t)$.
The
%3July
compact
%REASON ?
metrizable
dynamical system $(T,X)$,
%gl++
where $T=T^1$, is obtained as the orbit
closure $X={\cls} \{T^n\om:n\in \Z\}$ for a carefully constructed
(kite-like) element $\om \in \Om$ (see also the figure in Section
\ref{LE-HAE}).
%gl++
The fact that the function $\om :\R \to 2^I$ is
a Lipschitz function implies that each member of $X$
is Lipschitz as well with the same constant so
that $X$ as a family of functions is equicontinuous.
The compactness of $X$ follows from Arzela-Ascoli theorem.
We next sum up some of the salient facts we have
about $(T,X)$.

\begin{enumerate}
\item[(a)]
For every $x\in X$ there is a unique interval $[a,b]
\subset [0,1]$ such that:
\begin{enumerate}
\item[(i)]
$x(t)\subset [a,b],\ \forall t\in \R$,
\item[(ii)]
there exists a sequence $t_l\in \R$ with
$\lim\ x(t_l)=[a,b]$.
\end{enumerate}
We denote
$$
\mathbf{N}(x)=[a,b].
$$
\item[(b)]
The function $x\mapsto \mathbf{N}(x)$ is
lower semicontinuous; i.e.
$\lim_{\nu}x_{\nu}=x \Rightarrow
\liminf_{\nu}\mathbf{N}(x_{\nu})\supset \mathbf{N}(x)$.
\item[(c)]
Let us call intervals $[a,b]\subset [0,1]$ of the form
$\mathbf{N}(x),\ x\in X$,
{\em admissible\/} intervals. Then for every
admissible $[a,b]\subset [0,1]$
there exists a unique element $\om_{ab}\in X$
with $\mathbf{N}(\om_{ab})=\om_{ab}(0)=[a,b]$.
(In particular $\om_{01}=\om$.)
\item[(d)]
Let $J=\{\om_{ab}\in X: 0 \le a \le  b \le 1\}$. Then $J$ is a
closed subset of $X$ and $\mathbf{N}: J \to \{(a,b): 0 \le a \le b
\le 1\} \subset [0,1]\times [0,1]$ is a homeomorphism onto the set
of admissible intervals. (Not every subinterval of $[0,1]$ is
%me
admissible. For example neither $[0,9/10]$ nor any degenerate
interval with $9/10 < a=b \le 1$ is attained.)
\item[(e)]
Denoting $X_{ab}={\OCT}(\om_{ab})$ we have
$x\in X_{ab}$ iff $\mathbf{N}(x)\subset [a,b]$.
\item[(f)]
For each admissible interval $[a,b]\subset [0,1]$ the subsystem
$(T,X_{ab})$ is AE, with $Eq(X_{ab}) =\{x\in X:
\mathbf{N}(x)=[a,b]\}$.
\end{enumerate}

\br

These facts, perhaps
with the exception of item (b), are either
stated explicitly and proved in \cite{GW} or can be
easily deduced from it. For completeness we provide
a proof for (b).

\begin{proof}[Proof of (b)]
With no loss in generality we assume that
$\liminf_{\nu}\mathbf{N}(x_{\nu})
= \lim_{\nu}\mathbf{N}(x_{\nu}) =[a,b]$
and we then have to show that $[a,b] \supset \mathbf{N}(x)$.
There exists a sequence $m_i$ such that
$\lim_i T^{m_i}x(0)=\mathbf{N}(x)$. Therefore, given $\ep>0$,
there exists an $i$ with
\begin{equation}\label{1}
d(T^{m_i}x(0),\mathbf{N}(x))< \ep.
\end{equation}
Next choose $\nu$ such that
\begin{equation}\label{2}
d(T^{m_i}x_\nu(0),T^{m_i}x(0))< \ep
\end{equation}
and
\begin{equation}\label{3}
d(\mathbf{N}(x_\nu),[a,b])< \ep.
\end{equation}
Now, by \eqref{3} we have
$$
[a-\ep,b+\ep]\supset \mathbf{N}(x_\nu)\supset T^{m_i}x_\nu(0),
$$
hence by \eqref{1} and \eqref{2}
$$
[a-3\ep,b+3\ep]\supset \mathbf{N}(x).
$$
Since $\ep>0$ is arbitrary we conclude that indeed $[a,b]
\supset \mathbf{N}(x)$.
\end{proof}

\br

Of course this list implies the LE property of $(T,X)$. However,
we are after the stronger property HAE.
For this end consider now an arbitrary
closed invariant nonempty subset $Y$ of $X$. Let $J_Y$ be the
subset of $Y$ which consists of those elements $y\in J\cap Y$ for
which $\mathbf{N}(y)=y(0)$ is maximal; that is,
if $z\in Y$ and $\mathbf{N}(z)\supset \mathbf{N}(y)$
then $\mathbf{N}(z)=\mathbf{N}(y)$.

\br

\noindent{\it Claim 1:}\ The restriction of $\mathbf{N} \rest_Y: Y
\to [0,1]\times [0,1]$ is continuous at points of $J_Y$.

\begin{proof}
Suppose $Y \ni y_n \to y \in J_Y$. By the lower semicontinuity
of $\mathbf{N}$,
$$
[a,b] = \liminf_n \mathbf{N}(y_n) \supset \mathbf{N}(y).
$$
Choose a subsequence $n_i$ such that $\mathbf{N}(y_{n_i})
\to [a,b]$. Then for some sequence $m_i$ we have
$T^{m_i}y_{n_i}(0) \to [a,b]$. By compactness we can assume
with no loss in generality that $T^{m_i}y_{n_i} \to z$
for some $z\in Y$.
Now, $T^{m_i}y_{n_i}(0) \to z(0)=[a,b]\supset \mathbf{N}(y)$,
whence $\mathbf{N}(y)=[a,b]$. It follows easily that
$\lim_n \mathbf{N}(y_n) = \mathbf{N}(y)$.
\end{proof}

\br

In item (d) of the above list we noted that $J$ is a closed subset
of $X$ and $\mathbf{N}: J \to [0,1]\times [0,1]$ is a
homeomorphism into. Set $K=\mathbf{N}(J\cap Y)$ and
let $K_0\subset K$ be the subset of maximal elements in $K$;
i.e. $[a,b]\in K_0$ iff $[a,b]\in K$ and $K \ni [c,d] \supset [a,b]$
implies $[c,d]=[a,b]$. Clearly $K_0$ is a closed subset
of the closed set $K$ and for every $[c,d]\in K$ there exists
some $[a,b]\in K_0$ with $[c,d]\subset [a,b]$.

\br

\noindent{\it Claim 2:}\
$K_0=\mathbf{N}(J_Y)$.

\begin{proof}
Let $[a,b]$ be an element of $K_0$, then $[a,b]=\mathbf{N}(y)$ for
some $y\in J\cap Y$. If $[c,d]=\mathbf{N}(z)\supset [a,b]$ for
some $z \in Y$, then for some $z'\in {\OCT}(z)\subset Y$ we have
$z'(0)=[c,d]=\mathbf{N}(z')$. In particular $z'\in J\cap Y$ and
$\mathbf{N}(z')=[c,d]\in K$. Hence $[c,d]=[a,b]$ and it follows
that $y\in J_Y$.

Conversely, if $y\in J_Y$ with $y(0)=[a,b]=\mathbf{N}(y)$
and $\mathbf{N}(z)=z(0)=[c,d]\supset [a,b]$ for $z\in Y$, then
$[c,d]=[a,b]$ and $[a,b]\in K_0$.
\end{proof}

\br

\noindent{\it Claim 3:}\
$J_Y$ is closed and nonempty; in fact
$Y={\cls}\{T^n J_Y: n\in \Z\}$.

\begin{proof}
The fact that $J_Y$ is closed and nonempty is a
direct consequence of Claim 2. Clearly
$\mathbf{N}(Y)= \mathbf{N}(J\cap Y)= K$ and it follows
that every $[a,b] = \mathbf{N}(y)\in \mathbf{N}(Y)$ is a subset of some
$[c,d] =  \mathbf{N}(\om_{ab}) \in K_0$.
By item (e) we have $y\in X_{ab}={\OCT}(\om_{ab})$ and
our claim follows.
\end{proof}

\br

\noindent{\it Claim 4:}\ Every $\om_{ab} \in J_Y$ with $a< b$ is
in $Eq(Y)$.

\begin{proof}
The key fact in proving the inclusion $J_Y \setminus$ \{constant
functions\} $\subset Eq(Y)$ is a certain uniformity of the
function $\ep'=\ep'(\ep,b-a)$ provided by Lemma 3.5 of \cite{GW}.
In essence, as can be seen by combining Lemmas 3.5, 3.6 and 1.1 of
\cite{GW}, this function is the equicontinuity modulus function
for $D(z,w)={\sup}_{n\in \Z}d(T^nz,T^nw)$ on orbit closures in
$(T,X)$; i.e. given a point $x\in X$ with $\mathbf{N}(x)= [a,b]$
and $\ep>0$, the $\ep'$-neighborhood of $x$,
$B_{\ep'}(x) \cap  {\OCT}(x)$, in ${\OCT}(x)$ is $(\ep,D)$-small.
The point is,
that the $\ep'=\ep'(\ep,b-a)$ provided by Lemma 3.5 of \cite{GW}
is uniform in $x$ as long as $b-a$ is bounded away from zero.

Therefore, given a point $\om_{ab}\in J_Y$, with $a< b$,
and $\ep > 0$ we can choose a point $\om_{a'b'}\in J$
with $a' < a < b < b'$ so that  $a-a', b'-b$ are sufficiently
small to ensure that $\om_{ab}\in B_{\ep'}(\om_{a'b'})$.
Of course by (e) we have $\om_{ab}\in {\OCT}(\om_{a'b'})$.

By Claim 1, $\om_{ab}$ is a continuity point for the restriction of
the map $\mathbf{N}$ to $Y$ and it follows that
there exists a neighborhood $V$ of $\om_{ab}$ such that
$\mathbf{N}(y)\subset [a',b']$ for every $y\in V$,
hence $y \in {\OCT}(\om_{a'b'})$.
We now conclude that $B_{\ep'}(\om_{a'b'})\cap V$ is
an $(\ep,D)$-small neighborhood of $\om_{ab}$ {\em in
the subsystem\/} $Y$ and the proof that $\om_{ab}$ is
an equicontinuity point of the system $(T,Y)$ is complete.
\end{proof}

%gl++

\br

We next observe that $T$ acts as the identity on the
open subset
$$
U=Y\setminus {\cls}\{T^n \om_{ab}: \om_{ab}\in J_Y, a < b, n\in \Z\}
$$
(when non-empty) and thus every point in $U$ is an equicontinuity
point. This observation together with Claims 3 and 4 show that the
set $Eq(Y)$ of equicontinuity points is dense in $Y$. That is, $(T,Y)$
is an AE system and our proof of the HAE property of $(T,X)$ is
complete.
\end{proof}

\br

%\br
%
%%3July
%By Claim 4 and
%%Assertion 3 of
%Theorem \ref{hae=rn=v} we can conclude that $(T,X)$ is an HAE
%system.

\section{The mincenter of an RN system} \label{s:mincenter}

\sk

Unlike the case of
%eli
%topologically
transitive WAP systems, where the
{\it mincenter} (i.e. the closure of the union of the minimal
subsets of $X$) consists of a single minimal equicontinuous
subsystem, the mincenter of a
%eli
%topologically
transitive RN system
need not be minimal. In the G-W examples the mincenter consists of
a continuum of fixed points;
moreover, as was shown in \cite{GW} a slight modification of
the construction there will yield examples of HAE systems whose
mincenter consists of uncountably many nontrivial minimal
equicontinuous subsystems all isomorphic to a single
circle rotation. However, in Section \ref{LE-HAE} we will
present a more sophisticated modification which produces an example
of an LE system with a mincenter containing
uncountably many non-isomorphic rotations.
In the present section we obtain some
information about the mincenter of
%eli
%metrizable
RN systems.
This will be used in the next section to draw
a sharp distinction between LE and HAE systems.
For simplicity we deal with metrizable systems. Recall that for
such systems RN is the same as HAE.

The {\em prolongation\/} relation ${\Prol}(X)\subset X\times X$ of
a compact dynamical system $(G,X)$ is defined as follows:
%me I think "sequences" instead of "nets" is enough here (metric systems)
\begin{gather*}
{\Prol}(X)=\\
\{(x,x'): \exists \ {\text{nets $g_\nu\in G$ and $x_\nu\in X$ such
that $\lim_\nu x_\nu=x$ and $\lim_\nu g_\nu x_\nu=x'$}}\}.
\end{gather*}
It is easy to verify that ${\Prol}(X)$ is a closed symmetric and
$G$-invariant relation. For $x_0\in X$ we let
$$
{\Prol}[x_0]=\{x\in X: (x_0,x)\in {\Prol}(X)\}.
$$
Note that always $\OCG(x)\subset {\Prol}[x]$, and if
$x_0\in\OCG(x)$ then $x\in {\Prol}[x_0]$. For closed invariant
sets $A\subset B\subset X$ we say that $A$ is {\em capturing in
$B$\/} if $x\in B$ and $\OCG(x)\cap A\ne \emptyset$ imply $x\in A$
(see \cite{AuG}).

\begin{lem}\label{lem-1}
\begin{enumerate}
\item
Let $(X,d)$ be a metric $G$-system, $x_0\in Eq(X)$ and $x\in
{\Prol}[x_0]$, then $x\in \OCG(x_0)$. Hence,
$$
{\Prol}[x_0]=\OCG(x_0).
$$
\item
If $x_0\in Eq(X)$ and $x_0\in \OCG(x)$, then $x\in Eq(X)$ and
$x\in \OCG(x_0)$; that is, $Eq(X)$ is a capturing subset of $X$.
\end{enumerate}
\end{lem}
\begin{proof}
1. Given $\ep>0$ there exists $\del>0$ such that $z\in
B_\del(x_0)$ implies $d_G(x_0,z)<\ep$. There are nets $g_\nu\in G$
and $x_\nu\in X$ such that $\lim_\nu x_\nu=x_0$ and $\lim_\nu
g_\nu x_\nu=x$. For sufficiently large $\nu$ we have $x_\nu\in
B_\del(x_0)$ and $d(g_\nu x_\nu, x)< \ep$, hence
$$
d(g_\nu x_0,x) \le d(g_\nu x_0, g_\nu x_\nu) + d(g_\nu x_\nu, x) <
2 \ep,
$$
hence $x \in \OCG(x_0)$. Thus ${\Prol}[x_0]\subset \OCG(x_0)$. The
inclusion ${\Prol}[x_0] \supset \OCG(x_0)$ is always true.

2. Given $\ep>0$ there exists a $\del>0$ such that $d_G(x_0,z) <
\ep$ for every $z\in B_\del(x_0)$. There exists $g\in G$ with
$gx\in B_\del(x_0)$ and therefore an $\eta>0$ with
$gB_\eta(x)\subset B_\del(x_0)$. Now for every $h\in G$ and $w\in
B_\eta(x)$ we have
$$
d(hgx,hgw)< d(hgx,hx_0) + d(hgw,hx_0) < 2\ep.
$$
Thus also $x\in Eq(X)$. By assumption $x_0\in \OCG(x)$ hence $x\in
{\Prol}[x_0]$ and by part 1, $x\in \OCG(x_0)$.
\end{proof}

\begin{prop} \label{mincenter}
Let $(X,d)$ be a metrizable RN $G$-system, $M$ its mincenter. Then
$Eq(M)$ is a disjoint union of minimal equicontinuous systems,
each a capturing subset of $M$.
\end{prop}

\begin{proof} Our system $X$ is HAE by Theorem \ref{hae=rn=v}.
Therefore the subsystem $(G,M)$ is AE. Let $x_0\in M$ be an
equicontinuity point of $M$. Given $\ep>0$ there exists a $0<\del<
\ep$ such that $x\in B_\del(x_0)\cap M$ implies $d(g x_0,g x)<
\ep$ for every $g\in G$. Let $x'\in B_\del(x_0)$ be a minimal
point. It then follows that $S=\{g\in G: g x'\in  B_\del(x_0)\}$
is a syndetic subset of $G$ (i.e. $FS=G$ for some finite subset
$F$ of $G$). Collecting these estimations we get, for every $g\in
S$,
$$
d(g x_0,x_0)\le d(g x_0,g x') + d(g x',x_0) \le 2\ep.
$$
Thus for each $\ep>0$ the set $N(x_0,B_\ep(x_0))= \{g\in G: d(g
x_0,x_0) \le \ep\}$ is syndetic, whence $x_0$ is minimal.

Thus every equicontinuity point $x_0$ of $M$ is minimal and we
apply Lemma \ref{lem-1} to conclude that $Eq(M)$ is a capturing
subset of $M$.
\end{proof}

\begin{cor}
The mincenter $Z$ of a metrizable RN system $(G,X)$ is transitive iff
$Z$ is minimal and equicontinuous.
\end{cor}

\begin{rmk}\label{non-wandering}
The {\em Birkhoff center\/} $Y$ of a compact metrizable
$\Z$-dynamical system $(T,X)$ can be defined as the closure of its
recurrent points. A nonempty open set $U\subset X$ such that $T^j
U\cap U= \emptyset$ for all $j\in \Z\setminus \{0\}$ is called a
{\em wandering set\/}. The complement of the union of all
wandering sets is a closed invariant subsystem $Z_1\subset X$
which contains $Y$. Repeating this process (countably many times)
we get by transfinite induction a countable ordinal $\eta$ such
that $Z_\eta=Y$. Since an isolated transitive point of any compact
metric system is always an equicontinuity point it follows easily
that the system $(T,X)$ is LE iff its Birkhoff center $(T,Y)$ is
LE.
%eli
The same statement does not hold for RN systems. An example
of a compact sensitive system $(T,X)$ whose Birkhoff center
%gl+
consists of fixed points was shown to us by
E. Akin (private communication).
\end{rmk}

\br

\section{A recurrent transitive LE but not HAE system}
\label{LE-HAE}

\sk

As promised in Section \ref{s:ex} we will sketch in the present
section a modification of the G-W construction
that will yield a recurrent-transitive system which is LE but not HAE.
The possibility of introducing such a modification
(in order to achieve another goal) occurred to
the authors of \cite{GW} already at the time when this paper
was written. The first author (E.G.)
would like to thank B. Weiss for his help in checking the
details of the modified construction.

\begin{thm}\label{LEnotHAE}
There exists a recurrent-transitive LE but not HAE system.
\end{thm}

\begin{proof}
In the original construction the basic ``frames" $\al_n$
were defined by the formula
$$
\al_n(t) = \al_0\left(\frac t {p_n}\right),\quad n=1,2,\dots
$$
where $\al_0$ is the original periodic kite-like
function,
\br

\centerline{
\epsfxsize=15truecm \epsfbox{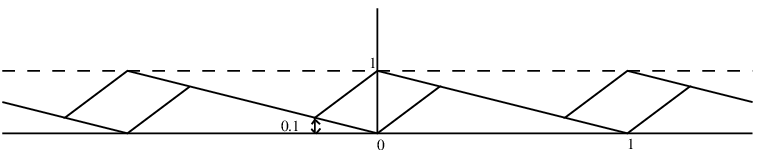}}
\centerline{The kite-like function $\al_0$}

\br
\noindent and the sequence $p_k$ is defined by $p_0=1$ and
$p_{n+1}=10 k_np_n$ for a sequence of integers $k_n\nearrow\infty$
such that
$$
\sum_{n=1}^\infty \frac{p_n}{p_{n+1}}=
\sum_{n=1}^\infty \frac{1}{10k_n}<\infty.
$$
In the modified construction the kite-like parts of $\al_n$
will not be changed but the lines between consecutive kites
will contain larger and larger segments in which the original
straight line will be replaced by graphs of
functions of the form
\begin{equation}\label{f}
f_\tet : t\mapsto \sin (2\pi\tet t),
\end{equation}
properly scaled so that they fit into our strip $\R \times [0,1]$.
At the outset the sequence $k_n$ will be chosen to grow
sufficiently fast in order to leave room for the insertion
of the sine functions.
The parameters $\tet$ will be constructed inductively as a binary tree
of irrational numbers
$\{\tet_\ep: \ep \in \{0,1\}^n\},\ n=1,2,\dots$, where at the $n+1$
stage $\tet_{\ep 0}=\tet_\ep$ and $\tet_{\ep 1}$ is a new point
in $[0,1]$. The numbers $\tet_\ep$ will satisfy
inequalities of the form
\begin{equation}\label{dio}
\|p_n\tet_\ep\| << \frac{1}{n^n},\quad \forall\ \ep\in
\cup_{k=1}^\infty \{0,1\}^k,
\end{equation}
where $\|\la\|$ denotes the distance of the real number $\la$
from the closest integer. The points on the circle which satisfy
the inequality \eqref{dio} at the stage $n+1$, form a union
of finitely many disjoint open intervals and the ``neighbor"
$\tet_{\ep 1}$ of $\tet_{\ep 0}=\tet_\ep$ will be chosen
in that same interval which already contains $\tet_{\ep 0}$.
When the construction is finished we end up with
a Cantor set $\La\subset \T$ consisting of the closure of the set
$\{\tet_\ep: \ep \in \cup_{k=1}^\infty \{0,1\}^k\}$.
At stage $n$ there will be finitely many functions $f_\tet$ with
parameters $\tet_\ep, \ \ep\in\cup_{k=1}^n \{0,1\}^k$
and they will replace segments of the straight lines
connecting the kites of $\al_n$. Each of these
functions will grow in amplitude very gradually from zero
to say $1/100$ and then after running for a long time with
maximal amplitude $1/100$ will symmetrically diminish in
amplitude till it becomes again a straight line. Each
function will appear once and their occurrences will
be separated by very long stretches of the straight line.
Of course this picture will be repeated periodically
between any two consecutive kites of $\al_n$.
Apart from these changes the construction of the functions
$\beta_n$ will be repeated unmodified as in \cite{GW}.

We claim that the construction sketched above, when carefully
carried out, will yield an element $\om\in \Om$ whose orbit
closure $X={\cls} \{T^n\om:n\in \Z\}$ will be, like the original
system, a recurrent-transitive LE system. However, unlike the old
system whose minimal sets were all fixed points, our new system
will have, for each $\tet\in \La$, a minimal subset isomorphic to
the irrational rotation $(R_\tet,\T)$. We will not verify here
these claims, whose proofs parallel the proofs of the original
construction in \cite{GW}. We will though demonstrate the fact
that $(T,X)$ is not HAE. Indeed, this is a direct consequence of
the following proposition. (A second proof will be given in Remark
\ref{r:last}.)

\begin{prop}\label{La}
Let $(T,X)$ be a compact metric cascade and suppose
that there exists an uncountable subset $\La \subset \T$
with the property that for each $\la\in \La$ there exists a subsystem
$Y_\la \subset X$ such that the system $(T,Y_\la)$ is isomorphic
to the rotation $(R_\la,\T)$ on the torus $\T=\R/\Z$.
Then $(T,X)$ is not HAE.
\end{prop}

\begin{proof}
Suppose to the contrary that $(T,X)$ is HAE and let
$Y={\cls}\left(\bigcup\{Y_\la: \la\in \La\}\right)$. By assumption
the system $(T,Y)$ is also HAE and clearly $Y$ coincides with its
mincenter $Y=M(Y)$. Let $A_0$ be a subset of $Y$ such that for
each $\la\in \La$ there is exactly one point in the intersection
$A_0\cap Y_\la$, and let $A= \cup\{ T^n A_0: n\in \Z\}$. If
$\{U_m\}_{m=1}^\infty$ is a countable basis for open sets in $Y$
then the set $O=\cup\{U_m: {\card}(U_m\cap A)\le \aleph_0\}$ is
open and it meets at most countably many $Y_\la$'s. Omitting, at
the outset, this countable set from $\La$ we can and shall assume
that $U_m\cap A$ is uncountable for every $m$. By the AE property
the set $Y_0=Eq(Y)$ of equicontinuity points is a dense $G_\del$
subset of $Y$, and by Proposition \ref{mincenter} each point of
$Y_0$ belongs to a minimal set. Since the set
%nov3 the notation Y_f is occupied
%$Y_f$
 $fix(Y)$ of fixed points in $Y$ is closed, it has an empty
interior and it follows that the set $Y_1=Y_0\setminus fix(Y)$ is
also a dense $G_\del$ subset of $Y$.

Choose a point $z_0\in Y_1$, then $z_0\in Z$ for some nontrivial
minimal set $Z$.
Now the system $Z$ can admit at most a countable set
of eigenvalues and therefore can be {\em not disjoint}
from at most countably many of the systems $Y_\la$. We
can therefore choose an infinite sequence $\{\la_n\} \subset \La$
and a sequence of points $y_n \in Y_{\la_n}$ such that
(i) $\lim_{n\to\infty}y_n=z_0$, (ii)
the set $\{\la_n: n=1,2,\dots\}$ is independent over
the rational numbers $\Q$, and (iii) $Z$ is disjoint from
the minimal system $\prod_{n=1}^\infty (R_{\la_n},\T)$.
Thus the dynamical system
$$
(T,\Om)= (T,Z)\times\prod_{n=1}^\infty (R_{\la_n},\T)
$$
is minimal and in particular for some sequence $m_i$ we have
\begin{gather*}
\lim_{i\to \infty} T^{m_i}y_n = y_n,\ {\text{ for}}\ n=1,2,\dots,\\
{\text{ while}}\quad \lim_{i\to \infty} T^{m_i}z_0 = z_1 \ne z_0.
\end{gather*}
Since $\lim_{n\to\infty}y_n=z_0$,
this contradicts the fact that $z_0$ is an
equicontinuity point and the proof of the lemma is complete.
\end{proof}

This also concludes the proof of Theorem \ref{LEnotHAE}.
\end{proof}

\br

\section{An enveloping semigroup characterization of HNS}
\label{envelop}

\sk

In this section we give an enveloping semigroup characterizations
of Asplund functions and HNS systems in terms of fragmented families
(Definition \ref{d:fr-family}).
In addition to fragmentability,
our approach essentially uses Namioka's theorem.
First we recall this fundamental result and an
auxiliary definition. A topological space $X$ is said to be
{\it \v{C}ech-complete} if $X$ is a $G_{\delta}$ subset in some
%gl-
%its
compact Hausdorff space.
%gl-
If $X$ is either locally compact
Hausdorff or a complete metric space then $X$ is
\v{C}ech-complete. We need the following version of Namioka's
theorem.

\begin{thm} \label{t:N-jct}
\cite{N-jct} (Namioka's Joint Continuity Theorem) Let $w: K \times
X \to M$ be a separately continuous function where $M$ is a metric
space, $K$ is compact and $X$ is \v{C}ech-complete. Then
there exists a dense $G_{\delta}$ set $X_0$ in $X$ such that $w$
is jointly continuous at every point of $K \times X_0$.
\end{thm}

Let $E=E(X)$ be the enveloping semigroup of a compact $G$-system
$X$. Recall that
$$
E^f:=\{p_f: X \to \R \}_{p \in E}, \quad p_f(x)=f(px).
$$
is a pointwise compact subset of $\R^X$, being a continuous image
of $E$ under the map
%aref
 $$q_f: E \to E^f, \hskip 0.2cm q_f(p)=p_f$$
(see Section \ref{Sec-BFT}).

For every $f \in C(X)$ define the map
$$w_f: E \times X \to \R,
\hskip 0.3cm w_f(p,x):=f(px).$$
%23JUne
%It induces two natural mappings:
%$$ E^f \times X \to \R, \hskip 0.2cm (p_f,x) \mapsto f(px)$$
% and
%gl+++
%7July (I see here your sign gl+++ but I did not find what is the change)
In turn $w_f$ induces the mapping $E^f \times X_f \to \R, \hskip
0.2cm (p_f, f_{\sharp}(x)) \mapsto
 f(px).$
Observe that by the proof of Proposition \ref{p:generalX_f}.2
(with $f_{\sharp}=\psi: \beta_G(X)=X \to X_f$) we have
$\psi(x_1)=\psi(x_2) \hskip 0.3cm \text{iff} \ f(gx_1)=f(gx_2)
\hskip 0.3cm \forall g \in G.$
 It follows that
$\psi(x_1)=\psi(x_2) \ \text{iff} \
 f(px_1)=f(px_2) \ \forall p \in E.$
Hence, $E^f \times X_f \to \R$ and the following commutative
diagram are well defined.
%end
\begin{equation*}
\xymatrix{
E \ar@<-2ex>[d]_{q_{f}} \times X
\ar@<2ex>[d]^{f_{\sharp}}\ar[r]  & X \ar[d]^{f} \\
E^f \times X_f \ar[r]  &  \R }
\end{equation*}

We are now ready to prove the following result.

\begin{thm} \label{t:env-asp}
Let $X$ be a compact $G$-system. The following are equivalent:
\ben
\item
$f \in Asp(X)$.
\item
$E^f$ is a fragmented family.
\item
$E^f$ is a barely continuous family.
\item
For every closed ($G$-invariant) subset $Y \subset X$ there exists
a dense $G_{\delta}$ subset $Y_0$ of $Y$ such that  the induced
map $ p_f: Y_0 \to \R, \hskip 0.2cm p_f(y)=f(py)$ is continuous
for every member $p$ of the enveloping semigroup $E$.
 \een
\end{thm}
\begin{proof}
1 $\Rightarrow$ 2: By Theorem \ref{thm-faces} the family
$\breve{G}^f:=\{\breve{g}_f: X \to \R \}_{g \in G}$ is fragmented.
Then the family $E^f$, being the pointwise closure of
$\breve{G}^f$, is also fragmented (Lemma \ref{l:fr-f-cls}).

2 $\Leftrightarrow$ 3: See Definition \ref{d:fr-family}.2.

2 $\Rightarrow$ 4: Since $E^f$ is a fragmented family then for
every closed nonempty subset $Y \subset X$ the family of
restrictions $E^f_Y:=\{p_f \rest_Y: Y \to \R \}$ is (locally)
fragmented. Now by Proposition \ref{fr-sep-1} (see also Definition
\ref{d:fr-family}.1) there exists a dense $G_{\delta}$ subset $Y_0
\subset Y$ such that every $y_0 \in Y_0$ is a point of
equicontinuity of the family $E^f_Y$. Clearly this implies that
$ p_f: Y_0 \to \R$ is continuous for every $p \in E$.

4 $\Rightarrow$ 1: We have to show by Theorem \ref{thm-faces} that
the $G$-map $f_{\sharp}: X \to RUC(G)$ is norm-fragmented.
%oct24
The action of $G$ on $RUC(G)$ preserves the norm.
Therefore, in this case $\mu_G=\mu$ holds, where $\mu$ is the
uniform structure generated by the norm.
By Lemma \ref{l-=}.1 it suffices to check that
$f_{\sharp}\rest_Y: Y \to (RUC(G), \mu)$
%%%
is locally fragmented for every
closed nonempty $G$-subset $Y$ in $X$.

By our assumption we can pick a dense $G_{\delta}$ subset $Y_0$ of
$Y$ such that  the induced map $ p_f: Y_0 \to \R, \quad
p_f(y)=f(py) $ is continuous for every $p \in E(X)$. It follows
that
$$
w_f\rest_{E\times Y_0}: E \times Y_0 \to \R, \quad w_f(p,y)=f(py)
$$
is separately continuous. Since $Y_0$ is \v{C}ech-complete, by
Namioka's theorem there exists a dense subset $Y_1$ of $Y_0$ such
that $w_f\rest_{E\times Y_0}$ is jointly continuous at every
$(p,y_1) \in E \times Y_1$.
%23June I tried to give here more details
Our aim is to prove that $f_{\sharp}\rest_Y: Y \to RUC(G)$ is
continuous at every $y_1 \in Y_1$. In fact we have to show that
every $y_1 \in Y_1$ is a point of equicontinuity of the family of
maps $\{{}_gf\rest_Y : Y \to \R \}_{g \in G}$. By the compactness
of $E$ and the inclusion $\breve{G} \subset E$
%end
it is sufficient to check that the map
$$
w_f\rest_{E\times Y}: E \times Y \to \R
$$
is continuous at each $(p,y_1) \in E \times Y_1$. In order to
check the latter condition fix $\eps >0$. By the joint continuity
of $w_f\rest_{E\times Y_1}: E \times Y_1 \to \R$, one can choose
an open neighborhood $U$ of $p$ in $E$ and an open neighborhood
$O$ of $y_1$ in the space $Y$ such that
$$
|f(py_1)-f(qy)|<\frac{\eps}{3}
$$
holds for every $q\in U$ and $y \in O \cap Y_1$. We claim that
$
 |f(py_1)-f(qz)| < \eps
$ for every $(q,z) \in U \times O$. Fix such a pair $(q,z)$ and
choose $g:=g_{q,z} \in G$ such that the corresponding
$g$-translation $\breve{g}: X \to X$ belongs to $U$ and satisfies
$$
|f(gz)-f(qz)| <\frac{\eps}{3}.
$$
Since $Y_1$ is dense in $Y$ and $\breve{g}: X \to X$ is
continuous, one can pick $a \in Y_1 \cap O$ such that
$$
|f(ga)-f(gz)|<\frac{\eps}{3}.
$$
Putting these estimations together we obtain the desired
inequality $ |f(py_1)-f(qz)| < \eps. $ Thus, we have shown that $
f_{\sharp}\rest_Y: Y \to RUC(G)$ is continuous at every $y_1 \in
Y_1$. Since $Y_1$ is dense in $Y$, we can conclude by Lemma
\ref{simple-fr}.2 that $f_{\sharp}\rest_Y$ is locally fragmented.
\end{proof}

As a corollary we obtain the following enveloping semigroup
characterization of metric RN systems. It certainly can be derived
also from Theorem \ref{hae=rn=v} and the result of
Akin-Auslander-Berg mentioned above (see Theorem \ref{aees}).

\begin{cor} \label{t:env-RN-m}
Let $X$ be a compact metric $G$-system. The following are
equivalent:
 \ben
\item
$(G,X)$ is RN.
%$X$ is RN (equivalently, HAE) $G$-system.
\item
For every closed ($G$-invariant) subspace $Y \subset X$ there
exists a dense
$G_{\delta}$ subset $Y_0$ of $Y$ such that for
 every $p \in E$ the induced map
$ p: Y_0 \to X, \hskip 0.2cm p(y):=py $ is continuous.
 \een
 \end{cor}
\begin{proof} $\Rightarrow$ (1) follows by Theorem \ref{t:env-asp}. Now we prove
that (1) $\Rightarrow$ (2). Since $X$ is a metric compact space we
can choose a countable dense subset $\{f_n: n \in \N \}$ in
$C(X)$. By Theorem \ref{thm1}.4, $C(X)=Asp(X)$. By Theorem
\ref{t:env-asp} for a given closed ($G$-invariant) subset $Y
\subset X$ and every $n \in \N$ there exists a dense $G_{\delta}$
subset $Y_n$ of $Y$ such that for $p \in E$ the induced map $
p_{f_n}: Y_n \to \R $ is continuous. Then it is easy to see that
$Y_0:= \cap_{n \in \N} Y_n$ is the desired subset of $Y$.
\end{proof}

\begin{defn} \label{d:F-semigroup}
 We say that a compact right topological semigroup $S$
 is an {\it $\mathcal F$-semigroup}
 if the family of maps
$\{\lambda_p: S \to S \}_{p \in S},$ where $\lambda_p(s)=ps$, is a
fragmented family.
%nov3
By Definitions \ref{d:fr-family}.1 and \ref{def:fr}.1 it is
equivalent to say that $S^f:=\{p_f: S \to \R \}_{p \in S} \
(\text{where} \hskip 0.3cm p_f(x)=f(px))$ is a fragmented family
for every $f \in C(S)$. Yet another way to formulate the
definition is to require that for every nonempty closed subset $A
\subset S$, every $f \in C(S)$ and $\eps >0$ there exists an open
subset $O \subset S$ such that $A \cap O$ is nonempty and the
subset $f(p(A \cap O))$ is $\eps$-small in $\R$ for every $p \in
S$. %%
\end{defn}

Every compact semitopological semigroup is an $\mathcal
F$-semigroup. The verification is easy applying Namioka's theorem
to the map $S \times A \to \R, \ (s,a) \mapsto f(sa),$ where $A$
is a closed non-empty subset of $S$.

%3July

\begin{thm} \label{t:equifragmented}
Let $X$ be a compact $G$-system. Consider the following conditions:
\begin{itemize}
\item [(a)]
$X$ is HNS (equivalently, $\rm{RN_{app}}$).
\item [(b)]
$\breve{G}:=\{\breve{g}: X \to X \}_{g \in G}$ is a fragmented family.
\item [(c)]
$E(X)=\{p: X \to X\}_{p \in E(X)}$ is a fragmented family.
\item [(d)]
$(G,E(X))$ is HNS (equivalently, $\rm{RN_{app}}$).
\item [(e)]
$E(X)$ is an $\mathcal F$-semigroup.
\end{itemize}
Then we have:
\ben
\item
Always, $(a) \Leftrightarrow (b) \Leftrightarrow (c) \Rightarrow
(d) \Leftrightarrow (e)$.
\item
If $X$ is point transitive then $(a) \Leftrightarrow (b)
\Leftrightarrow (c) \Leftrightarrow (d) \Leftrightarrow (e)$.
\een
\end{thm}
\begin{proof} 1.
%7July  (So is simpler)
%$(a) \Leftrightarrow (c)$:
$(a) \Leftrightarrow (b)$:
%gl+++
%$X$ is $\rm{RN_{app}}$ iff $C(X)=Asp(X)$  (Theorem \ref{thm1}.4).
%On the other hand, $C(X)=Asp(X)$ iff $E^f$ is a fragmented family
%for every $f \in C(X)$ (Theorem \ref{t:env-asp}). Using Lemma
%\ref{simple-fr}.5 it is easy to see that $E^f$ is fragmented for
%every $f \in C(X)$ iff $E=\{p: X \to X \}$ is a fragmented family
The proof follows from Theorem \ref{RN-app=HNS}.
%7July
%and Lemma \ref{l:fr-f-cls}.

$(b) \Leftrightarrow (c)$: Use Lemma \ref{l:fr-f-cls}.

$(a) \Rightarrow (d)$: By the definition $(G,E)$ is a
$G$-subsystem of $X^X$.  Since $\rm{RN_{app}}$ is closed under
subdirect products we get that $E$ is also in $\rm{RN_{app}}$.

$(d) \Leftrightarrow (e)$: $E(X)$ is an $\mathcal F$-semigroup iff
$\{\lambda_p: E \to E \}_{p \in E}$ is a fragmented family iff the
subfamily $\{\lambda_g: E \to E \}_{g \in G}$ is a fragmented
family (use once again Lemma \ref{l:fr-f-cls}). The latter
condition is equivalent to the assertion (d) as it follows by the
equivalence $(a) \Leftrightarrow (b)$ (applied to the system
$(G,E)$).

%By the equivalence $(a) \Leftrightarrow (b)$ (for the system
%$(G,E)$) we obtain that $\{\lambda_g: E \to E \}_{g \in G}$ is a
%fragmented family. Then its pointwise closure $\{\lambda_p: E \to
%E \}_{p \in E}$ is also a fragmented family (Lemma
%\ref{l:fr-f-cls})).

2. $(d) \Rightarrow (a)$: If $x_0$ is a transitive point of $X$
then the map $E \to X, \ p \mapsto px_0$ is a continuous onto
$G$-map. Since $\rm{RN_{app}}$ is closed under quotients we get
%gl+++
that $X$ also belongs to $\rm{RN_{app}}$.
\end{proof}

\begin{cor}
$G^{Asp}$ is an $\mathcal F$-semigroup for every topological group $G$.
\end{cor}
\begin{proof}
The compact $G$-system $X:=G^{Asp}$ is $\rm{RN_{app}}$ by Theorem
\ref{thm1}.6. Therefore, Theorem \ref{t:equifragmented} implies
that the enveloping semigroup $E(G^{Asp})$ is an $\mathcal
F$-semigroup. Since $(G^{Asp}, u_A(e))$ is point-universal
(Proposition \ref{r:ASP}), by Proposition \ref{univ} there exists
a $G$-isomorphism $\phi: (E(G^{Asp}), i(e)) \to (G^{Asp},u_A(e))$
of pointed $G$-systems. In fact this map is an isomorphism of
(right topological) semigroups because $u_A(G)$ is dense in
$G^{Asp}$ and $i(G)$ is dense in $E(G^{Asp})$.
\end{proof}

\begin{cor} \label{hae>E1}
Let $(G,X)$ be a compact HNS system. Then $p: X \to X$ is fragmented
(equivalently, Baire class 1, when $X$ is metric) for every $p \in
E(X)$.
\end{cor}
\begin{proof}
Use Theorem \ref{t:equifragmented} (and Proposition
\ref{Baire-1}.2).
\end{proof}

%end

For the definition of Rosenthal compacts see Section
\ref{Sec-BFT}.

\begin{thm} \label{t:B1}
\ben
\item
Let $X$ be a compact metrizable $G$-space. For every $f \in
Asp(X)$ the compact space $E^f \subset \R^X$ is a Rosenthal
compact.
\item
Let $X$ be a metrizable compact RN $G$-space. Then
the enveloping semigroup $E$ is a (separable) Rosenthal compact
with cardinality $\leq 2^{\aleph_0}$ (in particular, no
subspace of
%gl+
$E$ can be homeomorphic to $\beta \N$).
 \een
\end{thm}
\begin{proof} 1. Since $f \in Asp(X)$, by Theorem \ref{t:env-asp},
$E^f=\{p_f: X \to \R \}_{p \in E}$ is a fragmented family. In
particular, each map $p_f: X \to \R$ is fragmented. Since $X$ is
compact and metrizable we can apply Proposition \ref{Baire-1}.
Hence, each function $p_f \in E^f$ is of Baire class 1 (on the
Polish space $X$). Therefore, $E^f$ is a Rosenthal compact.

2. $C(X)=Asp(X)$ by Theorem \ref{thm1}.4. It follows by the first
assertion that $E^f$ is a Rosenthal compact for every $f \in
C(X)$. An application of the dynamical version of the BFT theorem,
Theorem \ref{D-BFT}, concludes the proof.
\end{proof}

\begin{remark} \label{r:last}
%6July
%\ben
%\item
%In fact the following stronger version of Theorem \ref{t:B1}.1 is
%true. If $\breve{G}$ is uniformly Lindel\"of (in particular, if
%$X$ is metrizable) then for every compact, not necessarily metric,
%system $X$ and every $f \in Asp(X)$ the space $E^f$ is a Rosenthal
%compact. Indeed, using Proposition \ref{p:generalX_f} it is easy
%to see that $E^f$ is naturally homeomorphic to
%%23June
%$E^{F_e}$ where for the maps $\psi=f_{\sharp}: X \to X_f$,
%$F_e: X_f \to \R$ we have $F_e \circ f_{\sharp}=f$.
%%end
%%3July
%We can apply Theorem \ref{t:B1} because $X_f$ is metrizable by
%Proposition \ref{prop:blocks} and $F_e \in Asp(X_f)$ by Remark
%\ref{r:asp-def}.2.
%\item
%end
Theorem \ref{t:B1}.2 can be used to obtain an alternative proof of
Proposition \ref{La}. In fact, as can be seen
%3July HOW ??
from Proposition \ref{envel},
%end
the enveloping semigroup of the system $(T,X)$ in Proposition
\ref{La} has cardinality $2^{2^{\aleph_0}}$.
%\een
\end{remark}
%%oct21
%\item
%Weakly almost periodic functions also admit an enveloping
%semigroup characterization. Namely, using Grothendieck's theorem
%(see the note before Theorem \ref{semi}), one can show that $f \in
%WAP(X)$ iff the mapping $E^f \times X_f \to \R, \hskip 0.2cm (p_f,
%f_{\sharp}(x)) \mapsto f(px)$ is separately continuous (iff the
%function induced on $X_f$ by $f \circ p$ is continuous for every
%$p \in E$). In particular, this ``local characterization" covers
%Ellis' result (see Theorem \ref{semi}). This characterization
%provides also a new (independent) proof of the inclusion $WAP(X)
%\subset Asp(X)$ (use assertion 4 of Theorem \ref{t:env-asp}).

%end

Our next example is of a metric minimal cascade $(T,X)$
which is not RN yet its enveloping semigroup
$E=E(T,X)$ (a) is a separable Rosenthal compact of cardinality
$2^{\aleph_0}$ and (b) has the property that each $p \in E$
is of Baire class 1.
Thus this example shows that the converse of Theorem
\ref{t:B1}.2 does not hold and neither does that of
Corollary \ref{hae>E1}.

\begin{exa}\label{two-arr}
Let $\T=\R/\Z$ be the one-dimensional torus, and let $\al\in \R$
be a fixed irrational number and $T_\al: \T \to \T$ is the
rotation by $\al$,\ $T_\al \beta = \beta  +\al \pmod{1}$. We
define a topological space $X$ and a continuous map $\pi: X \to
\T$ as follows. For $\beta\in \T \setminus \{n\al: n \in \Z\}$ the
preimage $\pi^{-1}(\beta)$ will be a singleton $x_\beta$. On the
other hand for each $n\in \Z$, $\pi^{-1}(n\al)$ will consist of
exactly two points $x^{-}_{n\al}$ and $x^{+}_{n\al}$. For
convenience we will use the notation $\beta^{\pm}$,\ ($\beta \in
\T$) for points of $X$, where $(n\al)^{-}=x^{-}_{n\al}$,\
$(n\al)^{+}=x^{+}_{n\al}$ and $\beta^{-}=\beta^{+}=x_\beta$ for
$\beta\in \T \setminus \{n\al: n \in \Z\}$. A basis for the
topology at a point of the form $x_\beta,\ \beta\in \T \setminus
\{n\al: n \in \Z\}$, is the collection of sets
$\pi^{-1}(\beta-\ep,\beta +\ep),\ \ep > 0$. For $(n\al)^{-}$ a
basis will be the collection of sets of the form $\{(n\al)^{-}\}
\cup \pi^{-1}(n\al-\ep,n\al)$, where $\ep>0$. Finally for
$(n\al)^{+}$ a basis will be the collection of sets of the form
$\{(n\al)^{+}\} \cup \pi^{-1}(n\al,n\al+\ep)$. It is not hard to
check that this defines a compact metrizable zero dimensional
topology on $X$ (in fact $X$ is homeomorphic to the Cantor set)
with respect to which $\pi$ is continuous. Next define $T: X\to X$
by the formula $T\beta^{\pm} = (\beta + \al)^{\pm}$. Again it is
not hard to see that $\pi: (T,X) \to (R_\al,\T)$ is a homomorphism
of dynamical systems and that $(T,X)$ is minimal and not
equicontinuous (in fact it is
%nov3 not defined in the paper (so maybe some reference is good here)
almost-automorphic; see e.g. Veech \cite{V65}).
In particular $(T,X)$ is not RN.

We now define for each $\ga \in \T$ two distinct maps
$p_\ga^{\pm}: X \to X$ by the formulas
%gl++
$$
p_\ga^+(\beta^{\pm})= (\beta + \ga)^+,
\qquad
p_\ga^-(\beta^{\pm})= (\beta + \ga)^-.
$$
We leave the verification of the following claims
as an exercise.
\begin{enumerate}
\item
For every $\ga\in \T$ and every {\em sequence}, $n_i
\nearrow \infty$ with $\lim_{i\to\infty}n_i\al = \ga,\
{\text{and}}\ \forall i,\ n_i\al < \ga$, we have
$\lim_{i\to\infty}T^{n_i} = p_\ga^{-}$ in $E(T,X)$.
An analogous statement holds for $p_\ga^{+}$.
\item
$E(T,X) = \{T^n: n\in \Z\} \cup
\{p_\ga^{\pm}: \ga \in \T\}$
\item
The subspace $\{T^n: n\in \Z\}$ inherits from
$E$ the discrete topology.
\item
The subspace $E(T,X) \setminus \{T^n: n\in \Z\} = \{p_\ga^{\pm}:
\ga \in \T\}$ is homeomorphic to the ``two arrows" space of
Alexandroff and Urysohn (see \cite[page 212]{Eng}, and also Ellis'
example
%me ``two circles" has also another meaning
%``two circles" example
%me 5.29 and not 5.59
\cite[Example 5.29]{E}).
It thus follows that $E$ is a separable
Rosenthal compact of cardinality $2^{\aleph_0}$.
\item
For each $\ga\in \T$ the complement of the set
$C(p_\ga^{\pm})$ of continuity points of $p_\ga^{\pm}$ is
the countable set $\{\beta^{\pm}: \beta + \ga = n\al,
\ {\text{for some}}\ n\in \Z\}$.
In particular each element of $E$ is of Baire class 1.
\end{enumerate}
\end{exa}

%gl+

\br

\section{A dynamical version of
Todor\u{c}evi\'{c}' theorem}\label{Tod}

\sk

A surprising result of Todor\u{c}evi\'{c} asserts that a Rosenthal
compact $X$ which is not metrizable obeys the following
alternative: either $X$ contains an uncountable discrete
subspace or it is an at most two-to-one continuous
preimage of a compact metric space
(\cite[Theorem 3]{T}).
%23June
%An application of a
%standard trick then yields the following result.
We present here the following dynamical version.

\begin{prop}\label{Todor}
(A dynamical Todor\u{c}evi\'{c} dichotomy)
%23June
%$\Ga$ be a countable discrete group.
Let $G$ be a uniformly Lindel\"of group and $(G,X)$ a compact
%point transitive
system with the property that $X$ is a Rosenthal compact. Then
either $X$ contains an uncountable discrete subspace or there
exists a metric dynamical system $(G,Y)$ and a $G$-factor $\pi:
(G,X) \to (G,Y)$ such that $|\pi^{-1}(y)| \le 2$ for every $y\in
Y$.
\end{prop}

\begin{proof}
%If we rule out the first alternative in Todor\u{c}evi\'{c}'s
%theorem then it follows by that theorem that there exists a
%compact metric space $Z$ and a continuous map $\phi: X \to Z$ with
%the property that $|\phi^{-1}(z)|\le 2$ for every $z\in Z$. Define
%a map $\pi: X \to Z^\Ga$ by $\pi(x)_\ga=\phi(\ga x),\ x\in X,\
%\ga\in \Ga$ and set $Y=\pi(X)$. It is easy to check that $\pi:
%(\Ga,X) \to (\Ga,Y)$ is a homomorphism of dynamical systems, where
%$\Ga$ acts on $Z^\Ga$ by permuting the coordinates
%$(\ga(y))_\del=y_{\del\ga},\ y\in Z^\Ga,\ \ga,\del\in \Ga$.
%Clearly $Y$ is a metrizable compact space and finally
%$\pi(x)=\pi(x')$ implies in particular that
%$\phi(x)=\pi(x)_e=\pi(x')_e=\phi(x')$, hence also $|\pi^{-1}(y)|
%\le 2$ for every $y\in Y$.

If we rule out the first alternative in Todor\u{c}evi\'{c}'s
theorem then it follows by that theorem that there exists a
compact metric space $Z$ and a continuous map $\phi: X \to Z$ with
the property that $|\phi^{-1}(z)|\le 2$ for every $z\in Z$. By
\cite[Theorem 2.11]{Me-sing} there exist a compact metric
$G$-space $Y$, a continuous onto $G$-map $f_1: X \to Y$ and a
continuous map $f_2: Y \to Z$ such that $\phi=f_2 \circ f_1$.
Clearly, $|f_1^{-1}(y)|\le 2$ for every $y \in Y$.
%end
\end{proof}

We do not know whether Theorem \ref{t:B1}.2
can be strengthened to the statement that
the enveloping semigroup of any compact metric RN system is
in fact metric. However, Proposition \ref{Todor}
yields the following.

\begin{cor}\label{Todor-cor}
Let $X$ be a metric RN $G$-system,
%23June
%where $G$ is a countable discrete group,
where $G$ is an arbitrary topological group. Then either the
enveloping semigroup
%3July
$E=E(X)$ contains an uncountable discrete subspace or it admits a
metric $G$-factor
%gl+
$\pi: (G,E) \to (G,Y)$ such that $|\pi^{-1}(y)| \le 2$ for every
$y\in Y$.
%end
\end{cor}
\begin{proof}
This follows directly from Theorem \ref{t:B1}.2 and Proposition
\ref{Todor} because
%3July
%we can assume that the acting group is second countable
the natural restriction ${\breve G}$ (see Section \ref{Sec-BFT})
is second countable
%end
(and hence, uniformly Lindel\"of).
\end{proof}
%end

%3July
%
%\sk
%
%\begin{remarks}\label{second}
%\ben
%\item
%Theorem \ref{t:B1}.2 can be used to obtain an alternative proof of
%Proposition \ref{La}. In fact, as can be seen
%%3July HOW ??
%from Proposition \ref{envel},
%%end
%the enveloping semigroup of the system $(T,X)$ in Proposition
%\ref{La} has cardinality $2^{2^{\aleph_0}}$.
%%oct21
%\item
%Weakly almost periodic functions also admit an enveloping
%semigroup characterization. Namely, using Grothendieck's theorem
%(see the note before Theorem \ref{semi}), one can show that $f \in
%WAP(X)$ iff the mapping $E^f \times X_f \to \R, \hskip 0.2cm (p_f,
%f_{\sharp}(x)) \mapsto f(px)$ is separately continuous (iff the
%function induced on $X_f$ by $f \circ p$ is continuous for every
%$p \in E$). In particular, this ``local characterization" covers
%Ellis' result (see Theorem \ref{semi}). This characterization
%provides also a new (independent) proof of the inclusion $WAP(X)
%\subset Asp(X)$ (use assertion 4 of Theorem \ref{t:env-asp}). \een
%\end{remarks}
%end

\begin{problem}
By Theorem \ref{t:B1}.2 the enveloping semigroup of the G-W
example is a separable Rosenthal compact (of cardinality
$2^{\aleph_0}$). We do not have a concrete description of this
enveloping semigroup and do not even know whether it is metrizable
%3July
%is it true that at least E, as the natural cascade, is also RN
%(we know of course that it is RN-appr) ?
%GENERAL Q: Let $X$ be an RN METRIC system. Is it true that $(G,E)$ is
%also RN (what if G=Z) ?
%CLOSELY RELATED QUESTION: Is it true that there exists a dense subset
%$D$ of $X$
%such that every $p \in E$ is continuous on $D$.
%Can we guarantee that $D$ "separates" $E$. That is whether
%for distinct $p,q \in E$ one can find $x \in D$ such that $px$
%and $qx$ are different ?
%end
or if it contains an uncountable discrete subspace.
\end{problem}

\br

\bibliographystyle{amsplain}

\end{document}